\documentclass[runningheads]{llncs}
\usepackage[T1]{fontenc}
\usepackage{graphicx}
\usepackage[dvipsnames]{xcolor}
\usepackage{graphicx}
\usepackage{amssymb}
\usepackage{amsfonts}
\usepackage{mathtools}
\usepackage{mathrsfs}
\usepackage{mathabx}
\usepackage{color}
\usepackage{transparent}
\usepackage{scalerel}
\usepackage{tikz}
\usepackage{trimclip}
\usepackage{relsize}
\usepackage{yfonts}
\newlength{\hrawidth}
\newlength{\hraheight}
\newlength{\minuswidth}
\newcommand{\threelines}{%
    \mathrel{%
        \clipbox{0 0 .55\hrawidth{} 0}{$\Rrightarrow$}%
    }
}
\newcommand{\LRarrow}{\Lleftarrow\mkern-8mu\Rrightarrow}
\newcommand{\Rlongarrow}{\threelines\mkern-4mu\Rrightarrow}

\renewcommand{\l}{\mathtt{l}}

\newcommand{\s}{\mathtt{s}}

\renewcommand{\L}{\mathtt{L}}
\newcommand{\M}{\mathtt{M}}
\newcommand{\m}{\mathtt{m}}
\renewcommand{\P}{\mathtt{P}}

\newcommand{\Q}{\mathtt{Q}}
\newcommand{\I}{\mathtt{I}}
\newcommand{\J}{\mathtt{J}}
\newcommand{\K}{\mathtt{K}}
\newcommand{\V}{\mathtt{V}}
\newcommand{\U}{\mathtt{U}}
\newcommand{\W}{\mathtt{W}}
\newcommand{\A}{\mathtt{A}}
\newcommand{\B}{\mathtt{B}}
\newcommand{\C}{\mathtt{C}}
\newcommand{\D}{\mathtt{D}}
\newcommand{\E}{\mathtt{E}}
\newcommand{\G}{\mathtt{G}}
\renewcommand{\H}{\mathtt{H}}
\newcommand{\F}{\mathtt{F}}

\newcommand{\X}{\mathtt{X}}
\newcommand{\Y}{\mathtt{Y}}
\newcommand{\Z}{\mathtt{Z}}
\newcommand{\nnl}{\vspace{0.15cm}\\}
\newcommand{\1}{\mathtt{1}}
\newcommand{\2}{\mathtt{2}}
\newcommand{\3}{\mathtt{3}}
\newcommand{\4}{\mathtt{4}}
\newcommand{\5}{\mathtt{5}}
\newcommand{\6}{\mathtt{6}}
\newcommand{\7}{\mathtt{7}}
\newcommand{\8}{\mathtt{8}}
\newcommand{\9}{\mathtt{9}}

\usetikzlibrary{shapes.geometric}
\usetikzlibrary{decorations.pathreplacing,calligraphy}
\usetikzlibrary{decorations.pathreplacing,tikzmark}
\usetikzlibrary{decorations.markings}
\usetikzlibrary{decorations.text,positioning,fit,quotes} 
\usetikzlibrary{arrows.meta}
\usetikzlibrary{calc,intersections,through,backgrounds,patterns,patterns.meta}

\begin{document}
\title{Manifold-based Proving Methods in Projective Geometry}
%
%
\author{Michael Martin Katzenberger \and
Jürgen Richter-Gebert}

\author{
{Michael Martin Katzenberger\textsuperscript{1}, Jürgen Richter-Gebert\textsuperscript{2}}\\
\vspace{2pt}
\textsuperscript{1} {\small\tt michael.martin.katzenberger@tum.de}\small{, Orcid: 0009-0004-2270-2502}\\
\textsuperscript{2} {\small\tt richter@tum.de}\small{, Orcid: 0009-0000-0049-9743}\\
}

\authorrunning{M. Katzenberger and J. Richter-Gebert}
%
\institute{Department of Mathematics, CIT, Technical University of Munich}
\maketitle              
\begin{abstract}
This article compares different proving methods for projective incidence theorems. In particular, a technique using quadrilateral tilings recently introduced  by Sergey Fomin and Pavlo Pylyavskyy is shown to be at most as strong as proofs using bi-quadratic final polynomials and thus, also proofs using Ceva-Menelaus-tilings. Furthermore, we demonstrate the equivalence between quadrilateral-tiling-proofs and proofs using exclusively Menelaus configurations. We exemplify the transition between the proofs in several examples in 2D and in 3D.

\end{abstract}
\section{Introduction}
Projective geometry studies properties that are invariant under projective transformations. We here focus on the  most basic of such properties: Incidence relations between points and lines (and in 3D also planes). They are natural projective invariants since they are stable under projective transformations. The typical blueprint of a projective incidence theorem reads like:
\[
  \textit{\small some incidences hold}
 +
  \textit{\small  some incidences do not hold} 
\implies
  \textit{\small another incidence holds.}
\]

Incidence theorems only including points and lines  are an important basic form of planar geometric theorems. By suitable tools many other geometric statements can ultimately be reformulated as theorems of this form and thus be reduced to statements involving points/lines/incidences only.

There is an extensive literature on the automatic generation of proofs in such a context, reaching from very general algebraic tools like Gröbner Bases~\cite{KS}, and Wu's Method~\cite{Wu} via still general more invariant geometric approaches like {\it final polynomials}~\cite{BS} to more specialised methods that exploit specific cancelation patterns in terms of bracket algebra and unveil some of the geometric reasons ``why'' some theorems do  hold. We here are interested in the comparison of the latter class of proving methods.

In this article we will compare (bi-quadratic) final polynomials~\cite{BS,RG2},
proofs by manifold structures on  Ceva-/Menelaus-configurations~\cite{ARG,RG1} and the recently introduced quadrilateral tiling  method by
 Fomin and  Pylyavskyy~\cite{FP,PS}. All three methods rely on cancellation patterns on the bracket level. The last two methods additionally have in common that the cancellation patterns can be organised on an orientable 2-manifold such that each cancellation corresponds to two patches of the manifold meeting along an edge. The 
statement that a theorem holds corresponds to the fact that the underlying manifold does not have a boundary. In other words: One has a theorem if everything closes up nicely.

We will quickly recall the established proving methods, present quadrilateral-tiling-proofs and then show that they are at most as powerful as the other techniques.

\section{Binomial Proofs and Ceva-Menelaus-Proofs}
\subsection{Binomial Proofs}
As usual we will abbreviate the determinant of three planar points in homogeneous coordinates by $[A,B,C]:=\det{(p_A,p_B,p_C)}$.
Binomial proofs have fist been introduced in the context of oriented matroid theory~\cite{BRG,RG2} and are by now pretty well understood. They utilize bi-quadratic final polynomials coming from Grassmann-Plücker-relations. For any points $\A,\B,\C,\D,\E$ in the projective plane the following equation holds:
\begin{equation}
        [\A,\B,\C][\A,\D,\E] - [\A,\B,\D][\A,\C,\E] + [\A,\B,\E][\A,\C,\D] = 0
\end{equation}
This means, that collinearity of $\A,\B$ and $\C$ immediately implies
\begin{equation}
        [\A,\B,\D][\A,\C,\E] = [\A,\B,\E][\A,\C,\D].
        \label{eq: grass2}
\end{equation}
And on the other hand, if equation (\ref{eq: grass2}) holds and we know that $\A,\D,\E$ are not collinear, collinearity of $\A,\B,\C$ is implied. 
In the binomial proving technique incidences of the hypotheses as well as the conclusion
are translated in equations of type (\ref{eq: grass2}). A proof consists of a collection of such equations such that the equation that represents the conclusion is obtained from hypotheses equations by multiplying all terms on the left and all all terms on the right and cancelling terms that occur on both sides. The non-incidence statements are used as non-degeneracy conditions to make sure that no division by zero arises along the cancellation process. While this method is not known to work in general it is known for producing short, readable and insightful proofs if it works.
Actually to the best of our knowledge it still might be
the case that if one allows to add auxiliary points every incidence theorem that holds over arbitrary fields admits a binomial proof. An equivalent question was posed in \cite{PS} as Problem 3.3. No proving strategy or counterexample currently seems to be in reach.

\subsection{Ceva-Menelaus-Proofs}

This proving method introduced in~\cite{RG1} uses multiple instances of Ceva's and Menelaus's theorems to prove incidences. The structure of the proof is organised along a triangulated orientable manifold. Each triangle carries a Ceva or a Menelaus configuration. Edge points of adjacent triangles are identified. If the corresponding incidence structure of the configuration holds for all but one of the triangles then it will automatically hold for the last triangle, hence a certain incidence theorem holds. Here the length ratios along all edges 
\begin{figure}[h]
    \begin{minipage}[c]{0.5\textwidth}
        \centering
        \vspace{0.5cm}
        \begin{tikzpicture}[]
            \coordinate [label=below:{$\A$}] (A) at (0,0);
            \coordinate [label=below:{$\B$}] (B) at (2,0);
            \coordinate [label=below:{}] (XY) at (3,0);
            \coordinate [label=above:{$\C$}] (C) at (2/3,2.5*2/3);
        
            \path [name path=A--B] (A) -- (B);
            \draw[very thick] (A) -- (B);
            \path [name path=A--C] (A) -- (C);
            \draw[very thick] (A) -- (C);
            \path [name path=C--B] (C) -- (B);
            \draw[very thick] (C) -- (B);
        
            \draw [very thick] (B) -- ($ (A)!2!(B) $) coordinate [label=below:$\Z$] (Z);
            \path [name path=B--X] (B) -- (Z);
            \draw [] (Z) -- ($ (A)!0.5!(C) $) coordinate [label=left:$\Y$] (Y);
            \path [name path=Z--Y] (Z) -- (Y);
        
            \path [name intersections={of=C--B and Z--Y,by={[label=above right:$\X$]X}}];
            
            \foreach \point in {A,B,C} {
                \fill [black,opacity=1] (\point) circle (2pt);
            }
            \foreach \point in {X,Y,Z} {
                \path[fill=white,draw=black, thick] (\point) circle[radius=1.5pt];
            }
        \end{tikzpicture}
        \vspace{0.22cm}
        \begin{equation*}
            \frac{\overrightarrow{|\A\Z|}}{\overrightarrow{|\Z\B|}}\cdot\frac{\overrightarrow{|\B\X|}}{\overrightarrow{|\X\C|}}\cdot\frac{\overrightarrow{|\C\Y|}}{\overrightarrow{|\Y\A|}}=-1
            \label{menelaus equation}
        \end{equation*}
    \end{minipage}
    \hfill
    \begin{minipage}[c]{0.45\textwidth}
        \centering
        \begin{tikzpicture}[]
            \coordinate [label=left:{$\A$}] (A) at (0,0);
            \coordinate [label=right:{$\B$}] (B) at (3,0);
            \coordinate [label=above:{$\C$}] (C) at (1,2.5);
            \coordinate [label={[label distance=0.15cm]right:$\D$}] (D) at (1.25,1);
            \path [name path=A--B] (A) -- (B);
            \draw[very thick] (A) -- (B);
            \path [name path=A--C] (A) -- (C);
            \draw[very thick] (A) -- (C);
            \path [name path=C--B] (C) -- (B);
            \draw[very thick] (C) -- (B);
        
            \draw [] (A) -- ($ (A)!1.5!(D) $) coordinate (E);
            \path [name path=A--E] (A) -- (E);
            \draw [] (B) -- ($ (B)!1.4!(D) $) coordinate (F);
            \path [name path=B--F] (B) -- (F);
            \draw [] (C) -- ($ (C)!1.7!(D) $) coordinate (G);
            \path [name path=C--G] (C) -- (G);
            
            \path [name intersections={of=C--B and A--E,by={[label=above right:$\X$]X}}];
            \path [name intersections={of=A--C and B--F,by={[label=above left:$\Y$]Y}}];
            \path [name intersections={of=A--B and C--G,by={[label=below:$\Z$]Z}}];
            \foreach \point in {A,B,C} {
                \fill [black,opacity=1] (\point) circle (2pt);
            }
            \path[fill=white,draw=black, thick] (D) circle[radius=2pt];
            \foreach \point in {X,Y,Z} {
                \path[fill=white,draw=black, thick] (\point) circle[radius=1.5pt];
            }
        \end{tikzpicture}
        \vspace{-0.1cm}
        \begin{equation*}
            \frac{\overrightarrow{|\A\Z|}}{\overrightarrow{|\Z\B|}}\cdot\frac{\overrightarrow{|\B\X|}}{\overrightarrow{|\X\C|}}\cdot\frac{\overrightarrow{|\C\Y|}}{\overrightarrow{|\Y\A|}}=1
        \end{equation*}
    \end{minipage}
    \caption{Menelaus's and Ceva's theorem and their corresponding equations.}\label{fig:enter-label}
\end{figure}
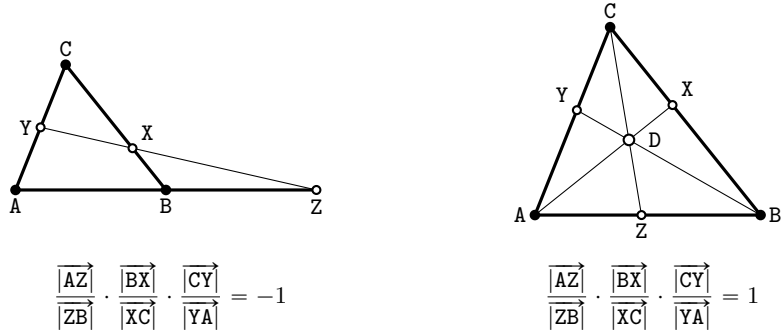
between two assigned triangles cancel  such that only the length ratios along the boundary of the unassigned triangle remain. This implies that the final triangle must also correspond to either a Ceva or Menelaus configuration forcing an additional incidence (the conclusion) to hold. For more details, see~\cite{RG1}. In~\cite{ARG} a procedure to translate a binomial proof into a Ceva-Menelaus-proof ({\it CM-proof}, for short) and vice versa is described, showing that these proving techniques are equally powerful, as long as it is allowed to add two auxiliary generic points to the configuration.

\section{Quad proofs á la Fomin and  Pylyavskyy}
In the  recent paper ``Incidences and Tilings''~\cite{FP},  Fomin and  Pylyavskyy introduce another proving technique that also relies on a manifold based approach. Unlike in the Ceva-Menelaus technique in this approach the basic building blocks
are quadrilateral tiles, each of them encoding a kind of cross ratio. The idea is to have
a variable cancel out whenever two tiles meet along an edge and the conclusion being 
the last tile to be filled. This is very similar to  the Ceva-Menelaus approach.
In a certain sense the quadrilateral building blocks are of more elementary nature: The vertices of each quadrilateral face (we will call them quadrangles) represent two points $\P,\Q$ and two lines
 $\l,\m$ not incident to each other and the additional 
\begin{figure}[ht]
    \centering
    \begin{minipage}{0.3\linewidth}
        \begin{tikzpicture}[>=Stealth]
            \node (A) at (0cm,0cm) {$\P$};
            \node (l) at (2cm,0cm) {$\l$};
            \node (B) at (2cm,-2cm) {$\Q$};
            \node (m) at (0cm,-2cm) {$\m$};
            \draw[-] (A) -- (l);
            \draw[-] (B) -- (l);
            \draw[-] (A) -- (m);
            \draw[-] (B) -- (m);
        \end{tikzpicture}
    \end{minipage}
    \begin{minipage}{0.3\linewidth}
        \begin{tikzpicture}[>=Stealth]
            \coordinate [label=left:{$\P$}] (A) at (0,0);
            \coordinate [label=right:{$\Q$}] (B) at (4,0.5);
            \coordinate [label=above right:{$\l$}] (C) at (0,1);
            \coordinate [label=below:{}] (D) at (4,-0.5);
            \coordinate [label= below right:{$\m$}] (E) at (2.5,1.5);        
        
            \path [name path=A--B] (A) -- (B);
            \draw[very thick] (A) -- (B);
            \path [name path=D--C] (D) -- (C);
            \draw[-] (D) -- (C);

            \path [name intersections={of=A--B and D--C,by={[label=above right:]X}}];

            \draw[-] (E) -- ($ (E)!2!(X) $);

            \foreach \point in {A,B} {
                \path[fill=white,draw=black, thick] (\point) circle[radius=1.5pt];
            }

            \path[fill=white,draw=black, thick] (X) circle[radius=2pt];
        \end{tikzpicture}
    \end{minipage}
    \caption{A coherent quadrangle and the underlying incidence structure.}\label{fig: quadrangle}
\end{figure}
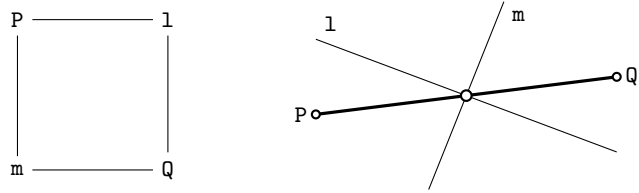
 \textit{``coherency''-condition} that $\P,\Q$ and $\l\land\m$ are collinear.  
  The collinearity of these three points, or --- dually --- the concurrency of the $3$ lines $\l,\m$ and $\P\lor\Q$, is equivalent to the, also self-dual, algebraic expression
\begin{equation}
    \frac{\langle\P,\l\rangle\cdot \langle\Q,\m\rangle}{\langle\P,\m\rangle\cdot\langle\Q,\l\rangle} = 1.
    \label{fomin eq}
\end{equation}
If one traverses such a quadrangle counter-clockwise, the scalar products corresponding to edges from {\it line} to {\it point} appear in the numerator of the fraction and the determinants corresponding to edges from {\it point} to {\it line} appear in the denominator of the fraction. This means gluing two counter-clockwise oriented quadrangles along an edge will result in the scalar products coming from this edge appearing once in the numerator and once in the denominator and thus cancel out. So, similar to Ceva-Menelaus-proofs, if in a tiling of a closed oriented surface with these quadrangles, all tiles but one are coherent, the remaining tile will be coherent, too. Thus, every quadrilateral tiling encodes a proof (we will refer to such proofs as \textit{quad-proofs}) of a projective incidence theorem. In~\cite{FP} many examples of quadrilateral tiling proofs for projective incidence theorems are given. 
One might guess that this more elementary proving technique provides more general proofs and can prove larger classes of theorems compared to the Ceva-Menelaus techniques. However, it turns out that every quad-proof can be turned to an equivalent proof using Menelaus theorems only.

\subsection{Extracting a binomial proof from a quad-proof}\label{section: extract bin}

As the coherency of a quadrangle is equivalent to an equation being satisfied, it seems natural trying to relate this proving method to binomial proofs. And in fact the process of extracting a binomial proof from a quad-proof is quite straight-forward. By choosing two points $\L_1,\L_2$ spanning $\l$, two points $\M_1,\M_2$ spanning $\m$ and defining the incidence point $\I = \l\land \m$ each coherent quadrangle actually encodes $3$ collinearities, yielding the equations
\begin{align*}
    \L_1,\L_2,\I\ \text{col.}&\quad\Longleftrightarrow\quad[\L_1,\L_2,\P][\L_1,\I,\Q]=[\L_1,\L_2,\Q][\L_1,\I,\P],\\
    \M_1,\M_2,\I\ \text{col.}&\quad\Longleftrightarrow\quad[\M_1,\M_2,\Q][\M_1,\I,\P]=[\M_1,\M_2,\P][\M_1,\I,\Q]\text{ and}\\
    \I,\P,\Q\ \text{col.}&\quad\Longleftrightarrow\quad[\I,\P,\L_1][\I,\Q,\M_1]=[\I,\P,\M_1][\I,\Q,\L_1].
\end{align*}
The multiplying left and right sides of these equations yields
\begin{equation}
    [\P,\L_1,\L_2][\Q,\M_1,\M_2]=[\P,\M_1,\M_2][\Q,\L_1,\L_2],
\end{equation}
which is equivalent to equation (\ref{fomin eq}). Applying this procedure to each quadrangle thus immediately gives us a binomial proof, so quad-proofs can only be ``as powerful as'' binomial proofs --- and hence also Ceva-Menelaus-proofs.
In other words: each quad-proof can be translated into a bi-quadratic proof and from there into a CM-proof.

\section{Relations between the Tiling-Based Methods}
\subsection{Quad-proofs and Ceva-Menelaus-proofs}

\begin{figure}[ht]
    \centering
    \hspace{-0.3cm}
    \begin{minipage}[c]{0.37\linewidth}
        \begin{tikzpicture}[>=Stealth]
            \node (A) at (0,0) {$\P_2$};
            \node (C) at (-3,0) {$\P_1$};
            \node (V) at (-1.5,1.73*1.5) {$\P_3$};
            \node (W) at (-1.5,-1.73*1.5) {$\P_4$};

            \node (X2) at ($ (C)!0.5!(W) $) {};
            \node (Y) at ($ (C)!0.5!(V) $) {};

            \node (A2) at ($ (A)!0.5!(C) $) {};
            \node (A5) at ($ (A)!0.5!(V) $) {};
            \node (A6) at ($ (A)!0.5!(W) $) {};

            \node (l1) at (barycentric cs:A=1,C=1,V=1) {$\l_5$};
            \node (l2) at (barycentric cs:A=1,C=1,W=1) {$\l_6$};

            \node (l3) at ($ (l1)!2!(Y) $) {$\l_1$};
            \node (l4) at ($ (l2)!2!(X2) $) {$\l_2$};
            \node (l5) at ($ (l2)!2!(A6) $) {$\l_3$};
            \node (l6) at ($ (l1)!2!(A5) $) {$\l_4$};
            
            \draw[-] (A) -- (l1);
            \draw[-] (A) -- (l2);
            \draw[-] (C) -- (l1);
            \draw[-] (C) -- (l2);
            \draw[-] (V) -- (l1);
            \draw[-] (W) -- (l2);
            \draw[-] (C) -- (l3);
            \draw[-] (C) -- (l4);
            \draw[-] (A) -- (l5);
            \draw[-] (A) -- (l6);
            \draw[-] (V) -- (l3);
            \draw[-] (W) -- (l4);
            \draw[-] (W) -- (l5);
            \draw[-] (V) -- (l6);
            \draw[-] (l3) -- ($ (Y)!1.85!(l3) $);
            \draw[-] (l4) -- ($ (X2)!1.85!(l4) $);
            \draw[-] (l5) -- ($ (A6)!1.85!(l5) $);
            \draw[-] (l6) -- ($ (A5)!1.85!(l6) $);
            \draw[dash pattern=on 1pt off 3pt] ($ (Y)!1.85!(l3) $) -- ($ (Y)!2.4!(l3) $);
            \draw[dash pattern=on 1pt off 3pt] ($ (X2)!1.85!(l4) $) -- ($ (X2)!2.4!(l4) $);
            \draw[dash pattern=on 1pt off 3pt] ($ (A6)!1.85!(l5) $) -- ($ (A6)!2.4!(l5) $);
            \draw[dash pattern=on 1pt off 3pt] ($ (A5)!1.85!(l6) $) -- ($ (A5)!2.4!(l6) $);
            \draw[-] (A) -- ($ (A2)!1.4!(A) $);
            \draw[-] (C) -- ($ (A2)!1.4!(C) $);
            \draw[dash pattern=on 1pt off 3pt] ($ (A2)!1.4!(A) $) -- ($ (A2)!1.6!(A) $);
            \draw[dash pattern=on 1pt off 3pt] ($ (A2)!1.4!(C) $) -- ($ (A2)!1.6!(C) $);
        \end{tikzpicture}
    \end{minipage}
    \hspace{0.5cm}
    \begin{minipage}[c]{0.1\linewidth}
        {\vspace{-0.97cm}\Huge \begin{equation*}
            \leadsto
        \end{equation*}}
    \end{minipage}
    \begin{minipage}[c]{0.4\linewidth}
    \centering
        \begin{tikzpicture}[>=Stealth]
            \node (A) at (0,0) {$\P_2$};
            \node (C) at (-3,0) {$\P_1$};
            \node (V) at (-1.5,1.73*1.5) {$\P_3$};
            \node (W) at (-1.5,-1.73*1.5) {$\P_4$};

            \node (X2) at ($ (C)!0.5!(W) $) {$\E_5$};
            \node (Y) at ($ (C)!0.5!(V) $) {$\E_3$};

            \node (A2) at ($ (A)!0.5!(C) $) {$\E_1$};
            \node (A5) at ($ (A)!0.5!(V) $) {$\E_2$};
            \node (A6) at ($ (A)!0.5!(W) $) {$\E_4$};

            \node (l1) at (barycentric cs:A=1,C=1,V=1) {\color{gray}$\l_5$};
            \node (l2) at (barycentric cs:A=1,C=1,W=1) {\color{gray}$\l_6$};

            \node (l3) at ($ (l1)!2!(Y) $) {\color{gray}$\l_1$};
            \node (l4) at ($ (l2)!2!(X2) $) {\color{gray}$\l_2$};
            \node (l5) at ($ (l2)!2!(A6) $) {\color{gray}$\l_3$};
            \node (l6) at ($ (l1)!2!(A5) $) {\color{gray}$\l_4$};
            
            \draw[-, color=red] (A) -- (A2);
            \draw[-, color=red] (A) -- (A5);
            \draw[-, color=red] (A) -- (A6);
            \draw[-, color=red] (A2) -- (C);
            \draw[-, color=red] (A5) -- (V);
            \draw[-, color=red] (A6) -- (W);
            \draw[-, color=red] (C) -- (X2);
            \draw[-, color=red] (C) -- (Y);
            \draw[-, color=red] (V) -- (Y);
            \draw[-, color=red] (W) -- (X2);
            \draw[dash pattern=on 2pt off 3pt,color=gray] (A) -- (l1);
            \draw[dash pattern=on 2pt off 3pt,color=gray] (A) -- (l2);
            \draw[dash pattern=on 2pt off 3pt,color=gray] (C) -- (l1);
            \draw[dash pattern=on 2pt off 3pt,color=gray] (C) -- (l2);
            \draw[dash pattern=on 2pt off 3pt,color=gray] (V) -- (l1);
            \draw[dash pattern=on 2pt off 3pt,color=gray] (W) -- (l2);
            \draw[dash pattern=on 2pt off 3pt,color=gray] (C) -- (l3);
            \draw[dash pattern=on 2pt off 3pt,color=gray] (C) -- (l4);
            \draw[dash pattern=on 2pt off 3pt,color=gray] (A) -- (l5);
            \draw[dash pattern=on 2pt off 3pt,color=gray] (A) -- (l6);
            \draw[dash pattern=on 2pt off 3pt,color=gray] (V) -- (l3);
            \draw[dash pattern=on 2pt off 3pt,color=gray] (W) -- (l4);
            \draw[dash pattern=on 2pt off 3pt,color=gray] (W) -- (l5);
            \draw[dash pattern=on 2pt off 3pt,color=gray] (V) -- (l6);
            \draw[-, color=red] (C) -- ($ (W)!1.22!(C) $);
            \draw[-, color=red] (C) -- ($ (V)!1.22!(C) $);
            \draw[-, color=red] (A) -- ($ (W)!1.22!(A) $);
            \draw[-, color=red] (A) -- ($ (V)!1.22!(A) $);
            \draw[dash pattern=on 1pt off 3pt, color=red] ($ (W)!1.22!(C) $) -- ($ (W)!1.4!(C) $);
            \draw[dash pattern=on 1pt off 3pt, color=red] ($ (V)!1.22!(C) $) -- ($ (V)!1.4!(C) $);
            \draw[dash pattern=on 1pt off 3pt, color=red] ($ (W)!1.22!(A) $) -- ($ (W)!1.4!(A) $);
            \draw[dash pattern=on 1pt off 3pt, color=red] ($ (V)!1.22!(A) $) -- ($ (V)!1.4!(A) $);
            \draw[-, color=red] (V) -- (-2.15,1.73*1.5);
            \draw[-, color=red] (V) -- (-0.85,1.73*1.5);
            \draw[-, color=red] (W) -- (-2.15,-1.73*1.5);
            \draw[-, color=red] (W) -- (-0.85,-1.73*1.5);
            \draw[dash pattern=on 1pt off 3pt, color=red] (-2.15,1.73*1.5) -- (-2.7,1.73*1.5);
            \draw[dash pattern=on 1pt off 3pt, color=red] (-0.85,1.73*1.5) -- (-0.3,1.73*1.5);
            \draw[dash pattern=on 1pt off 3pt, color=red] (-2.15,-1.73*1.5) -- (-2.7,-1.73*1.5);
            \draw[dash pattern=on 1pt off 3pt, color=red] (-0.85,-1.73*1.5) -- (-0.3,-1.73*1.5);
            \draw[dash pattern=on 2pt off 3pt,color=gray] (l3) -- ($ (Y)!1.85!(l3) $);
            \draw[dash pattern=on 2pt off 3pt,color=gray] (l4) -- ($ (X2)!1.85!(l4) $);
            \draw[dash pattern=on 2pt off 3pt,color=gray] (l5) -- ($ (A6)!1.85!(l5) $);
            \draw[dash pattern=on 2pt off 3pt,color=gray] (l6) -- ($ (A5)!1.85!(l6) $);
            \draw[dash pattern=on 2pt off 3pt,color=gray] (A) -- ($ (A2)!1.5!(A) $);
            \draw[dash pattern=on 2pt off 3pt,color=gray] (C) -- ($ (A2)!1.5!(C) $);
        \end{tikzpicture}
    \end{minipage}
    \caption{Translating a quad-proof into a pure Menelaus-proof.}\label{fig: quad mene}
\end{figure}
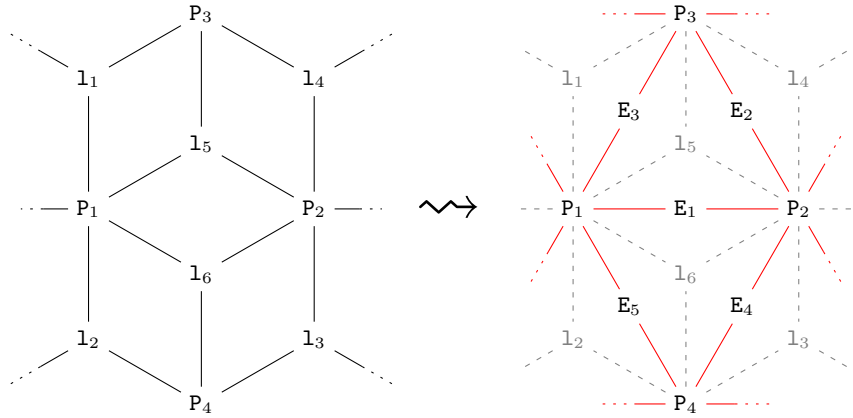

One way to directly translate quad-proofs into Ceva-Menelaus-proofs is actually already implicitly provided  in~\cite{FP}. There, a translation process to proofs using ``triangulations of closed oriented surfaces'' is described, that essentially converts a quad-proof into a Ceva-Menelaus-proof using \textit{only} Menelaus-triangles, however without mentioning the theorem of Menelaus explicitly. This is done by first iteratively inserting additional quadrangles in the places of ``line-vertices'' contained in more than $3$ quadrangles by splitting them into $2$ vertices,
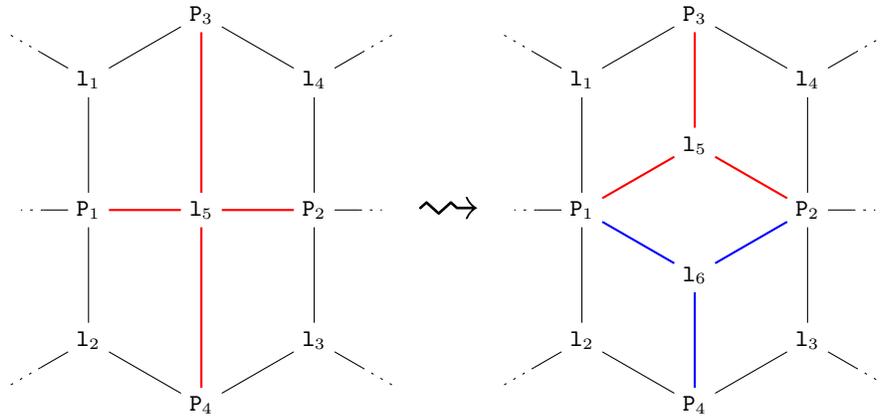
\begin{figure}[ht]
    \centering
    \hspace{-0.3cm}
    \begin{minipage}[c]{0.37\linewidth}
    \centering
        \begin{tikzpicture}[>=Stealth]
            \node (A) at (0,0) {$\P_2$};
            \node (C) at (-3,0) {$\P_1$};
            \node (V) at (-1.5,1.73*1.5) {$\P_3$};
            \node (W) at (-1.5,-1.73*1.5) {$\P_4$};

            \node (X2) at ($ (C)!0.5!(W) $) {};
            \node (Y) at ($ (C)!0.5!(V) $) {};

            \node (A2) at ($ (A)!0.5!(C) $) {};
            \node (A5) at ($ (A)!0.5!(V) $) {};
            \node (A6) at ($ (A)!0.5!(W) $) {};

            \node (l1) at (barycentric cs:A=1,C=1,V=1) {};
            \node (l2) at (barycentric cs:A=1,C=1,W=1) {};
            \node (lx) at (barycentric cs:A=1,C=1) {$\l_5$};

            \node (l3) at ($ (l1)!2!(Y) $) {$\l_1$};
            \node (l4) at ($ (l2)!2!(X2) $) {$\l_2$};
            \node (l5) at ($ (l2)!2!(A6) $) {$\l_3$};
            \node (l6) at ($ (l1)!2!(A5) $) {$\l_4$};
            \draw[thick, color = red] (A) -- (lx);
            \draw[thick, color = red] (C) -- (lx);
            \draw[thick, color = red] (V) -- (lx);
            \draw[thick, color = red] (W) -- (lx);
            \draw[-] (C) -- (l3);
            \draw[-](C) -- (l4);
            \draw[-] (A) -- (l5);
            \draw[-] (A) -- (l6);
            \draw[-] (V) -- (l3);
            \draw[-] (W) -- (l4);
            \draw[-] (W) -- (l5);
            \draw[-] (V) -- (l6);
            \draw[-] (l3) -- ($ (Y)!1.85!(l3) $);
            \draw[-] (l4) -- ($ (X2)!1.85!(l4) $);
            \draw[-] (l5) -- ($ (A6)!1.85!(l5) $);
            \draw[-] (l6) -- ($ (A5)!1.85!(l6) $);
            \draw[dash pattern=on 1pt off 3pt] ($ (Y)!1.85!(l3) $) -- ($ (Y)!2.4!(l3) $);
            \draw[dash pattern=on 1pt off 3pt] ($ (X2)!1.85!(l4) $) -- ($ (X2)!2.4!(l4) $);
            \draw[dash pattern=on 1pt off 3pt] ($ (A6)!1.85!(l5) $) -- ($ (A6)!2.4!(l5) $);
            \draw[dash pattern=on 1pt off 3pt] ($ (A5)!1.85!(l6) $) -- ($ (A5)!2.4!(l6) $);
            \draw[-] (A) -- ($ (lx)!1.4!(A) $);
            \draw[-] (C) -- ($ (lx)!1.4!(C) $);
            \draw[dash pattern=on 1pt off 3pt] ($ (lx)!1.4!(A) $) -- ($ (lx)!1.6!(A) $);
            \draw[dash pattern=on 1pt off 3pt] ($ (lx)!1.4!(C) $) -- ($ (lx)!1.6!(C) $);
        \end{tikzpicture}
    \end{minipage}
    \hspace{0.5cm}
    \begin{minipage}[c]{0.1\linewidth}
        {\vspace{-0.97cm}\Huge \begin{equation*}
            \leadsto
        \end{equation*}}
    \end{minipage}
    \begin{minipage}[c]{0.4\linewidth}
        \begin{tikzpicture}[>=Stealth]
            \node (A) at (0,0) {$\P_2$};
            \node (C) at (-3,0) {$\P_1$};
            \node (V) at (-1.5,1.73*1.5) {$\P_3$};
            \node (W) at (-1.5,-1.73*1.5) {$\P_4$};

            \node (X2) at ($ (C)!0.5!(W) $) {};
            \node (Y) at ($ (C)!0.5!(V) $) {};

            \node (A2) at ($ (A)!0.5!(C) $) {};
            \node (A5) at ($ (A)!0.5!(V) $) {};
            \node (A6) at ($ (A)!0.5!(W) $) {};

            \node (l1) at (barycentric cs:A=1,C=1,V=1) {$\l_5$};
            \node (l2) at (barycentric cs:A=1,C=1,W=1) {$\l_6$};
            \node (lx) at (barycentric cs:A=1,C=1) {};

            \node (l3) at ($ (l1)!2!(Y) $) {$\l_1$};
            \node (l4) at ($ (l2)!2!(X2) $) {$\l_2$};
            \node (l5) at ($ (l2)!2!(A6) $) {$\l_3$};
            \node (l6) at ($ (l1)!2!(A5) $) {$\l_4$};
            \draw[thick, color = red] (A) -- (l1);
            \draw[thick, color=blue] (A) -- (l2);
            \draw[thick, color = red] (C) -- (l1);
            \draw[thick, color=blue] (C) -- (l2);
            \draw[thick, color = red] (V) -- (l1);
            \draw[thick, color=blue] (W) -- (l2);
            \draw[-] (C) -- (l3);
            \draw[-](C) -- (l4);
            \draw[-] (A) -- (l5);
            \draw[-] (A) -- (l6);
            \draw[-] (V) -- (l3);
            \draw[-] (W) -- (l4);
            \draw[-] (W) -- (l5);
            \draw[-] (V) -- (l6);
            \draw[-] (l3) -- ($ (Y)!1.85!(l3) $);
            \draw[-] (l4) -- ($ (X2)!1.85!(l4) $);
            \draw[-] (l5) -- ($ (A6)!1.85!(l5) $);
            \draw[-] (l6) -- ($ (A5)!1.85!(l6) $);
            \draw[dash pattern=on 1pt off 3pt] ($ (Y)!1.85!(l3) $) -- ($ (Y)!2.4!(l3) $);
            \draw[dash pattern=on 1pt off 3pt] ($ (X2)!1.85!(l4) $) -- ($ (X2)!2.4!(l4) $);
            \draw[dash pattern=on 1pt off 3pt] ($ (A6)!1.85!(l5) $) -- ($ (A6)!2.4!(l5) $);
            \draw[dash pattern=on 1pt off 3pt] ($ (A5)!1.85!(l6) $) -- ($ (A5)!2.4!(l6) $);
            \draw[-] (A) -- ($ (lx)!1.4!(A) $);
            \draw[-] (C) -- ($ (lx)!1.4!(C) $);
            \draw[dash pattern=on 1pt off 3pt] ($ (lx)!1.4!(A) $) -- ($ (lx)!1.6!(A) $);
            \draw[dash pattern=on 1pt off 3pt] ($ (lx)!1.4!(C) $) -- ($ (lx)!1.6!(C) $);
        \end{tikzpicture}
    \end{minipage}
    \caption{Splitting a line-vertex by inserting a new quadrangle.}
\end{figure}
extending the quad-proof to a proof of a more general projective incidence theorem (This more general theorem collapses to the original theorem in case the $2$ lines corresponding to the line-vertices coincide). Then, if all line-vertices are contained in exactly $3$ quadrangles, the quadrilateral tiling can be translated into a triangulation by splitting all quadrangles in half and gluing triples of halves containing the same line-vertex to triangles, as visualized in Figure~\ref{fig: quad mene}. It is easy to see that the underlying incidence structure of these triangles is exactly that of the theorem of Menelaus.
This procedure can simply be reversed to get from any pure Menelaus-proof to a quad-proof, so quad-proofs and pure Menelaus-proofs are equally powerful.

\newpage
\subsection{Ceva vs Menelaus}
\begin{figure}[h]
    \centering
    {\begin{minipage}[c]{0.28\linewidth}
    \centering
    \begin{tikzpicture}[]
        \coordinate [label=left:{$\A$}] (A) at (0,0);
        \coordinate [label=right:{$\B$}] (B) at (3,0);
        \coordinate [label=above:{$\C$}] (C) at (1,2.5);
        \coordinate [label={[label distance=0.15cm]right:$\D$}] (D) at (1.25,1);
    
        \path [name path=A--B] (A) -- (B);
        \draw[very thick] (A) -- (B);
        \path [name path=A--C] (A) -- (C);
        \draw[very thick] (A) -- (C);
        \path [name path=C--B] (C) -- (B);
        \draw[very thick] (C) -- (B);
    
        \draw [] (A) -- ($ (A)!1.5!(D) $) coordinate (E);
        \path [name path=A--E] (A) -- (E);
        \draw [] (B) -- ($ (B)!1.4!(D) $) coordinate (F);
        \path [name path=B--F] (B) -- (F);
        \draw [] (C) -- ($ (C)!1.7!(D) $) coordinate (G);
        \path [name path=C--G] (C) -- (G);
        
        \path [name intersections={of=C--B and A--E,by={[label=above right:$\X$]X}}];
        \path [name intersections={of=A--C and B--F,by={[label=above left:$\Y$]Y}}];
        \path [name intersections={of=A--B and C--G,by={[label=below:$\Z$]Z}}];
        \foreach \point in {A,B,C} {
            \fill [black,opacity=1] (\point) circle (2pt);
        }
        \path[fill=white,draw=black, thick] (D) circle[radius=2pt];
        \foreach \point in {X,Y,Z} {
            \path[fill=white,draw=black, thick] (\point) circle[radius=1.5pt];
        }
    \end{tikzpicture}
    \end{minipage}
    {\huge \textbf{=}}
    \begin{minipage}[c]{0.28\linewidth}
    \centering
    \begin{tikzpicture}[]
        \coordinate [label=left:{$\A$}] (A) at (0,0);
        \coordinate [label=right:{$\B$}] (B) at (3,0);
        \coordinate [label=above:{$\C$}] (C) at (1,2.5);
        \coordinate [label={[label distance=0.15cm]right:$\D$}] (D) at (1.25,1);
    
        \path [name path=A--B] (A) -- (B);
        \draw[very thick] (A) -- (B);
        \path [name path=A--C] (A) -- (C);
        \draw[very thick] (A) -- (C);
    
        \draw [] (B) -- ($ (B)!1.4!(D) $) coordinate (F);
        \path [name path=B--F] (B) -- (F);
        \draw [very thick] (C) -- ($ (C)!1.7!(D) $) coordinate (G);
        \path [name path=C--G] (C) -- (G);
        
        \path [name intersections={of=A--B and C--G,by={[label=below:$\Z$]Z}}];
        \path [name intersections={of=A--C and B--F,by={[label=above left:$\Y$]Y}}];
    
        \foreach \point in {A,Z,C} {
            \fill [black,opacity=1] (\point) circle (2pt);
        }
        \foreach \point in {D,Y,B} {
            \path[fill=white,draw=black, thick] (\point) circle[radius=1.5pt];
        }
    
    \end{tikzpicture}
    \end{minipage}
    {\huge \textbf{+}}
    \begin{minipage}[c]{0.3\linewidth}
    \centering
    \begin{tikzpicture}[]
        \coordinate [label=left:{$\A$}] (A) at (0,0);
        \coordinate [label=right:{$\B$}] (B) at (3,0);
        \coordinate [label=above:{$\C$}] (C) at (1,2.5);
        \coordinate [label={[label distance=0.15cm]right:$\D$}] (D) at (1.25,1);
    
        \path [name path=A--B] (A) -- (B);
        \draw[very thick] (A) -- (B);
        \path [name path=C--B] (C) -- (B);
        \draw[very thick] (C) -- (B);
    
        \draw [] (A) -- ($ (A)!1.5!(D) $) coordinate (E);
        \path [name path=A--E] (A) -- (E);
        \draw [very thick] (C) -- ($ (C)!1.7!(D) $) coordinate (G);
        \path [name path=C--G] (C) -- (G);
        
        \path [name intersections={of=C--B and A--E,by={[label=above right:$\X$]X}}];
        \path [name intersections={of=A--B and C--G,by={[label=below:$\Z$]Z}}];
    
        \foreach \point in {Z,B,C} {
            \fill [black,opacity=1] (\point) circle (2pt);
        }
        \foreach \point in {D,X,A} {
            \path[fill=white,draw=black, thick] (\point) circle[radius=1.5pt];
        }
    
    \end{tikzpicture}
    \end{minipage}}
    \caption{\centering{Proof of Ceva's theorem using Menelaus's theorem twice.}}\label{fig: Ceva = 2 Mene}
\end{figure}
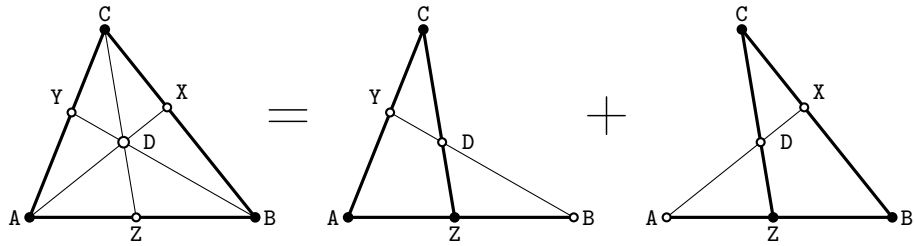

Naturally the question arises, if pure Menelaus-proofs are already as powerful as Ceva-Menelaus-proofs. It turns out, that Ceva’s theorem can actually be combined from two applications
of Menelaus's theorem, visualized in Figure~\ref{fig: Ceva = 2 Mene}.
However, a simple argument shows that replacing single Ceva-triangles in a Ceva-Menelaus-proof
by any number Menelaus-triangles is not possible:
\begin{itemize}
    \item[1.] The total number of triangles is even: Each triangle has $3$ edges and each edge is contained in $2$ triangles, so $3 \cdot \# triangles = 2 \cdot \# edges$.
    \item[2.] The total number of Menelaus-triangles (thus by 1.\ also Ceva-triangles) is even: Each Menelaus-triangle adds a factor ``$-1$'', but the product of all equations
    must be $1$.
    \item[3.] Replacing a single Ceva-triangle by Menelaus-triangles makes the number of Ceva-triangles odd.
\end{itemize}
Nevertheless, some replacements are possible. If a Ceva-Menelaus-proof contains $2$ adjacent Ceva-triangles, both can be split simultaneously into $4$ adjacent Menelaus-triangles sharing the former edge point
of the edge between the two Ceva-triangles as a vertex point, as shown below.
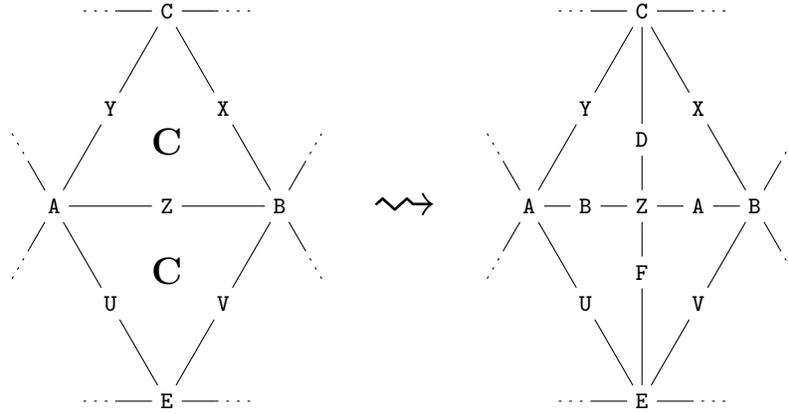
\begin{figure}[h]
    \centering
    \begin{minipage}[c]{0.4\linewidth}
        \centering
        \begin{tikzpicture}[>=Stealth]
            \node (A) at (0,0) {$\B$};
            \node (C) at (-3,0) {$\A$};
            \node (V) at (-1.5,1.73*1.5) {$\C$};
            \node (W) at (-1.5,-1.73*1.5) {$\E$};

            \node (X2) at ($ (C)!0.5!(W) $) {$\U$};
            \node (Y) at ($ (C)!0.5!(V) $) {$\Y$};

            \node (A2) at ($ (A)!0.5!(C) $) {$\Z$};
            \node (A5) at ($ (A)!0.5!(V) $) {$\X$};
            \node (A6) at ($ (A)!0.5!(W) $) {$\V$};

            \node (A7) at (barycentric cs:A=1,C=1,V=1) {$\mathlarger{\mathlarger{\mathlarger{\mathlarger{\mathbf{C}}}}}$};
            \node (A7) at (barycentric cs:A=1,C=1,W=1) {$\mathlarger{\mathlarger{\mathlarger{\mathlarger{\mathbf{C}}}}}$};
            
            \draw[-] (A) -- (A2);
            \draw[-] (A) -- (A5);
            \draw[-] (A) -- (A6);
            \draw[-] (A2) -- (C);
            \draw[-] (A5) -- (V);
            \draw[-] (A6) -- (W);
            \draw[-] (C) -- (X2);
            \draw[-] (C) -- (Y);
            \draw[-] (V) -- (Y);
            \draw[-] (W) -- (X2);
            \draw[-] (C) -- ($ (W)!1.22!(C) $);
            \draw[-] (C) -- ($ (V)!1.22!(C) $);
            \draw[-] (A) -- ($ (W)!1.22!(A) $);
            \draw[-] (A) -- ($ (V)!1.22!(A) $);
            \draw[dash pattern=on 1pt off 3pt] ($ (W)!1.22!(C) $) -- ($ (W)!1.4!(C) $);
            \draw[dash pattern=on 1pt off 3pt] ($ (V)!1.22!(C) $) -- ($ (V)!1.4!(C) $);
            \draw[dash pattern=on 1pt off 3pt] ($ (W)!1.22!(A) $) -- ($ (W)!1.4!(A) $);
            \draw[dash pattern=on 1pt off 3pt] ($ (V)!1.22!(A) $) -- ($ (V)!1.4!(A) $);
            \draw[-] (V) -- (-2.15,1.73*1.5);
            \draw[-] (V) -- (-0.85,1.73*1.5);
            \draw[-] (W) -- (-2.15,-1.73*1.5);
            \draw[-] (W) -- (-0.85,-1.73*1.5);
            \draw[dash pattern=on 1pt off 3pt] (-2.15,1.73*1.5) -- (-2.7,1.73*1.5);
            \draw[dash pattern=on 1pt off 3pt] (-0.85,1.73*1.5) -- (-0.3,1.73*1.5);
            \draw[dash pattern=on 1pt off 3pt] (-2.15,-1.73*1.5) -- (-2.7,-1.73*1.5);
            \draw[dash pattern=on 1pt off 3pt] (-0.85,-1.73*1.5) -- (-0.3,-1.73*1.5);
        \end{tikzpicture}
    \end{minipage}
    \begin{minipage}[c]{0.1\linewidth}
        {\vspace{-0.97cm}\Huge \begin{equation*}
            \leadsto
        \end{equation*}}
    \end{minipage}
    \begin{minipage}[c]{0.4\linewidth}
        \centering
        \begin{tikzpicture}[>=Stealth]
            \node (A) at (0,0) {$\B$};
            \node (C) at (-3,0) {$\A$};
            \node (V) at (-1.5,1.73*1.5) {$\C$};
            \node (W) at (-1.5,-1.73*1.5) {$\E$};

            \node (X2) at ($ (C)!0.5!(W) $) {$\U$};
            \node (Y) at ($ (C)!0.5!(V) $) {$\Y$};

            \node (A2) at ($ (A)!0.5!(C) $) {$\Z$};
            \node (A5) at ($ (A)!0.5!(V) $) {$\X$};
            \node (A6) at ($ (A)!0.5!(W) $) {$\V$};

            \node (B1) at ($ (A)!0.5!(A2) $) {$\A$};
            \node (B2) at ($ (A2)!0.5!(C) $) {$\B$};

            \node (A7) at ($ (V)!0.66!(A2) $) {$\D$};
            \node (A8) at ($ (W)!0.66!(A2) $) {$\F$};

            \draw[-] (A2) -- (A7);
            \draw[-] (A2) -- (A8);
            \draw[-] (V) -- (A7);
            \draw[-] (W) -- (A8);
            \draw[-] (A) -- (B1);
            \draw[-] (B1) -- (A2);
            \draw[-] (C) -- (B2);
            \draw[-] (B2) -- (A2);
            \draw[-] (A) -- (A5);
            \draw[-] (A) -- (A6);
            \draw[-] (A5) -- (V);
            \draw[-] (A6) -- (W);
            \draw[-] (C) -- (X2);
            \draw[-] (C) -- (Y);
            \draw[-] (V) -- (Y);
            \draw[-] (W) -- (X2);
            \draw[-] (C) -- ($ (W)!1.22!(C) $);
            \draw[-] (C) -- ($ (V)!1.22!(C) $);
            \draw[-] (A) -- ($ (W)!1.22!(A) $);
            \draw[-] (A) -- ($ (V)!1.22!(A) $);
            \draw[dash pattern=on 1pt off 3pt] ($ (W)!1.22!(C) $) -- ($ (W)!1.4!(C) $);
            \draw[dash pattern=on 1pt off 3pt] ($ (V)!1.22!(C) $) -- ($ (V)!1.4!(C) $);
            \draw[dash pattern=on 1pt off 3pt] ($ (W)!1.22!(A) $) -- ($ (W)!1.4!(A) $);
            \draw[dash pattern=on 1pt off 3pt] ($ (V)!1.22!(A) $) -- ($ (V)!1.4!(A) $);
            \draw[-] (V) -- (-2.15,1.73*1.5);
            \draw[-] (V) -- (-0.85,1.73*1.5);
            \draw[-] (W) -- (-2.15,-1.73*1.5);
            \draw[-] (W) -- (-0.85,-1.73*1.5);
            \draw[dash pattern=on 1pt off 3pt] (-2.15,1.73*1.5) -- (-2.7,1.73*1.5);
            \draw[dash pattern=on 1pt off 3pt] (-0.85,1.73*1.5) -- (-0.3,1.73*1.5);
            \draw[dash pattern=on 1pt off 3pt] (-2.15,-1.73*1.5) -- (-2.7,-1.73*1.5);
            \draw[dash pattern=on 1pt off 3pt] (-0.85,-1.73*1.5) -- (-0.3,-1.73*1.5);
        \end{tikzpicture}
    \end{minipage}
    \caption{Splitting $2$ adjacent Ceva-triangles.}\label{fig: 2 Ceva = 4 Mene}
\end{figure}
While this does not allow translating every Ceva-Menelaus-proof into a pure Menelaus-proof, it at atleast shows that pure Menelaus-proofs must be stronger than pure Ceva-proofs.

\subsection{Summary of the Results}
In total we get the following relations between the proving methods final polynomials ($\mathbf{P}$), binomial proofs ($\mathbf{B}$), Ceva-Menelaus-proofs ($\mathbf{CM}$),
quad-proofs ($\mathbf{Q}$), pure Menelaus-proofs ($\mathbf{M}$) and pure Ceva-proofs ($\mathbf{C}$), where ``$\Rlongarrow$'' means ``at least as powerful as''.
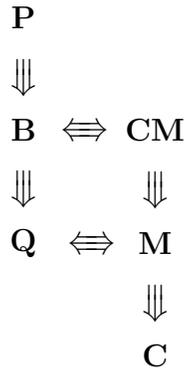
\begin{figure}[h!]
    \centering
    \begin{tikzpicture}
        \node (A) at (0,0) {$\mathlarger{\mathlarger{\mathlarger{\mathbf{P}}}}$};
        \node (B) at (0,-1.5) {$\mathlarger{\mathlarger{\mathlarger{\mathbf{B}}}}$};
        \node (C) at (1.75,-1.5) {$\mathlarger{\mathlarger{\mathlarger{\mathbf{CM}}}}$};
        \node (D) at (0,-3){$\mathlarger{\mathlarger{\mathlarger{\mathbf{Q}}}}$};
        \node (E) at (1.75,-3){$\mathlarger{\mathlarger{\mathlarger{\mathbf{M}}}}$};
        \node (F) at (1.75,-4.5){$\mathlarger{\mathlarger{\mathlarger{\mathbf{C}}}}$};
        \path [decorate, decoration={text along path, text format delimiters={|}{|}, text align={center}, raise=-3pt, text={|\footnotesize| ${ \mathlarger{\mathlarger{\mathlarger{\Rlongarrow}}} }$ }}] (A) -- (B);
        \path [decorate, decoration={text along path, text format delimiters={|}{|}, text align={center}, raise=-3pt, text={|\footnotesize| ${ \mathlarger{\mathlarger{\mathlarger{\LRarrow}}} }$ }}] (B) -- (C);
        \path [decorate, decoration={text along path, text format delimiters={|}{|}, text align={center}, raise=-3pt, text={|\footnotesize| ${ \mathlarger{\mathlarger{\mathlarger{\Rlongarrow}}} }$ }}] (B) -- (D);
        \path [decorate, decoration={text along path, text format delimiters={|}{|}, text align={center}, raise=-3pt, text={|\footnotesize| ${ \mathlarger{\mathlarger{\mathlarger{\LRarrow}}} }$ }}] (D) -- (E);
        \path [decorate, decoration={text along path, text format delimiters={|}{|}, text align={center}, raise=-3pt, text={|\footnotesize| ${ \mathlarger{\mathlarger{\mathlarger{\Rlongarrow}}} }$ }}] (C) -- (E);
        \path [decorate, decoration={text along path, text format delimiters={|}{|}, text align={center}, raise=-3pt, text={|\footnotesize| ${ \mathlarger{\mathlarger{\mathlarger{\Rlongarrow}}} }$ }}] (E) -- (F);
    \end{tikzpicture}
    \caption{Hierarchy of proving techniques.}
\end{figure}
\newpage

\section{Two illustrative examples}

We now will exemplify the translation between the proving methods by two simple, but suitably examples: Desargues' and Pappus' Theorems. We purposely chose the most elementary theorems, since there the connections between the approaches can be seen in a most transparent way.

\subsection{Desargues' Theorem}

Desargues' Theorem can be --- and has been --- shown using all presented techniques. In particular, it can be proven using:
\begin{itemize}
    \item[1.] A tiling of the sphere consisting of $6$ quadrangles (a cube) as shown in~\cite{FP},
    \item[2.] A tiling of the sphere consisting of $4$ Menelaus-triangles (a tetrahedron) as shown in~\cite{RG1} or
    \item[3.] $10$ binomial equations coming from $10$ collinear point triples and Grassmann-Plücker-relations also shown in~\cite{RG1,RG2}.
\end{itemize} 
\begin{figure}
    \centering
    \begin{minipage}[c]{0.45\linewidth}
        \centering
        \begin{tikzpicture}[>=Stealth]
            \node[regular polygon,regular polygon sides=4,minimum size=5.5cm] (r) at (0,0) {};
            \node (r1) at (r.corner 1) {$\C$};
            \node (r2) at (r.corner 2) {$\l_1$};
            \node (r3) at (r.corner 3) {$\A$};
            \node (r4) at (r.corner 4) {$\l_2$};

            \node[regular polygon,regular polygon sides=4,minimum size=2.5cm] (s) at (0,-0.5) {};
            \node (s1) at (s.corner 1) {$\l_3$};
            \node (s2) at (s.corner 2) {$\B$};
            \node (s3) at (s.corner 3) {$\l_4$};
            \node (s4) at (s.corner 4) {$\D$};
    
            \draw[-] (r1) -- (r2);
            \draw[-] (r2) -- (r3);
            \draw[-] (r3) -- (r4);
            \draw[-] (r4) -- (r1);

            \draw[-] (s1) -- (s2);
            \draw[-] (s2) -- (s3);
            \draw[-] (s3) -- (s4);
            \draw[-] (s4) -- (s1);
            
            \draw[-] (s1) -- (r1);
            \draw[-] (s2) -- (r2);
            \draw[-] (s3) -- (r3);
            \draw[-] (s4) -- (r4);

            \draw[color=red,thick,dash pattern=on 2pt off 3pt] (s2) -- (r1);
            \draw[color=red,thick,dash pattern=on 2pt off 3pt] (s2) -- (r3);
            \draw[color=red,thick,dash pattern=on 2pt off 3pt] (s4) -- (r1);
            \draw[color=red,thick,dash pattern=on 2pt off 3pt] (s4) -- (r3);
            \draw[color=red,thick,dash pattern=on 2pt off 3pt] (s4) -- (s2);
            \draw[color=red,thick,dash pattern=on 2pt off 3pt] (r1) -- (r3);
        \end{tikzpicture}
    \end{minipage}
    \hspace{0.5cm}
    \begin{minipage}[c]{0.45\linewidth}
        \centering
        \begin{tikzpicture}[]
            \coordinate [label=above right:{$\A$}] (A) at (0.8*1,0.8*4);
            \coordinate [label=left:{$\B$}] (B) at (0,0);
            \coordinate [label=right:{$\C$}] (C) at (0.8*5,0);
            \coordinate [label=above left:{$\D$}] (D) at (0.8*4,0.8*4);
        
            \path [name path=A--B] (A) -- (B);
            \draw[very thick] (A) -- (B);
            \path [name path=B--C] (B) -- (C);
            \draw[very thick] (B) -- (C);
            \path [name path=C--D] (C) -- (D);
            \draw[very thick] (C) -- (D);
            \path [name path=D--A] (D) -- (A);
            \draw[very thick] (D) -- (A);
        
            \draw [very thick] (C) -- ($ (C)!1.4!(A) $) coordinate (j);    
            \draw [] (C) -- ($ (C)!1.4!(A) $) coordinate [label=above right:$\Y$] (E);
            \path [name path=A--E] (A) -- (E);
            
            \draw [] (E) -- ($ (B)!0.5!(C) $) coordinate [label=below:$\X$] (I);
            \path [name path=E--I] (E) -- (I);
    
            \path [name intersections={of=A--B and E--I,by={[label=below left:$\Z$]H}}];

            \draw [] (E) -- ($ (D)!0.2!(C) $) coordinate [label=below right:$\W$] (G);
            \path [name path=E--G] (E) -- (G);
    
            \path [name intersections={of=D--A and E--G,by={[label=above:$\U$]F}}];
    
            \coordinate [] (k) at ($ (H)!2.7!(F) $);
            \path [name path=H--F] (H) -- (k);
            
            \coordinate [] (l) at ($ (I)!1.8!(G) $) ;
            \path [name path=I--G] (I) -- (l);
    
            \path [name intersections={of=H--F and I--G,by={[label=above left:$\V$]X}}];
            \draw[-] (X) -- (H);
            \draw[-] (X) -- (I);
            \draw [very thick] (B) -- (X); 
            
            \foreach \point in {A,B,C,D} {
                \fill [black,opacity=1] (\point) circle (2pt);
            }
    
            \foreach \point in {E,I,H,F,G,X} {
                \path[fill=white,draw=black, thick] (\point) circle[radius=1.5pt];
            }

            \coordinate [label=left:{$\l_1$}] (C) at ($ (E)!0.2!(I) $);
            \coordinate [label=above:{$\l_2$}] (C) at ($ (E)!0.25!(G) $);
            \coordinate [label=right:{$\l_3$}] (C) at ($ (X)!0.25!(I) $);
            \coordinate [label=above:{$\l_4$}] (C) at ($ (X)!0.3!(H) $);
        \end{tikzpicture}
    \end{minipage}
    \caption{Tiling proofs for Desargues' theorem.}\label{fig: Desargue}
\end{figure}
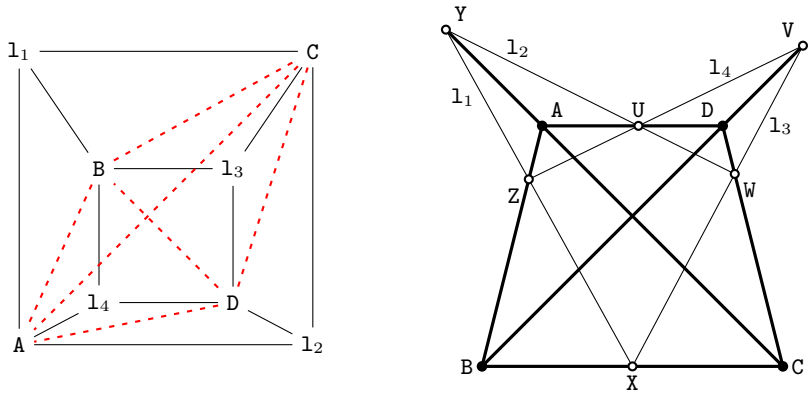
We start with a the quad-proof whose combinatorics corresponds to the cube (Figure~\ref{fig: Desargue}, left). The eight vertices of the cube have to be considered as the vertices of a bipartite graph. Half of them ($\A,\ldots,\D$) representing points and the other half ($\l_1,\ldots,\l_4$) representing lines. The consistency condition for each of the six faces corresponds to one incidence. For instance the face with vertices $\B,\C,\l_1,\l_3$ corresponds to the point $\X$ in the Figure~\ref{fig: Desargue} (right) at which the lines $\l_1,\l_3, \overline{\B\C}$ meet. In that sense every face of the cube gets associated with a point in Desargues' Theorem and $\A\ldots\B, \U\ldots \Z$ form the ten points of the theorem.
\begin{figure}[t]
    \centering
    \begin{tikzpicture}
        \draw (0, 0) node[inner sep=0] {\includegraphics[width=12cm]{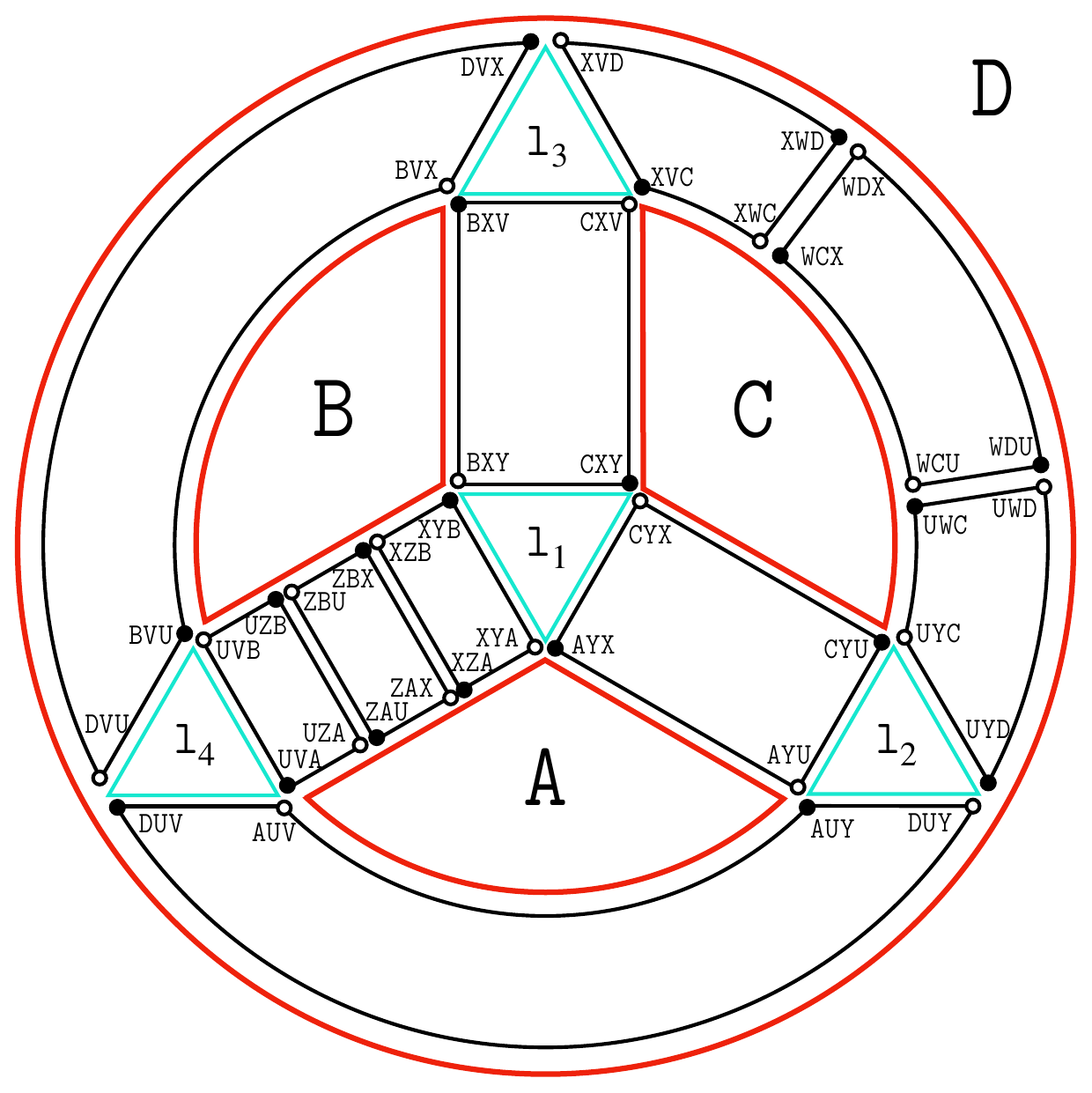}};
    \end{tikzpicture}
    \caption{The big picture.}\label{fig: the big picture}
\end{figure}

We now translate the coherent face conditions into a collection of binomial determinant expressions. For that we must express the lines as spans of two points. In our example for each line we take any two  of the three points incident to it. While the approach from Section~\ref{section: extract bin} would actually generate $18$ binomial equations from the quad-proof, the specific combinatorial structure allows us to 
``take a shortcut'' and translate four of the quadrangles using only {\it one} binomial equation each, by describing the lines directly as joins of the incidence points of adjacent quadrangles. 
We get
\begin{align*}
    \text{coherent quadrangle } \A,\l_2,\C,\l_1 \Longleftrightarrow &\ [\A,\Y,\X][\C,\Y,\U] = [\A,\Y,\U][\C,\Y,\X] \\
    \text{coherent quadrangle } \C,\l_3,\B,\l_1 \Longleftrightarrow &\ [\C,\X,\Y][\B,\X,\V] = [\C,\X,\V][\B,\X,\Y] \\
    \text{coherent quadrangle } \B,\l_3,\D,\l_4 \Longleftrightarrow &\ [\B,\V,\U][\D,\V,\X] = [\B,\V,\X][\D,\V,\U] \\
    \text{coherent quadrangle } \D,\l_2,\A,\l_4 \Longleftrightarrow &\ [\D,\U,\V][\A,\U,\Y] = [\D,\U,\Y][\A,\U,\V]
\end{align*}
for four of the faces. We collect the remaining brackets by translating the equations coming from the two quadrangles using the method from~\ref{section: extract bin}.
\begin{align*}
    \text{coherent quadrangle } \A,\l_1,\B,\l_4 \Longleftrightarrow &\ [\A,\U,\V][\B,\X,\Y] = [\A,\X,\Y][\B,\U,\V] \\
                                                \Longleftrightarrow &\ [\U,\V,\A][\U,\Z,\B] = [\U,\V,\B][\U,\Z,\A] \text{ and}\\
                                                                    &\ [\X,\Y,\B][\X,\Z,\A] = [\X,\Y,\A][\X,\Z,\B] \text{ and}\\
                                                                    &\ [\Z,\A,\U][\Z,\B,\X] = [\Z,\A,\X][\Z,\B,\U] \\
    \text{coherent quadrangle } \C,\l_2,\D,\l_3 \Longleftrightarrow &\ [\C,\X,\V][\D,\U,\Y] = [\C,\U,\Y][\D,\X,\V] \\
                                                \Longleftrightarrow &\ [\X,\V,\C][\X,\W,\D] = [\X,\V,\D][\X,\W,\C] \text{ and}\\
                                                                    &\ [\U,\Y,\D][\U,\W,\C] = [\U,\Y,\C][\U,\W,\D] \text{ and}\\
                                                                    &\ [\W,\C,\X][\W,\D,\U] = [\W,\C,\U][\W,\D,\X]
\end{align*}

This results in $10$ binomial equations cancelling each other out. In this collection 
of 10 binomial equations each equation encodes exactly one collinearity of the underlying Desargues' Theorem. Thus we have a binomial proof. The cancellation pattern
in the quad-proof directly translates to the proving structure in the binomial proof.

Figure~\ref{fig: the big picture} shows the situation in a more wholistic way.
There ten quadrilateral regions are shown (black). Each  vertex corresponds to the occurrence of one bracket in one binomial equation. Black vertices correspond to brackets occurring on the left. White vertices occur on the right of the equations. The vertices are labelled with the letters of the corresponding  bracket. The cancellation property corresponds to the fact that whenever a black and a white vertex meet they are labelled by the same letters.

Each black edge can be considered as  a ratio of the form $[\mathtt{OPQ}]/[\mathtt{OPR}]$ (dividing its black end by its white end). Thus each edge represents a length ratio. There are four red regions in the image labelled $\A,\ldots,\D$. They represent the {\it points} in the original quad-proof. 
Note that each vertex adjacent to such a region contains the corresponding letter. 
The blue triangles represent the {\it lines} in the quad-proof. The vertices adjacent to such a triangle all share two points. These two points span the corresponding line.

Finally, we can read off the corresponding Ceva-Menelaus proof. For this we only have to focus on the green triangles. They will be considered as the triangles in a Menelaus proof. Each pair of triangles is glued along two edges that are connected by a long black region (either consisting of one or of three quadrilateral cells). The two segments at the end of such a region represent the same length ratio. For instance the region connecting $\l_1$ and $\l_4$ connects to $\l_1$ by $[\mathtt{XYB}]/[\mathtt{XYA}]$ and connects to  $\l_4$ by $[\mathtt{UVA}]/[\mathtt{UVB}]$.  Both fractions represent the same length ratio  ${\overrightarrow{|\Z\A|}}/{\overrightarrow{|\Z\B|}}$ (one of them as a reciprocal of the other). This is the ratio in which the lines $\l_1$ and $\l_4$ cut the line $\overline{\mathtt{AB}}$.

If we  collect the products of the rations for the $4$ Menelaus triangles $\A\B\C$ with line $\l_1$, $\A\C\D$ with line $\l_2$, $\B\C\D$ with line $\l_3$ and $\A\B\D$ with line $\l_4$, then their Menelaus-equations are equivalent to the obvious identities
\begin{align*}
    \A\B\C \text{ with } \l_1
    &\Longleftrightarrow &
    \ \frac{[\X,\Y,\A]}{[\X,\Y,\B]}\cdot\frac{[\X,\Y,\B]}{[\X,\Y,\C]}\cdot\frac{[\X,\Y,\C]}{[\X,\Y,\A]} = 1 
 &\Longleftrightarrow&
   \frac {\overrightarrow{|\A\Z|}}{\overrightarrow{|\Z\B|}}\cdot
   \frac {\overrightarrow{|\B\X|}}{\overrightarrow{|\X\C|}}\cdot
   \frac {\overrightarrow{|\C\Y|}}{\overrightarrow{|\Y\A|}}=-1&,
    \\
    \A\C\D \text{ with } \l_2 
    &\Longleftrightarrow &
    \ \frac{[\U,\Y,\A]}{[\U,\Y,\C]}\cdot\frac{[\U,\Y,\C]}{[\U,\Y,\D]}\cdot\frac{[\U,\Y,\D]}{[\U,\Y,\A]} = 1
    & \Longleftrightarrow &
   \frac {\overrightarrow{|\A\Y|}}{\overrightarrow{|\Y\C|}}\cdot
   \frac {\overrightarrow{|\C\W|}}{\overrightarrow{|\W\D|}}\cdot
   \frac {\overrightarrow{|\D\U|}}{\overrightarrow{|\U\A|}}=-1&, \\
    \B\C\D \text{ with } \l_3& \Longleftrightarrow &\ \frac{[\X,\W,\C]}{[\X,\W,\B]}\cdot\frac{[\X,\W,\B]}{[\X,\W,\D]}\cdot\frac{[\X,\W,\D]}{[\X,\W,\C]} = 1 & \Longleftrightarrow &
   \frac {\overrightarrow{|\C\X|}}{\overrightarrow{|\X\B|}}\cdot
   \frac {\overrightarrow{|\B\V|}}{\overrightarrow{|\V\D|}}\cdot
   \frac {\overrightarrow{|\D\W|}}{\overrightarrow{|\W\C|}}=-1&, \\
    \A\B\D \text{ with } \l_4 
    &\Longleftrightarrow &
    \ \frac{[\U,\V,\B]}{[\U,\V,\A]}\cdot\frac{[\U,\V,\A]}{[\U,\V,\D]}\cdot\frac{[\U,\V,\D]}{[\U,\V,\B]} = 1 
    & \Longleftrightarrow &    
   \frac {\overrightarrow{|\B\Z|}}{\overrightarrow{|\Z\A|}}\cdot
   \frac {\overrightarrow{|\A\U|}}{\overrightarrow{|\U\D|}}\cdot
   \frac {\overrightarrow{|\D\V|}}{\overrightarrow{|\V\B|}}=-1&.
\end{align*}
Now these fractions contain {\it exactly} the same brackets in the numerator and denominator as the equations coming from the quadrangles on each side of the equality.
Thus we have a Menelaus proof. Hence Figure~\ref{fig: the big picture} visualizes the structure of all three proofs in one single picture.

\subsection{Pappus' Theorem}

In the last section we have seen, that binomial equations seem to be the ``building blocks'' for quad-proofs, in the sense that each quadrangle can be expressed by at most $3$ binomial equation. On the other hand we have seen, that binomial proofs and Ceva-Menelaus proofs are in a certain sense ``dual'' to each other, as one lives in the gaps of the other. We continue with one of the most fundamental projective incidence theorems: Pappus' Theorem.
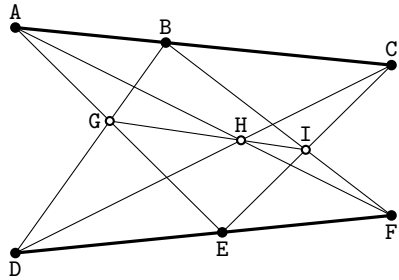
\begin{figure}[ht]
        \centering
        \begin{tikzpicture}[]
        \coordinate [label=above:{$\A$}] (A) at (0,3);
        \coordinate [label=above:{$\C$}] (C) at (5,2.5);
        \coordinate [label=below:{$\D$}] (D) at (0,0);
        \coordinate [label=below:{$\F$}] (F) at (5,0.5);

        \coordinate [label=above:{$\B$}] (B) at ($ (A)!0.4!(C) $);
        \coordinate [label=below:{$\E$}] (E) at ($ (D)!0.55!(F) $);
    
        \path [name path=A--C] (A) -- (C);
        \draw[very thick] (A) -- (C);
        \path [name path=D--F] (D) -- (F);
        \draw[very thick] (D) -- (F);
   
        \path [name path=A--E] (A) -- (E);
        \draw[] (A) -- (E);
        \path [name path=A--F] (A) -- (F);
        \draw[] (A) -- (F);
        \path [name path=B--D] (B) -- (D);
        \draw[] (B) -- (D);
        \path [name path=B--F] (B) -- (F);
        \draw[] (B) -- (F);
        \path [name path=C--D] (C) -- (D);
        \draw[] (C) -- (D);
        \path [name path=C--E] (C) -- (E);
        \draw[] (C) -- (E);

        \path [name intersections={of=A--E and B--D,by={[label=left:$\G$]G}}];
        \path [name intersections={of=A--F and C--D,by={[label=above:$\H$]H}}];
        \path [name intersections={of=B--F and C--E,by={[label=above:$\I$]I}}];

        \path [name path=G--I] (G) -- (I);
        \draw[] (G) -- (I);

        \foreach \point in {A,B,C,D,E,F} {
            \fill [black,opacity=1] (\point) circle (2pt);
        }

        \foreach \point in {G,H,I} {
            \path[fill=white,draw=black, thick] (\point) circle[radius=1.5pt];
        }
    \end{tikzpicture}
    \caption{The theorem of Pappus.}
    \label{fig: Drawing Pappus}
\end{figure}

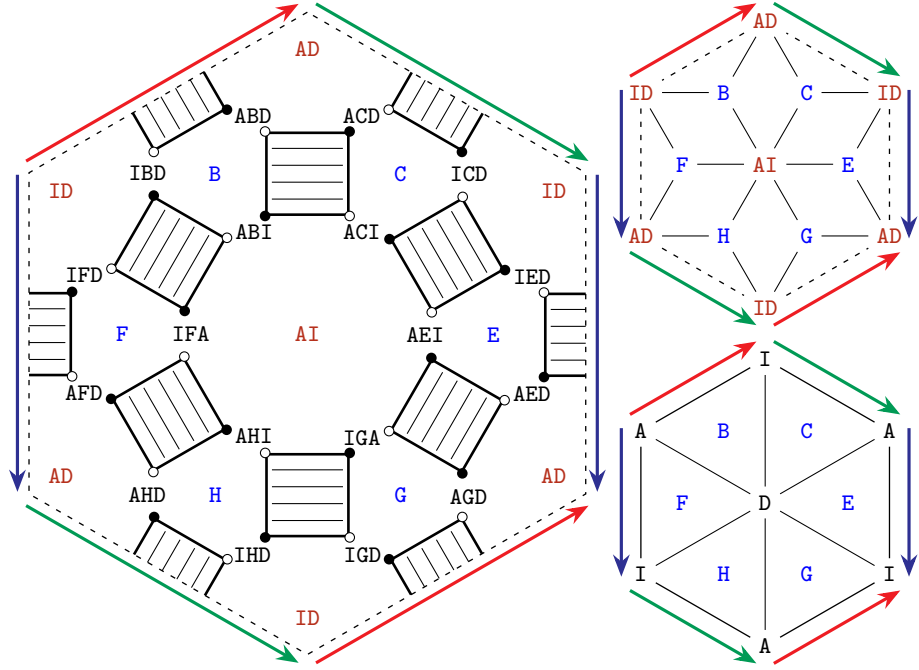
\begin{figure}[ht]
    \centering
\begin{minipage}[c]{0.64\linewidth}
    \def\factor{0.77}
    \begin{tikzpicture}[>=Stealth]
        \node[regular polygon, regular polygon sides=6, minimum width=0.77/0.65*7.2cm, rotate = 360/12] (hex) at (0,-\factor-\factor*1.73205-0.02) {};
        \node[regular polygon, regular polygon sides=6, minimum width=0.77/0.65*7.5cm, rotate = 360/12] (h) at (0,-\factor-\factor*1.73205-0.02) {};

        \draw[dash pattern=on 2pt off 3pt, line width=0.5pt] (hex.corner 1) -- (hex.corner 2) -- (hex.corner 3) -- (hex.corner 4) -- (hex.corner 5) -- (hex.corner 6) -- (hex.corner 1);
        
        \node (hex1) at (h.corner 1) {};
        \node (hex2) at (h.corner 2) {};
        \node (hex3) at (h.corner 3) {};
        \node (hex4) at (h.corner 4) {};
        \node (hex5) at (h.corner 5) {};
        \node (hex6) at (h.corner 6) {};

        \node (AI) at (0,-\factor-\factor*1.73205-0.02) {\color{BrickRed}$\A\I$};
        \node (AD1) at ($(AI)!0.85!(hex1)$) {\color{BrickRed}$\A\D$};
        \node (AD2) at ($(AI)!0.85!(hex3)$) {\color{BrickRed}$\A\D$};
        \node (AD3) at ($(AI)!0.85!(hex5)$) {\color{BrickRed}$\A\D$};
        \node (ID1) at ($(AI)!0.85!(hex2)$) {\color{BrickRed}$\I\D$};
        \node (ID2) at ($(AI)!0.85!(hex4)$) {\color{BrickRed}$\I\D$};
        \node (ID3) at ($(AI)!0.85!(hex6)$) {\color{BrickRed}$\I\D$};

        \draw[->, very thick, color=Red] (hex2) -- (hex1);
        \draw[->, very thick, color=Red] (hex4) -- (hex5);
        \draw[->, very thick, color=ForestGreen] (hex1) -- (hex6);
        \draw[->, very thick, color=ForestGreen] (hex3) -- (hex4);
        \draw[->, very thick, color=Blue] (hex6) -- (hex5);
        \draw[->, very thick, color=Blue] (hex2) -- (hex3);

        \clip (hex.corner 1) -- (hex.corner 2) -- (hex.corner 3) --
        (hex.corner 4) -- (hex.corner 5) -- (hex.corner 6) -- cycle;

        \node[regular polygon,
                draw,
                regular polygon sides = 4, line width=1pt,
                minimum size = \factor*2cm] (p) at (0,0) {};
        \node[regular polygon,
                regular polygon sides = 4,
                minimum size = \factor*2*sqrt(2)*1cm] (q) at (0,0) {};

        \draw[-] ($ ($ (p.corner 2)!0.2!(p.corner 3) $)!0.08!($ (p.corner 1)!0.2!(p.corner 4) $)$) -- ($ ($ (p.corner 2)!0.2!(p.corner 3) $)!0.92!($ (p.corner 1)!0.2!(p.corner 4) $) $);
        \draw[-] ($ ($ (p.corner 2)!0.4!(p.corner 3) $)!0.08!($ (p.corner 1)!0.4!(p.corner 4) $)$) -- ($ ($ (p.corner 2)!0.4!(p.corner 3) $)!0.92!($ (p.corner 1)!0.4!(p.corner 4) $) $);
        \draw[-] ($ ($ (p.corner 2)!0.6!(p.corner 3) $)!0.08!($ (p.corner 1)!0.6!(p.corner 4) $)$) -- ($ ($ (p.corner 2)!0.6!(p.corner 3) $)!0.92!($ (p.corner 1)!0.6!(p.corner 4) $) $);
        \draw[-] ($ ($ (p.corner 2)!0.8!(p.corner 3) $)!0.08!($ (p.corner 1)!0.8!(p.corner 4) $)$) -- ($ ($ (p.corner 2)!0.8!(p.corner 3) $)!0.92!($ (p.corner 1)!0.8!(p.corner 4) $) $);

        \foreach \point in {p.corner 1, p.corner 3} {
                \fill [black,opacity=1] (\point) circle (2pt);
            }
    
        \foreach \point in {p.corner 2, p.corner 4} {
                \path[fill=white,draw=black] (\point) circle[radius=1.8pt];
        }

        \node (q1) at (q.corner 2) {$\ \A\B\D$};
        \node (q2) at (q.corner 3) {};
        \node (q3) at (q.corner 4) {$\A\C\I\ $};
        \node (q4) at (q.corner 1) {};
        \node (tm) at ($ (q.center)!1+sqrt(3)/3!($ (q1)!0.5!(q2) $) $) {\color{blue}$\B$}; 
        \node[regular polygon,
                rotate=90,
                regular polygon sides = 3,
                minimum size = \factor*sqrt(3)/3*4cm] (t) at (tm) {};

        \node (qm) at ($ ($ (t.corner 3)!0.5!(t.corner 1) $) !\factor*-1cm! (t.center)$) {}; 
        \node[regular polygon,
                draw,line width=1pt,
                rotate=30,
                regular polygon sides = 4,line width=1pt,
                minimum size = \factor*2cm] (p) at (qm) {};
        \node[regular polygon,
                regular polygon sides = 4,line width=1pt,
                rotate=30,
                minimum size = \factor*2*sqrt(2)*1cm] (q) at (qm) {};   
        \foreach \point in {p.corner 2, p.corner 4} {
                \fill [black,opacity=1] (\point) circle (2pt);
            }
    
        \foreach \point in {p.corner 3, p.corner 1} {
                \path[fill=white,draw=black] (\point) circle[radius=1.8pt];
        }   
        \draw[-] ($ ($ (p.corner 1)!0.2!(p.corner 2) $)!0.08!($ (p.corner 4)!0.2!(p.corner 3) $)$) -- ($ ($ (p.corner 1)!0.2!(p.corner 2) $)!0.92!($ (p.corner 4)!0.2!(p.corner 3) $) $);
        \draw[-] ($ ($ (p.corner 1)!0.4!(p.corner 2) $)!0.08!($ (p.corner 4)!0.4!(p.corner 3) $)$) -- ($ ($ (p.corner 1)!0.4!(p.corner 2) $)!0.92!($ (p.corner 4)!0.4!(p.corner 3) $) $);
        \draw[-] ($ ($ (p.corner 1)!0.6!(p.corner 2) $)!0.08!($ (p.corner 4)!0.6!(p.corner 3) $)$) -- ($ ($ (p.corner 1)!0.6!(p.corner 2) $)!0.92!($ (p.corner 4)!0.6!(p.corner 3) $) $);
        \draw[-] ($ ($ (p.corner 1)!0.8!(p.corner 2) $)!0.08!($ (p.corner 4)!0.8!(p.corner 3) $)$) -- ($ ($ (p.corner 1)!0.8!(p.corner 2) $)!0.92!($ (p.corner 4)!0.8!(p.corner 3) $) $);
        \node (qm) at ($ ($ (t.corner 1)!0.5!(t.corner 2) $) !\factor*-1cm! (t.center)$) {}; 
        \node[regular polygon,
                draw,line width=1pt,
                rotate=60,
                regular polygon sides = 4,line width=1pt,
                minimum size = \factor*2cm] (p) at (qm) {};
        \node[regular polygon,
                regular polygon sides = 4,line width=1pt,
                rotate=60,
                minimum size = \factor*2*sqrt(2)*1cm] (q) at (qm) {};
        \draw[-] ($ ($ (p.corner 2)!0.2!(p.corner 3) $)!0.08!($ (p.corner 1)!0.2!(p.corner 4) $)$) -- ($ ($ (p.corner 2)!0.2!(p.corner 3) $)!0.92!($ (p.corner 1)!0.2!(p.corner 4) $) $);
        \draw[-] ($ ($ (p.corner 2)!0.4!(p.corner 3) $)!0.08!($ (p.corner 1)!0.4!(p.corner 4) $)$) -- ($ ($ (p.corner 2)!0.4!(p.corner 3) $)!0.92!($ (p.corner 1)!0.4!(p.corner 4) $) $);
        \draw[-] ($ ($ (p.corner 2)!0.6!(p.corner 3) $)!0.08!($ (p.corner 1)!0.6!(p.corner 4) $)$) -- ($ ($ (p.corner 2)!0.6!(p.corner 3) $)!0.92!($ (p.corner 1)!0.6!(p.corner 4) $) $);
        \draw[-] ($ ($ (p.corner 2)!0.8!(p.corner 3) $)!0.08!($ (p.corner 1)!0.8!(p.corner 4) $)$) -- ($ ($ (p.corner 2)!0.8!(p.corner 3) $)!0.92!($ (p.corner 1)!0.8!(p.corner 4) $) $);

        \foreach \point in {p.corner 1, p.corner 3} {
                \fill [black,opacity=1] (\point) circle (2pt);
            }
    
        \foreach \point in {p.corner 2, p.corner 4} {
                \path[fill=white,draw=black] (\point) circle[radius=1.8pt];
        }
        \node (q1) at (q.corner 2) {$\I\F\D\ $};
        \node (q2) at (q.corner 3) {$\I\F\A$};
        \node (q3) at (q.corner 4) {$\ \A\B\I$};
        \node (q4) at (q.corner 1) {$\I\B\D$};

        \node (tm) at ($ (q.center)!1+sqrt(3)/3!($ (q1)!0.5!(q2) $) $) {\color{blue}$\F$}; 
        \node[regular polygon,
                rotate=150,line width=1pt,
                regular polygon sides = 3,
                minimum size = \factor*sqrt(3)/3*4cm] (t) at (tm) {};
        
        \node (qm) at ($ ($ (t.corner 3)!0.5!(t.corner 1) $) !\factor*-1cm! (t.center)$) {}; 
        \node[regular polygon,
                draw,
                rotate=90,
                regular polygon sides = 4,line width=1pt,
                minimum size = \factor*2cm] (p) at (qm) {};
       \draw[-] ($ ($ (p.corner 1)!0.2!(p.corner 2) $)!0.08!($ (p.corner 4)!0.2!(p.corner 3) $)$) -- ($ ($ (p.corner 1)!0.2!(p.corner 2) $)!0.92!($ (p.corner 4)!0.2!(p.corner 3) $) $);
        \draw[-] ($ ($ (p.corner 1)!0.4!(p.corner 2) $)!0.08!($ (p.corner 4)!0.4!(p.corner 3) $)$) -- ($ ($ (p.corner 1)!0.4!(p.corner 2) $)!0.92!($ (p.corner 4)!0.4!(p.corner 3) $) $);
        \draw[-] ($ ($ (p.corner 1)!0.6!(p.corner 2) $)!0.08!($ (p.corner 4)!0.6!(p.corner 3) $)$) -- ($ ($ (p.corner 1)!0.6!(p.corner 2) $)!0.92!($ (p.corner 4)!0.6!(p.corner 3) $) $);
        \draw[-] ($ ($ (p.corner 1)!0.8!(p.corner 2) $)!0.08!($ (p.corner 4)!0.8!(p.corner 3) $)$) -- ($ ($ (p.corner 1)!0.8!(p.corner 2) $)!0.92!($ (p.corner 4)!0.8!(p.corner 3) $) $);
        \node[regular polygon,
                regular polygon sides = 4,line width=1pt,
                rotate=90,
                minimum size = \factor*2*sqrt(2)*1cm] (q) at (qm) {};   
        \foreach \point in {p.corner 2, p.corner 4} {
                \fill [black,opacity=1] (\point) circle (2pt);
            }
    
        \foreach \point in {p.corner 3, p.corner 1} {
                \path[fill=white,draw=black] (\point) circle[radius=1.8pt];
        }   

        \node (qm) at ($ ($ (t.corner 1)!0.5!(t.corner 2) $) !\factor*-1cm! (t.center)$) {}; 
        \node[regular polygon,
                draw,line width=1pt,
                rotate=120,
                regular polygon sides = 4,line width=1pt,
                minimum size = \factor*2cm] (p) at (qm) {};
        \draw[-] ($ ($ (p.corner 2)!0.2!(p.corner 3) $)!0.08!($ (p.corner 1)!0.2!(p.corner 4) $)$) -- ($ ($ (p.corner 2)!0.2!(p.corner 3) $)!0.92!($ (p.corner 1)!0.2!(p.corner 4) $) $);
        \draw[-] ($ ($ (p.corner 2)!0.4!(p.corner 3) $)!0.08!($ (p.corner 1)!0.4!(p.corner 4) $)$) -- ($ ($ (p.corner 2)!0.4!(p.corner 3) $)!0.92!($ (p.corner 1)!0.4!(p.corner 4) $) $);
        \draw[-] ($ ($ (p.corner 2)!0.6!(p.corner 3) $)!0.08!($ (p.corner 1)!0.6!(p.corner 4) $)$) -- ($ ($ (p.corner 2)!0.6!(p.corner 3) $)!0.92!($ (p.corner 1)!0.6!(p.corner 4) $) $);
        \draw[-] ($ ($ (p.corner 2)!0.8!(p.corner 3) $)!0.08!($ (p.corner 1)!0.8!(p.corner 4) $)$) -- ($ ($ (p.corner 2)!0.8!(p.corner 3) $)!0.92!($ (p.corner 1)!0.8!(p.corner 4) $) $);
        \node[regular polygon,
                regular polygon sides = 4,line width=1pt,
                rotate=120,
                minimum size = \factor*2*sqrt(2)*1cm] (q) at (qm) {};
        \foreach \point in {p.corner 1, p.corner 3} {
                \fill [black,opacity=1] (\point) circle (2pt);
            }
    
        \foreach \point in {p.corner 2, p.corner 4} {
                \path[fill=white,draw=black] (\point) circle[radius=1.8pt];
        }
        \node (q1) at (q.corner 2) {$\A\H\D$};
        \node (q2) at (q.corner 3) {$\ \A\H\I$};
        \node (q3) at (q.corner 4) {};
        \node (q4) at (q.corner 1) {$\A\F\D\ $};

        \node (tm) at ($ (q.center)!1+sqrt(3)/3!($ (q1)!0.5!(q2) $) $) {\color{blue}$\H$}; 
        \node[regular polygon,
                rotate=210,
                regular polygon sides = 3,
                minimum size = \factor*sqrt(3)/3*4cm] (t) at (tm) {};

        \node (qm) at ($ ($ (t.corner 3)!0.5!(t.corner 1) $) !\factor*-1cm! (t.center)$) {}; 
        \node[regular polygon,
                draw,line width=1pt,
                rotate=150,
                regular polygon sides = 4,line width=1pt,
                minimum size = \factor*2cm] (p) at (qm) {};
        \node[regular polygon,
                regular polygon sides = 4,line width=1pt,
                rotate=150,
                minimum size = \factor*2*sqrt(2)*1cm] (q) at (qm) {}; 
        \draw[-] ($ ($ (p.corner 1)!0.2!(p.corner 2) $)!0.08!($ (p.corner 4)!0.2!(p.corner 3) $)$) -- ($ ($ (p.corner 1)!0.2!(p.corner 2) $)!0.92!($ (p.corner 4)!0.2!(p.corner 3) $) $);
        \draw[-] ($ ($ (p.corner 1)!0.4!(p.corner 2) $)!0.08!($ (p.corner 4)!0.4!(p.corner 3) $)$) -- ($ ($ (p.corner 1)!0.4!(p.corner 2) $)!0.92!($ (p.corner 4)!0.4!(p.corner 3) $) $);
        \draw[-] ($ ($ (p.corner 1)!0.6!(p.corner 2) $)!0.08!($ (p.corner 4)!0.6!(p.corner 3) $)$) -- ($ ($ (p.corner 1)!0.6!(p.corner 2) $)!0.92!($ (p.corner 4)!0.6!(p.corner 3) $) $);
        \draw[-] ($ ($ (p.corner 1)!0.8!(p.corner 2) $)!0.08!($ (p.corner 4)!0.8!(p.corner 3) $)$) -- ($ ($ (p.corner 1)!0.8!(p.corner 2) $)!0.92!($ (p.corner 4)!0.8!(p.corner 3) $) $);
        \foreach \point in {p.corner 2, p.corner 4} {
                \fill [black,opacity=1] (\point) circle (2pt);
            }
    
        \foreach \point in {p.corner 3, p.corner 1} {
                \path[fill=white,draw=black] (\point) circle[radius=1.8pt];
        } 
        
        \node (qm) at ($ ($ (t.corner 1)!0.5!(t.corner 2) $) !\factor*-1cm! (t.center)$) {}; 
        \node[regular polygon,
                draw,line width=1pt,
                rotate=180,
                regular polygon sides = 4,line width=1pt,
                minimum size = \factor*2cm] (p) at (qm) {};
        \draw[-] ($ ($ (p.corner 2)!0.2!(p.corner 3) $)!0.08!($ (p.corner 1)!0.2!(p.corner 4) $)$) -- ($ ($ (p.corner 2)!0.2!(p.corner 3) $)!0.92!($ (p.corner 1)!0.2!(p.corner 4) $) $);
        \draw[-] ($ ($ (p.corner 2)!0.4!(p.corner 3) $)!0.08!($ (p.corner 1)!0.4!(p.corner 4) $)$) -- ($ ($ (p.corner 2)!0.4!(p.corner 3) $)!0.92!($ (p.corner 1)!0.4!(p.corner 4) $) $);
        \draw[-] ($ ($ (p.corner 2)!0.6!(p.corner 3) $)!0.08!($ (p.corner 1)!0.6!(p.corner 4) $)$) -- ($ ($ (p.corner 2)!0.6!(p.corner 3) $)!0.92!($ (p.corner 1)!0.6!(p.corner 4) $) $);
        \draw[-] ($ ($ (p.corner 2)!0.8!(p.corner 3) $)!0.08!($ (p.corner 1)!0.8!(p.corner 4) $)$) -- ($ ($ (p.corner 2)!0.8!(p.corner 3) $)!0.92!($ (p.corner 1)!0.8!(p.corner 4) $) $);
        \node[regular polygon,
                regular polygon sides = 4,line width=1pt,
                rotate=180,
                minimum size = \factor*2*sqrt(2)*1cm] (q) at (qm) {};
        \foreach \point in {p.corner 1, p.corner 3} {
                \fill [black,opacity=1] (\point) circle (2pt);
            }
    
        \foreach \point in {p.corner 2, p.corner 4} {
                \path[fill=white,draw=black] (\point) circle[radius=1.8pt];
        }
        \node (q1) at (q.corner 2) {$\I\G\D\ $};
        \node (q2) at (q.corner 3) {};
        \node (q3) at (q.corner 4) {};
        \node (q4) at (q.corner 1) {$\ \I\H\D$};

        \node (tm) at ($ (q.center)!1+sqrt(3)/3!($ (q1)!0.5!(q2) $) $) {\color{blue}$\G$}; 
        \node[regular polygon,
                rotate=270,
                regular polygon sides = 3,
                minimum size = \factor*sqrt(3)/3*4cm] (t) at (tm) {};

        \node (qm) at ($ ($ (t.corner 3)!0.5!(t.corner 1) $) !\factor*-1cm! (t.center)$) {}; 
        \node[regular polygon,line width=1pt,
                draw,
                rotate=210,
                regular polygon sides = 4,line width=1pt,
                minimum size = \factor*2cm] (p) at (qm) {};
        \draw[-] ($ ($ (p.corner 1)!0.2!(p.corner 2) $)!0.08!($ (p.corner 4)!0.2!(p.corner 3) $)$) -- ($ ($ (p.corner 1)!0.2!(p.corner 2) $)!0.92!($ (p.corner 4)!0.2!(p.corner 3) $) $);
        \draw[-] ($ ($ (p.corner 1)!0.4!(p.corner 2) $)!0.08!($ (p.corner 4)!0.4!(p.corner 3) $)$) -- ($ ($ (p.corner 1)!0.4!(p.corner 2) $)!0.92!($ (p.corner 4)!0.4!(p.corner 3) $) $);
        \draw[-] ($ ($ (p.corner 1)!0.6!(p.corner 2) $)!0.08!($ (p.corner 4)!0.6!(p.corner 3) $)$) -- ($ ($ (p.corner 1)!0.6!(p.corner 2) $)!0.92!($ (p.corner 4)!0.6!(p.corner 3) $) $);
        \draw[-] ($ ($ (p.corner 1)!0.8!(p.corner 2) $)!0.08!($ (p.corner 4)!0.8!(p.corner 3) $)$) -- ($ ($ (p.corner 1)!0.8!(p.corner 2) $)!0.92!($ (p.corner 4)!0.8!(p.corner 3) $) $);
        \node[regular polygon,
                regular polygon sides = 4,line width=1pt,
                rotate=210,
                minimum size = \factor*2*sqrt(2)*1cm] (q) at (qm) {}; 
        \foreach \point in {p.corner 2, p.corner 4} {
                \fill [black,opacity=1] (\point) circle (2pt);
            }
    
        \foreach \point in {p.corner 3, p.corner 1} {
                \path[fill=white,draw=black] (\point) circle[radius=1.8pt];
        }

        \node (qm) at ($ ($ (t.corner 1)!0.5!(t.corner 2) $) !\factor*-1cm! (t.center)$) {}; 
        \node[regular polygon,
                draw,
                rotate=240,
                regular polygon sides = 4,line width=1pt,
                minimum size = \factor*2cm] (p) at (qm) {};
        \draw[-] ($ ($ (p.corner 2)!0.2!(p.corner 3) $)!0.08!($ (p.corner 1)!0.2!(p.corner 4) $)$) -- ($ ($ (p.corner 2)!0.2!(p.corner 3) $)!0.92!($ (p.corner 1)!0.2!(p.corner 4) $) $);
        \draw[-] ($ ($ (p.corner 2)!0.4!(p.corner 3) $)!0.08!($ (p.corner 1)!0.4!(p.corner 4) $)$) -- ($ ($ (p.corner 2)!0.4!(p.corner 3) $)!0.92!($ (p.corner 1)!0.4!(p.corner 4) $) $);
        \draw[-] ($ ($ (p.corner 2)!0.6!(p.corner 3) $)!0.08!($ (p.corner 1)!0.6!(p.corner 4) $)$) -- ($ ($ (p.corner 2)!0.6!(p.corner 3) $)!0.92!($ (p.corner 1)!0.6!(p.corner 4) $) $);
        \draw[-] ($ ($ (p.corner 2)!0.8!(p.corner 3) $)!0.08!($ (p.corner 1)!0.8!(p.corner 4) $)$) -- ($ ($ (p.corner 2)!0.8!(p.corner 3) $)!0.92!($ (p.corner 1)!0.8!(p.corner 4) $) $);
        \node[regular polygon,
                regular polygon sides = 4,line width=1pt,
                rotate=240,
                minimum size = \factor*2*sqrt(2)*1cm] (q) at (qm) {};
        \foreach \point in {p.corner 1, p.corner 3} {
                \fill [black,opacity=1] (\point) circle (2pt);
            }
    
        \foreach \point in {p.corner 2, p.corner 4} {
                \path[fill=white,draw=black] (\point) circle[radius=1.8pt];
        }
        \node (q1) at (q.corner 2) {$\ \A\E\D$};
        \node (q2) at (q.corner 3) {$\A\E\I$};
        \node (q3) at (q.corner 4) {$\I\G\A\ $};
        \node (q4) at (q.corner 1) {$\A\G\D$};

        \node (tm) at ($ (q.center)!1+sqrt(3)/3!($ (q1)!0.5!(q2) $) $) {\color{blue}$\E$}; 
        \node[regular polygon,
                rotate=330,line width=1pt,
                regular polygon sides = 3,
                minimum size = \factor*sqrt(3)/3*4cm] (t) at (tm) {};

        \node (qm) at ($ ($ (t.corner 3)!0.5!(t.corner 1) $) !\factor*-1cm! (t.center)$) {}; 
        \node[regular polygon,
                draw,
                rotate=270,
                regular polygon sides = 4,line width=1pt,
                minimum size = \factor*2cm] (p) at (qm) {};
        \draw[-] ($ ($ (p.corner 1)!0.2!(p.corner 2) $)!0.08!($ (p.corner 4)!0.2!(p.corner 3) $)$) -- ($ ($ (p.corner 1)!0.2!(p.corner 2) $)!0.92!($ (p.corner 4)!0.2!(p.corner 3) $) $);
        \draw[-] ($ ($ (p.corner 1)!0.4!(p.corner 2) $)!0.08!($ (p.corner 4)!0.4!(p.corner 3) $)$) -- ($ ($ (p.corner 1)!0.4!(p.corner 2) $)!0.92!($ (p.corner 4)!0.4!(p.corner 3) $) $);
        \draw[-] ($ ($ (p.corner 1)!0.6!(p.corner 2) $)!0.08!($ (p.corner 4)!0.6!(p.corner 3) $)$) -- ($ ($ (p.corner 1)!0.6!(p.corner 2) $)!0.92!($ (p.corner 4)!0.6!(p.corner 3) $) $);
        \draw[-] ($ ($ (p.corner 1)!0.8!(p.corner 2) $)!0.08!($ (p.corner 4)!0.8!(p.corner 3) $)$) -- ($ ($ (p.corner 1)!0.8!(p.corner 2) $)!0.92!($ (p.corner 4)!0.8!(p.corner 3) $) $);
        \node[regular polygon,
                regular polygon sides = 4,line width=1pt,
                rotate=270,
                minimum size = \factor*2*sqrt(2)*1cm] (q) at (qm) {}; 
        \foreach \point in {p.corner 2, p.corner 4} {
                \fill [black,opacity=1] (\point) circle (2pt);
            }
    
        \foreach \point in {p.corner 3, p.corner 1} {
                \path[fill=white,draw=black] (\point) circle[radius=1.8pt];
        }

        \node (qm) at ($ ($ (t.corner 1)!0.5!(t.corner 2) $) !\factor*-1cm! (t.center)$) {}; 
        \node[regular polygon,
                draw,
                rotate=300,
                regular polygon sides = 4,line width=1pt,
                minimum size = \factor*2cm] (p) at (qm) {};
        \draw[-] ($ ($ (p.corner 2)!0.2!(p.corner 3) $)!0.08!($ (p.corner 1)!0.2!(p.corner 4) $)$) -- ($ ($ (p.corner 2)!0.2!(p.corner 3) $)!0.92!($ (p.corner 1)!0.2!(p.corner 4) $) $);
        \draw[-] ($ ($ (p.corner 2)!0.4!(p.corner 3) $)!0.08!($ (p.corner 1)!0.4!(p.corner 4) $)$) -- ($ ($ (p.corner 2)!0.4!(p.corner 3) $)!0.92!($ (p.corner 1)!0.4!(p.corner 4) $) $);
        \draw[-] ($ ($ (p.corner 2)!0.6!(p.corner 3) $)!0.08!($ (p.corner 1)!0.6!(p.corner 4) $)$) -- ($ ($ (p.corner 2)!0.6!(p.corner 3) $)!0.92!($ (p.corner 1)!0.6!(p.corner 4) $) $);
        \draw[-] ($ ($ (p.corner 2)!0.8!(p.corner 3) $)!0.08!($ (p.corner 1)!0.8!(p.corner 4) $)$) -- ($ ($ (p.corner 2)!0.8!(p.corner 3) $)!0.92!($ (p.corner 1)!0.8!(p.corner 4) $) $);
        \node[regular polygon,
                regular polygon sides = 4,line width=1pt,
                rotate=300,
                minimum size = \factor*2*sqrt(2)*1cm] (q) at (qm) {};
        \foreach \point in {p.corner 1, p.corner 3} {
                \fill [black,opacity=1] (\point) circle (2pt);
            }
    
        \foreach \point in {p.corner 2, p.corner 4} {
                \path[fill=white,draw=black] (\point) circle[radius=1.8pt];
        }
        \node (q1) at (q.corner 2) {$\I\C\D$};
        \node (q2) at (q.corner 3) {};
        \node (q3) at (q.corner 4) {};
        \node (q4) at (q.corner 1) {$\ \I\E\D$};

        \node (tm) at ($ (q.center)!1+sqrt(3)/3!($ (q1)!0.5!(q2) $) $) {\color{blue}$\C$}; 
        \node[regular polygon,
                rotate=30,
                regular polygon sides = 3,
                minimum size = \factor*sqrt(3)/3*4cm] (t) at (tm) {};

        \node (qm) at ($ ($ (t.corner 3)!0.5!(t.corner 1) $) !\factor*-1cm! (t.center)$) {}; 
        \node[regular polygon,
                draw,
                rotate=330,
                regular polygon sides = 4,line width=1pt,
                minimum size = \factor*2cm] (p) at (qm) {};
        \draw[-] ($ ($ (p.corner 1)!0.2!(p.corner 2) $)!0.08!($ (p.corner 4)!0.2!(p.corner 3) $)$) -- ($ ($ (p.corner 1)!0.2!(p.corner 2) $)!0.92!($ (p.corner 4)!0.2!(p.corner 3) $) $);
        \draw[-] ($ ($ (p.corner 1)!0.4!(p.corner 2) $)!0.08!($ (p.corner 4)!0.4!(p.corner 3) $)$) -- ($ ($ (p.corner 1)!0.4!(p.corner 2) $)!0.92!($ (p.corner 4)!0.4!(p.corner 3) $) $);
        \draw[-] ($ ($ (p.corner 1)!0.6!(p.corner 2) $)!0.08!($ (p.corner 4)!0.6!(p.corner 3) $)$) -- ($ ($ (p.corner 1)!0.6!(p.corner 2) $)!0.92!($ (p.corner 4)!0.6!(p.corner 3) $) $);
        \draw[-] ($ ($ (p.corner 1)!0.8!(p.corner 2) $)!0.08!($ (p.corner 4)!0.8!(p.corner 3) $)$) -- ($ ($ (p.corner 1)!0.8!(p.corner 2) $)!0.92!($ (p.corner 4)!0.8!(p.corner 3) $) $);
        \node[regular polygon,
                regular polygon sides = 4,line width=1pt,
                rotate=330,
                minimum size = \factor*2*sqrt(2)*1cm] (q) at (qm) {}; 
        \foreach \point in {p.corner 2, p.corner 4} {
                \fill [black,opacity=1] (\point) circle (2pt);
            }
    
        \foreach \point in {p.corner 3, p.corner 1} {
                \path[fill=white,draw=black] (\point) circle[radius=1.8pt];
        }
        \node (q4) at (q.corner 3) {$\A\C\D\ $};

    \end{tikzpicture}
\end{minipage}
\begin{minipage}[c]{0.35\linewidth}
    \centering
    \def\ngon{6}
        \begin{tikzpicture}[>=Stealth]
            \node[regular polygon,regular polygon sides=\ngon,rotate=-30,minimum size=3.8cm] (r) at (0,0) {};
            \node (r1) at (r.corner 1) {\color{BrickRed}$\I\D$};
            \node (r2) at (r.corner 2) {\color{BrickRed}$\A\D$};
            \node (r3) at (r.corner 3) {\color{BrickRed}$\I\D$};
            \node (r4) at (r.corner 4) {\color{BrickRed}$\A\D$};
            \node (r5) at (r.corner 5) {\color{BrickRed}$\I\D$};
            \node (r6) at (r.corner 6) {\color{BrickRed}$\A\D$};
    
            \draw[dash pattern=on 2pt off 3pt, line width=0.5pt] (r1) -- (r2);
            \draw[dash pattern=on 2pt off 3pt, line width=0.5pt] (r4) -- (r5);
            \draw[dash pattern=on 2pt off 3pt, line width=0.5pt] (r6) -- (r1);
            \draw[dash pattern=on 2pt off 3pt, line width=0.5pt] (r3) -- (r4);
            \draw[dash pattern=on 2pt off 3pt, line width=0.5pt] (r2) -- (r3);
            \draw[dash pattern=on 2pt off 3pt, line width=0.5pt] (r5) -- (r6);
    
            \node (A1) at (barycentric cs:r2=1,r4=1,r6=1) {\color{BrickRed}$\A\I$};
    
            \node (l1) at (barycentric cs:r1=1,r2=1,A1=1) {\color{blue}$\C$};
            \node (l2) at (barycentric cs:r2=1,r3=1,A1=1) {\color{blue}$\B$};
            \node (l3) at (barycentric cs:r3=1,r4=1,A1=1) {\color{blue}$\F$};
            \node (l4) at (barycentric cs:r4=1,r5=1,A1=1) {\color{blue}$\H$};
            \node (l5) at (barycentric cs:r5=1,r6=1,A1=1) {\color{blue}$\G$};
            \node (l6) at (barycentric cs:r6=1,r1=1,A1=1) {\color{blue}$\E$}; 
    
            \draw[-] (r1) -- (l1);
            \draw[-] (r2) -- (l2);
            \draw[-] (r3) -- (l3);
            \draw[-] (r4) -- (l4);
            \draw[-] (r5) -- (l5);
            \draw[-] (r6) -- (l6);
            
            \draw[-] (r2) -- (l1);
            \draw[-] (r3) -- (l2);
            \draw[-] (r4) -- (l3);
            \draw[-] (r5) -- (l4);
            \draw[-] (r6) -- (l5);
            \draw[-] (r1) -- (l6);
    
            \draw[-] (A1) -- (l1);
            \draw[-] (A1) -- (l2);
            \draw[-] (A1) -- (l3);
            \draw[-] (A1) -- (l4);
            \draw[-] (A1) -- (l5);
            \draw[-] (A1) -- (l6);

            \node[regular polygon, regular polygon sides=6, minimum width=4.4cm, rotate = 360/12] (h) at (A1) {};

            \node (hex1) at (h.corner 1) {\ };
            \node (hex2) at (h.corner 2) {\ };
            \node (hex3) at (h.corner 3) {\ };
            \node (hex4) at (h.corner 4) {\ };
            \node (hex5) at (h.corner 5) {\ };
            \node (hex6) at (h.corner 6) {\ };
            \draw[->, very thick, color=Red] (hex2) -- (hex1);
            \draw[->, very thick, color=Red] (hex4) -- (hex5);
            \draw[->, very thick, color=ForestGreen] (hex1) -- (hex6);
            \draw[->, very thick, color=ForestGreen] (hex3) -- (hex4);
            \draw[->, very thick, color=Blue] (hex6) -- (hex5);
            \draw[->, very thick, color=Blue] (hex2) -- (hex3);

            \node[regular polygon,regular polygon sides=\ngon,rotate=-30,minimum size=3.8cm] (r) at (0,-4.5) {};
            \node (r1) at (r.corner 1) {$\A$};
            \node (r2) at (r.corner 2) {$\I$};
            \node (r3) at (r.corner 3) {$\A$};
            \node (r4) at (r.corner 4) {$\I$};
            \node (r5) at (r.corner 5) {$\A$};
            \node (r6) at (r.corner 6) {$\I$};
    
            \draw[ line width=0.5pt] (r1) -- (r2);
            \draw[ line width=0.5pt] (r4) -- (r5);
            \draw[ line width=0.5pt] (r6) -- (r1);
            \draw[ line width=0.5pt] (r3) -- (r4);
            \draw[line width=0.5pt] (r2) -- (r3);
            \draw[line width=0.5pt] (r5) -- (r6);
    
            \node (A1) at (barycentric cs:r2=1,r4=1,r6=1) {$\D$};
    
            \node (l1) at (barycentric cs:r1=1,r2=1,A1=1) {\color{blue}$\C$};
            \node (l2) at (barycentric cs:r2=1,r3=1,A1=1) {\color{blue}$\B$};
            \node (l3) at (barycentric cs:r3=1,r4=1,A1=1) {\color{blue}$\F$};
            \node (l4) at (barycentric cs:r4=1,r5=1,A1=1) {\color{blue}$\H$};
            \node (l5) at (barycentric cs:r5=1,r6=1,A1=1) {\color{blue}$\G$};
            \node (l6) at (barycentric cs:r6=1,r1=1,A1=1) {\color{blue}$\E$}; 
            \draw[-] (r1) -- (A1);
            \draw[-] (r2) -- (A1);
            \draw[-] (r3) -- (A1);
            \draw[-] (r4) -- (A1);
            \draw[-] (r5) -- (A1);
            \draw[-] (r6) -- (A1);

            \node[regular polygon, regular polygon sides=6, minimum width=4.4cm, rotate = 360/12] (h) at (A1) {};
        
        \node (hex1) at (h.corner 1) {};
        \node (hex2) at (h.corner 2) {};
        \node (hex3) at (h.corner 3) {};
        \node (hex4) at (h.corner 4) {};
        \node (hex5) at (h.corner 5) {};
        \node (hex6) at (h.corner 6) {};
        \draw[->, very thick, color=Red] (hex2) -- (hex1);
        \draw[->, very thick, color=Red] (hex4) -- (hex5);
        \draw[->, very thick, color=ForestGreen] (hex1) -- (hex6);
        \draw[->, very thick, color=ForestGreen] (hex3) -- (hex4);
        \draw[->, very thick, color=Blue] (hex6) -- (hex5);
        \draw[->, very thick, color=Blue] (hex2) -- (hex3);
        
        \end{tikzpicture}
\end{minipage}
\label{fig: binomial proof}
\caption{The proofs of Pappos's theorem visualized.}
\end{figure}
Again, both~\cite{FP} and~\cite{RG1} give tiling-proofs for it and we analyze how they relate to each other and reveal some additional structure, which was ``hidden'' in the proofs of Desargues' theorem. A similar analysis but without the quad-proofs
was already performed in \cite{ARG}.
We will start from the binomial proof:
\begin{align*}
        \A\B\C \quad\text{ collinear} &\quad\Longleftrightarrow\quad[\A,\B,\D][\A,\C,\I]=[\A,\B,\I][\A,\C,\D] \\
        \A\H\F \quad\text{ collinear} &\quad\Longleftrightarrow\quad[\A,\H,\D][\A,\F,\I]=[\A,\H,\I][\A,\F,\D] \\
        \A\E\G \quad\text{ collinear} &\quad\Longleftrightarrow\quad[\A,\E,\D][\A,\G,\I]=[\A,\E,\I][\A,\G,\D] \\   
        \D\C\H \quad\text{ collinear} &\quad\Longleftrightarrow\quad[\D,\C,\A][\D,\H,\I]=[\D,\C,\I][\D,\H,\A] \\
        \D\F\E \quad\text{ collinear} &\quad\Longleftrightarrow\quad[\D,\F,\A][\D,\E,\I]=[\D,\F,\I][\D,\E,\A] \\
        \D\G\B \quad\text{ collinear} &\quad\Longleftrightarrow\quad[\D,\G,\A][\D,\B,\I]=[\D,\G,\I][\D,\B,\A] \\
        \I\G\H \quad\text{ collinear} &\quad\Longleftrightarrow\quad[\I,\G,\D][\I,\H,\A]=[\I,\G,\A][\I,\H,\D] \\
        \I\B\F \quad\text{ collinear} &\quad\Longleftrightarrow\quad[\I,\B,\A][\I,\F,\D]=[\I,\B,\D][\I,\F,\A] \\
        \I\C\E \quad\text{ collinear} &\quad\Longleftrightarrow\quad[\I,\C,\D][\I,\E,\A]=[\I,\C,\A][\I,\E,\D]
\end{align*}

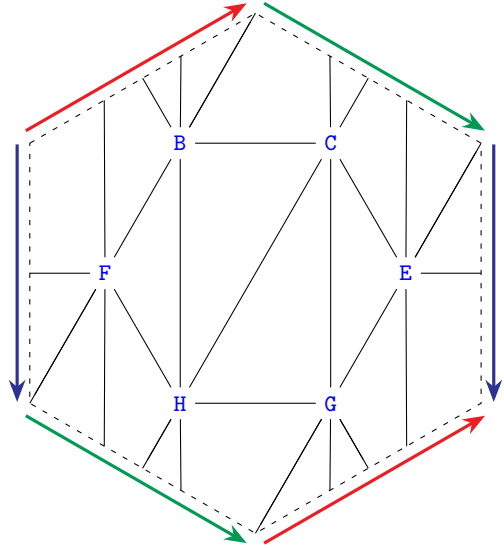
\begin{figure}[h]
\centering
    \begin{tikzpicture}[>=Stealth]
            \node (center) at (0,0) { };
            \node (F) at (2,0) {\color{blue}$\E$};
            \node (G) at (-2,0) {\color{blue}$\F$};
            \node (E) at (1,1.73) {\color{blue}$\C$};
            \node (I) at (1,-1.73) {\color{blue}$\G$};
            \node (H) at (-1,1.73) {\color{blue}$\B$};
            \node (D) at (-1,-1.73) {\color{blue}$\H$};

            \node (X) at ($ (F)!0.5!(E) $) {};
            \node (X2) at ($ (G)!0.5!(D) $) {};
            \node (Y) at ($ (G)!0.5!(H) $) {};
            \node (Y2) at ($ (F)!0.5!(I) $) {};
            \node (Z2) at ($ (E)!0.5!(H) $) {};
            \node (Z) at ($ (I)!0.5!(D) $) {};

            \node (A1) at ($ (center)!2!(X) $) {};
            \node (A2) at ($ (center)!2!(X2) $) {};
            \node (A3) at ($ (center)!2!(Y) $) {};
            \node (A4) at ($ (center)!2!(Y2) $) {};
            \node (A5) at ($ (center)!2!(Z) $) {};
            \node (A6) at ($ (center)!2!(Z2) $) {};

            \node (o) at (barycentric cs:G=0.7,D=1,H=1) {};
            \node (n) at (barycentric cs:F=0.7,I=1,E=1) {};
            \node (m) at ($ (center)!0.5!(H) $) {};
            \node (l) at ($ (center)!0.5!(I) $) {};

            \node (t) at ($ (A6)!0.5!(E) $) {};
            \node (u) at ($ (A5)!0.5!(D) $) {};

            \node (v) at ($ (A4)!0.5!(F) $) {};
            \node (help13) at ($ (A3)!0.8!(A6) $) {};
            \node (help14) at ($ (A5)!0.8!(A4) $) {};      
            \node (help15) at ($ (A4)!0.34!(A1) $) {};
            \node (help16) at ($ (A2)!0.34!(A3) $) {};
            
            \node (w) at ($ (A3)!0.5!(G) $) {};
            \node (help17) at ($ (A3)!0.8!(A6) $) {};
            \node (help18) at ($ (A5)!0.8!(A4) $) {};
            \node (help9) at ($ (A1)!0.2!(A6) $) {};
            \node (help10) at ($ (A5)!0.2!(A2) $) {};
            \node (help11) at ($ (A1)!0.8!(A6) $) {};
            \node (help12) at ($ (A5)!0.8!(A2) $) {};

            \node (p) at ($ (E)!0.37!(A1) $) {};

            \node (help1) at ($ (A1)!0.38!(A6) $) {};
            \node (help2) at ($ (A5)!0.38!(A2) $) {};

            \node (q) at ($ (D)!0.37!(A2) $) {};

            \node (help3) at ($ (A1)!0.61!(A6) $) {};
            \node (help4) at ($ (A5)!0.61!(A2) $) {};

            \node (r) at ($ (H)!0.37!(A3) $) {};

            \node (help5) at ($ (A3)!0.38!(A6) $) {};
            \node (help6) at ($ (A5)!0.38!(A4) $) {};

            \node (s) at ($ (I)!0.37!(A4) $) {};

            \node (help7) at ($ (A3)!0.61!(A6) $) {};

            \node (help8) at ($ (A5)!0.61!(A4) $) {};

            \node (B1) at ($ (G)!0.5!(D) $) {};
            \node (B2) at ($ (D)!0.5!(I) $) {};
            \node (B3) at ($ (I)!0.5!(F) $) {};
            \node (B4) at ($ (F)!0.5!(E) $) {};
            \node (B5) at ($ (E)!0.5!(H) $) {};
            \node (B6) at ($ (H)!0.5!(G) $) {};
            \node (B7) at ($ (H)!0.5!(D) $) {};
            \node (B8) at ($ (D)!0.5!(E) $) {};
            \node (B9) at ($ (E)!0.5!(I) $) {};
            \node (B15) at ($ (D)!0.5!(A2) $) {};
            \node (B16) at ($ (I)!0.5!(A4) $) {};
            \node (B17) at ($ (E)!0.5!(A1) $) {};
            \node (B18) at ($ (H)!0.5!(A3) $) {};
            \node (B12) at ($ (A5)!0.5!(A2) $) {};
            \node (B13) at ($ (A5)!0.5!(A4) $) {};
            \node (B14) at ($ (A1)!0.5!(A4) $) {};
            \node (B122) at ($ (A1)!0.5!(A6) $) {};
            \node (B132) at ($ (A3)!0.5!(A6) $) {};
            \node (B142) at ($ (A3)!0.5!(A2) $) {};

            \draw[dash pattern=on 2pt off 3pt] ($ (center)!2!(Y) $) -- ($ (center)!2!(Z2) $) -- ($ (center)!2!(X) $) -- ($ (center)!2!(Y2) $) -- ($ (center)!2!(Z) $) -- ($ (center)!2!(X2) $) -- ($ (center)!2!(Y) $);

            \draw[-] (H) -- (G);
            \draw[-] (D) -- (G);
            \draw[-] (H) -- (D);
            \draw[-] (D) -- (E);
            \draw[-] (H) -- (E);
            \draw[-] (D) -- (I);
            \draw[-] (E) -- (I);
            \draw[-] (E) -- (F);
            \draw[-] (F) -- (I);
            \draw[-] (F) -- ($ (A1)!0.5!(A4) $);
            \draw[-] (F) -- ($ (A1)!0!(A4) $);
            \draw[-] (G) -- ($ (A3)!0.5!(A2) $);
            \draw[-] (G) -- ($ (A3)!1!(A2) $);
            \draw[-] (I) -- (B13);
            \draw[-] (H) -- ($ (A3)!0.5!(A6) $);
            \draw[-] (H) -- ($ (A3)!1!(A6) $);
            \draw[-] (E) -- ($ (A1)!0.5!(A6) $);
            \draw[-] (D) -- (B12);
            \draw[-] (H) -- (A6);
            \draw[-] (A2) -- (G);
            \draw[-] (A1) -- (F);
            \draw[-] (A5) -- (I);
            \draw[-] (G) -- ($ (A3)!0.33!(A6) $);
            \draw[-] (I) -- ($ (A5)!0.33!(A4) $);
            \draw[-] (I) -- ($ (A5)!0.5!(A4) $);
            \draw[-] (I) -- ($ (A5)!0!(A4) $);
            \draw[-] (D) -- ($ (A5)!0.33!(A2) $);
            \draw[-] (D) -- ($ (A5)!0.5!(A2) $);
            \draw[-] (E) -- ($ (A1)!0.67!(A6) $);
            \draw[-] (G) -- ($ (A5)!0.67!(A2) $);
            \draw[-] (F) -- ($ (A1)!0.33!(A6) $);
            \draw[-] (H) -- ($ (A3)!0.67!(A6) $);
            \draw[-] (F) -- ($ (A5)!0.67!(A4) $);

            \node[regular polygon, regular polygon sides=6, minimum width=7.3cm, rotate = 360/12] (h) at ($ (E)!0.5!(D) $) {};

            \node (hex1) at (h.corner 1) {\ };
            \node (hex2) at (h.corner 2) {\ };
            \node (hex3) at (h.corner 3) {\ };
            \node (hex4) at (h.corner 4) {\ };
            \node (hex5) at (h.corner 5) {\ };
            \node (hex6) at (h.corner 6) {\ };
            \draw[->, very thick, color=Red] (hex2) -- (hex1);
            \draw[->, very thick, color=Red] (hex4) -- (hex5);
            \draw[->, very thick, color=ForestGreen] (hex1) -- (hex6);
            \draw[->, very thick, color=ForestGreen] (hex3) -- (hex4);
            \draw[->, very thick, color=Blue] (hex6) -- (hex5);
            \draw[->, very thick, color=Blue] (hex2) -- (hex3);

    \end{tikzpicture}
    \caption{Proof of Pappos's theorem using $12$ Menelaus triangles.}
    \label{fig: 12 mene}
\end{figure}

As before, any of the $9$ equations is implied by the other $8$. Now, in the same way as in Figure~\ref{fig: the big picture} we will depict the binomial equations as small rectangles, glued together whereever they share a bracket. While the proofs for Desargues' theorem lived on a sphere the resulting structure here lives on a torus which can easily be seen on the left in Figure~\ref{fig: binomial proof}. As we already saw in the last example, the holes correspond to points or lines, depending on how many points in the brackets are shared. For example, the top left triangle has the brackets $\A\B\D, \A\B\I$ and $\B\D\I$ as vertices, so all of them contain the point $\B$. On the other hand, all $6$ vertices of
the hexagon in the middle share the points $\A$ and $\I$, so this gives us the line spanned by them. This way we can immediately read of quadrilateral-tiling proof: We just need to create quadrangles with points and lines as vertices, as seen in Figure~\ref{fig: binomial proof} top right. Extracting a Ceva-Menelaus-Proof from this picture is just as easy, as the triangles whose share $1$ point correspond exactly to Ceva-triangles with the shared point as ``Ceva-point''. In particular, e.g. the top left triangle $\A\B\D \rightarrow \I\B\D \rightarrow \A\B\I$ corresponds to the Ceva-triangle $\I \rightarrow \A \rightarrow \D$ with the point $\B$ in the middle, yielding:
\begin{equation*}
    1 = \frac{[\A,\B,\D]}{[\I,\B,\D]}\frac{[\I,\B,\D]}{[\A,\B,\I]}\frac{[\A,\B,\I]}{[\A,\B,\D]},
\end{equation*}
which is equal to the length ratios along the edges of the triangle. Doing this for all $6$ triangles we end up with the pure Ceva-proof in the bottom right of Figure~\ref{fig: binomial proof}. In fact, these two proofs were already presented in \cite{FP} and \cite{RG1}, respectively and are now shown to be ``the same proof'' in the sense of coming from the same binomial cancellation pattern. 
\nnl
If we compare Figures \ref{fig: the big picture} and \ref{fig: binomial proof} again, we can see that in the first one we saw triangular holes where all brackets at the vertices shared $2$ points and in the second one we saw triangular holes where all brackets at the vertices vertices shared only $1$ point. We then saw that these holes can be interpreted as Menelaus-~respectively Ceva-triangles. We now focus on the non-triangular holes. They will lead to a  pure Menelaus-proof whose existence has been claimed in Section 4.
 If we take a close look at the $3$ hexagonal holes in \ref{fig: binomial proof} we can see, that all brackets at their vertices share the same $2$ points: $\A$ and $\I$, for the central one, $\A$ and $\D$ for the one in the top and $I$ and $\D$ for the one in the bottom. So actually, each of those can be triangulated using $4$ Menelaus-triangles resulting in another proof of Pappus's theorem using $12$ Menelaus triangles, which happens to be exactly one of Menelaus proofs we would get applying the translation procedure from Section 4.

This method in fact works for all binomial proofs! The proof in \cite{ARG} showing that Ceva-Menelaus-proofs and binomial proofs are equivalent, forms chains of binomial equations (the black rectangles in Figure \ref{fig: binomial proof}) and uses cycles of unmatched edges to fit Ceva- and Menelaus-triangles inside. In fact there are always $2$ possibilities to do this, 
yielding at least $2$  Ceva-Menelaus-proofs per binomial proof.

\section{Three three-dimensional examples}
One might wonder, how the different concepts generalize to higher dimensions.    
In~\cite{FP}
Fomin and  Pylyavskyy give several examples of quad-proofs that can be considered as proofs for 3-dimensional incidence theorems. Also in \cite{RG2} bi-quadratic final polynomial proofs are given for some 3D incidence theorems. We here exemplify how the Ceva/Menelaus Method can be interpreted in a 3-dimensional context. For that{\color{ForestGreen},} we first have to explore how the theorems of Ceva and Menelaus can be transformed to higher dimensions. We restrict ourselves to the 3-dimensional case.
\begin{figure}[h]
    \centering
    \begin{tikzpicture}
        \draw (0, 0) node[inner sep=0] {\includegraphics[width=10cm]{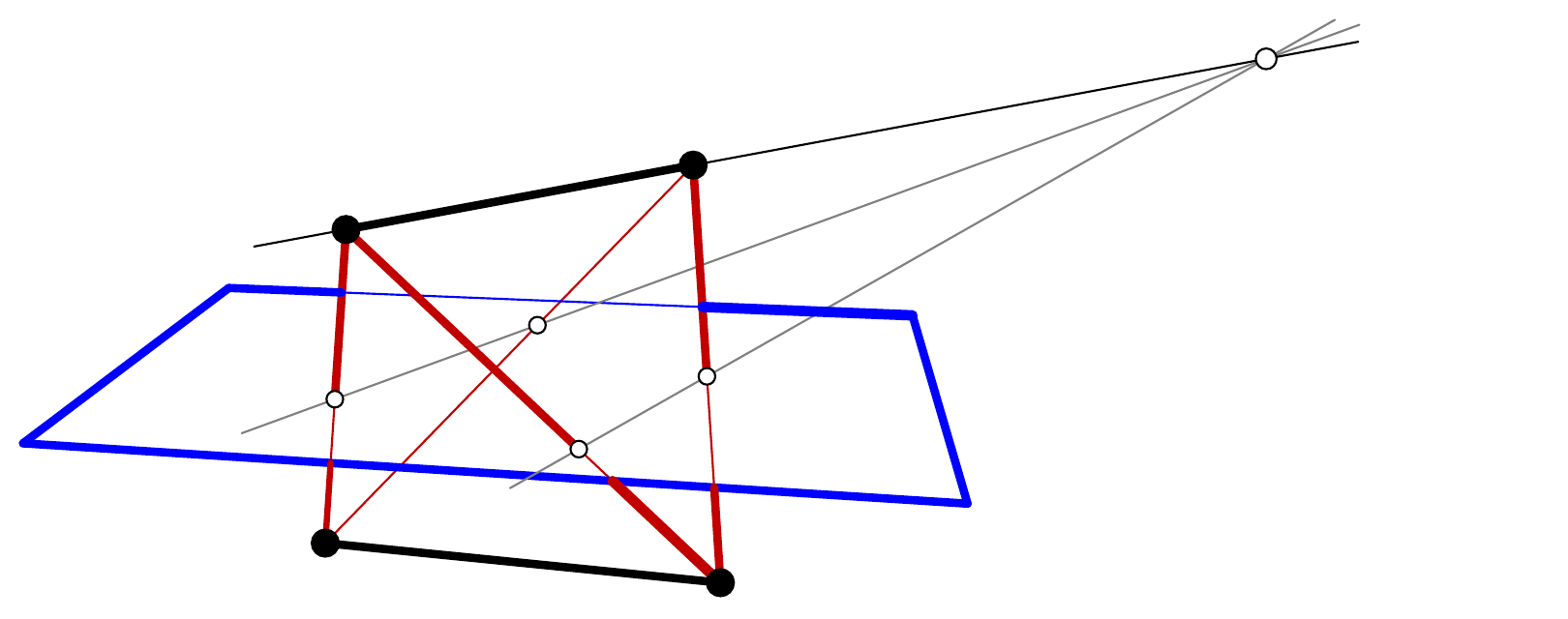}};
        
        \node at (-3,.7) {$p_1$};
        \node at (0,-1.7) {$p_2$};
        \node at (-.4,1.2) {$p_3$};
        \node at (-3.1,-1.7) {$p_4$};
        \node at (-3.1,-0.35) {$q_4$};
        \node at (-1.6,-0.85) {$q_1$};
        \node at (-0.7,-0.3) {$q_2$};
        \node at (-1.85,-0.06) {$q_3$};

    \end{tikzpicture}
    
    \caption{3D-Menelaus}\label{fig:3DM}
\end{figure}
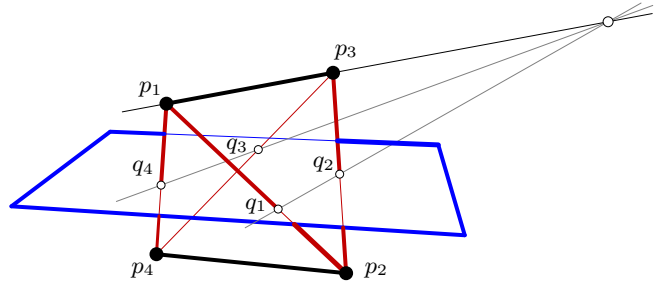

There are several ways to generalize Ceva's and Menelaus's theorem to higher dimensions. Here we will refer to the following two versions, whose proofs are immediate.
\begin{theorem} (3D Menelaus):
Let $p_1,\ldots,p_4$ be four non coplanar points in three-space. Consider a plane in general position $H$ and the intersections $q_i:=H\cap \overline{p_{i},p_{i+1}}$ (indices mod 4). Then we have (with oriented distances):
\[
{|p_1,q_1| \over |q_1,p_2|} \cdot
{|p_2,q_2| \over |q_2,p_3|} \cdot
{|p_3,q_3| \over |q_3,p_4|} \cdot
{|p_4,q_4| \over |q_4,p_1|} =1
\]
\end{theorem}

\medskip
\noindent
Similarly we have a version of Ceva's Theorem.
\begin{theorem} (3D Ceva):
Let $p_1,\ldots,p_4$ be four non coplanar points in three space. Consider one further point $a$.
For the line $\overline{p_{i},p_{i+1}}$ consider the plane $H_i:=a\vee p_{i+2}\vee p_{i+3}$ and the intersection point
 $q_i:=H_i\cap \overline{p_{i},p_{i+1}}$ (indices mod~4). Then we have:
\[
{|p_1,q_1| \over |q_1,p_2|} \cdot
{|p_2,q_2| \over |q_2,p_3|} \cdot
{|p_3,q_3| \over |q_3,p_4|} \cdot
{|p_4,q_4| \over |q_4,p_1|} =1
\]
\end{theorem}

The proofs of these two theorems follow exactly the same reasoning as in the 2D case presented in \cite{RG1}.
The situation of the 3D-Menelaus is shown in Figure~\ref{fig:3DM}. The blue plane plays the role of the plane $H$, cutting the edges of a 4-cycle in the tetrahedron. The image also indicates a simple argument for the multiratio equation. First observe that the three lines $q_1\vee q_2$,  $q_3\vee q_4$ and $p_1,p_3$ must meet in a point since they are the mutual intersections of three planes $H$, $p_1p_2p_3$, $p_1p_4p_3$.
Now consider the drawing literally as an orthogonal projection of the 3D situation. Then the 4-cycle can be considered as a version of two 2D Menelaus configurations glued  along the edge $p_1p_3$ along wich the joint factor cancels.
\smallskip 
From there it is also immediate to arrive at the corresponding 3D-Ceva theorem. For this observe that both Ceva and Menelaus  have the same multiratio{\color{Red}n} ($=1$) hence  the four points that split the cycle can as well be interpreted as Menelaus points or Ceva points. In Other words if $q_1,\ldots,q_4$ satisfy the Menelaus condition (they are coplaner) they also satisfy the Ceva condition (the planes $p_1p_2q_3$, $p_2p_3q_4$, $p_3p_4q_1$, $p_4p_1q_2$ meet in a point). The intersection  $p_1p_2q_3 \wedge p_3p_4q_1$ is the line $q_3q_1$ and the intersection $p_3p_4q_1 \wedge p_4p_1q_2$ is the line $q_2q_4$. Since the four points $q_i$ are coplanar these two lines (and hence the four planes) must intersect. The argument can easily be reversed.
\medskip 

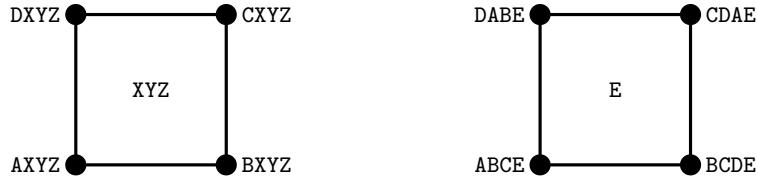
\begin{figure}[h]
    \begin{minipage}[c]{0.5\textwidth}
        \centering
        \begin{tikzpicture}[]
            \coordinate [label=left:{$\A\X\Y\Z\;$}] (A) at (0,0);
            \coordinate [label=right:{$\;\B\X\Y\Z$}] (B) at (2,0);
            \coordinate [label=right:{$\;\C\X\Y\Z$}] (C) at (2,2);
            \coordinate [label=left:{$\D\X\Y\Z\;$}] (D) at (0,2);
        
            \path [name path=A--B] (A) -- (B);
            \draw[very thick] (A) -- (B);
            \path [name path=B--C] (B) -- (C);
            \draw[very thick] (B) -- (C);
            \path [name path=C--D] (C) -- (D);
            \draw[very thick] (C) -- (D);
            \path [name path=D--A] (D) -- (A);
            \draw[very thick] (D) -- (A);
        
            \fill [black,opacity=1] (A) circle (4pt);
            \fill [black,opacity=1] (B) circle (4pt);
            \fill [black,opacity=1] (C) circle (4pt);
            \fill [black,opacity=1] (D) circle (4pt);
            \node (r42) at (1,1) {$\X\Y\Z$};

                   \end{tikzpicture}
                  
    \end{minipage}
     \begin{minipage}[c]{0.5\textwidth}
        \centering
        \begin{tikzpicture}[]
            \coordinate [label=left:{$\A\B\C\E\;$}] (A) at (0,0);
            \coordinate [label=right:{$\;\B\C\D\E$}] (B) at (2,0);
            \coordinate [label=right:{$\;\C\D\A\E$}] (C) at (2,2);
            \coordinate [label=left:{$\D\A\B\E\;$}] (D) at (0,2);
        
            \path [name path=A--B] (A) -- (B);
            \draw[very thick] (A) -- (B);
            \path [name path=B--C] (B) -- (C);
            \draw[very thick] (B) -- (C);
            \path [name path=C--D] (C) -- (D);
            \draw[very thick] (C) -- (D);
            \path [name path=D--A] (D) -- (A);
            \draw[very thick] (D) -- (A);
        
            \fill [black,opacity=1] (A) circle (4pt);
            \fill [black,opacity=1] (B) circle (4pt);
            \fill [black,opacity=1] (C) circle (4pt);
            \fill [black,opacity=1] (D) circle (4pt);
            \node (r42) at (1,1) {$\E$};

                   \end{tikzpicture}
                  
    \end{minipage}
    \caption{The labelings of 3D Ceva and Menelaus quadrangles.}\label{fig:3Dcycles}
\end{figure}

We will call these two configurations Menelaus and Ceva configurations, respectively. These CM 4-gons will be used in our manifold proofs. Figure~\ref{fig:3Dcycles} shows the labelings of these 4-gons if we use them in a base graph. The left picture shows a Menelaus quadrilateral. Where we assume that the vertices of the tetrahedron are $\A,\ldots, \D$ and the plane $H$ is is spanned by three points $\X\Y\Z$. The right picture shows the Ceva situation with a tetrahedron $\A,\ldots, \D$ and an interior point $\E$.

\subsection{The Sixteen Point Theorem}   
We now consider a spatial theorem that  is known under the name {\it Sixteen Point Theorem}.
Already in a very early preprint \cite{BRG2} on 
bi-quadratic final polynomials the following two incidence theorems, their bi-quadratic final polynomials, and their regular combinatorial properties was mentioned. We {\it revisit} these structures here
in order to emphasize how they can be translated into manifold proofs.

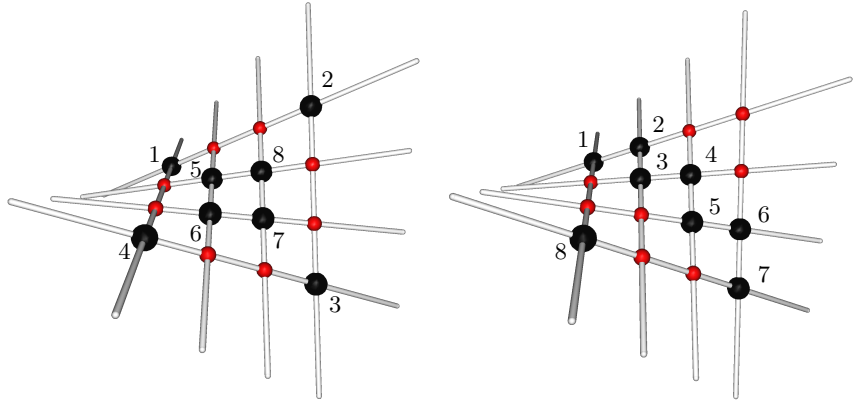
\begin{figure}[h]
    \centering
    \begin{tikzpicture}
        \draw (0, 0) node[inner sep=0] {\includegraphics[width=5.6cm]{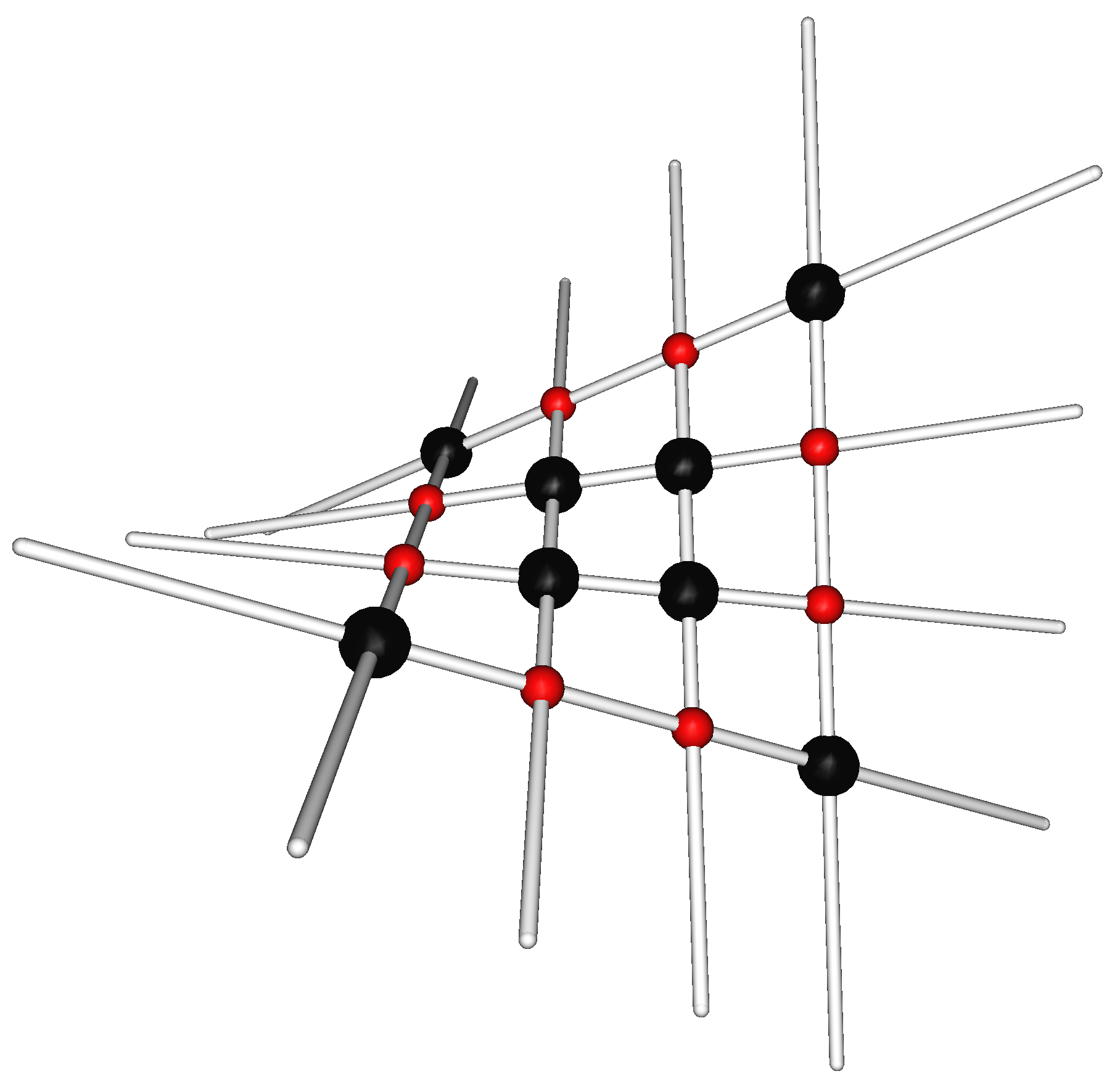}};
        \coordinate [label=left:{$1$}] (1) at (-.6,.6);
        \coordinate [label=left:{$4$}] (1) at (-1,-.7);
        \coordinate [label=left:{$3$}] (1) at (1.8,-1.4);
        \coordinate [label=left:{$2$}] (1) at (1.7,1.6);
        \coordinate [label=left:{$5$}] (1) at (-.05,.4);
        \coordinate [label=left:{$6$}] (1) at (-.05,-.45);
        \coordinate [label=left:{$7$}] (1) at (1.05,-.55);
        \coordinate [label=left:{$8$}] (1) at (1.05,.6);

        \draw (5.9, 0) node[inner sep=0] {\includegraphics[width=5.6cm]{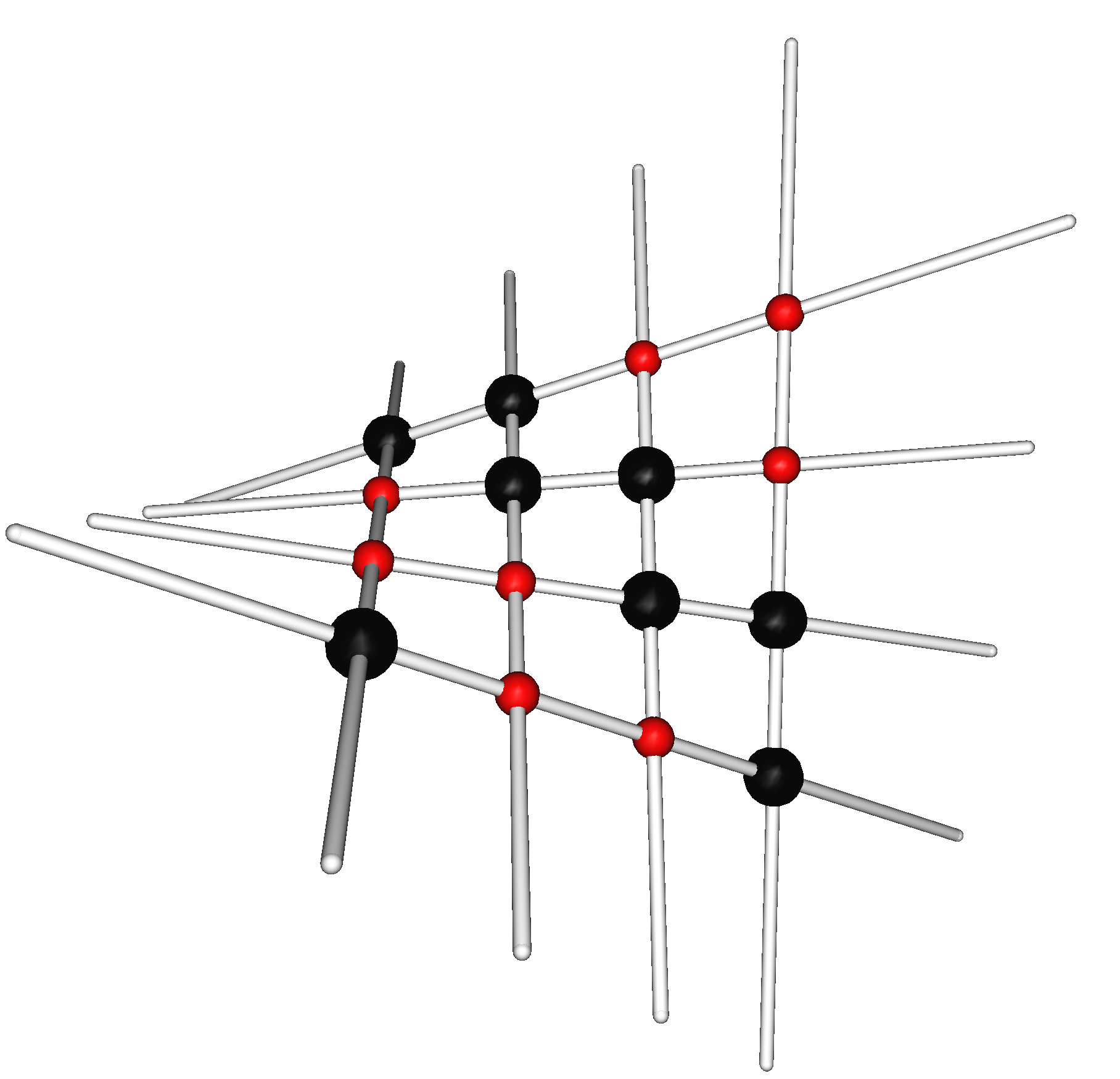}};
        
         \coordinate [label=left:{$1$}] (1) at (-.9+6,.8);
         \coordinate [label=left:{$2$}] (1) at (.1+6,1.0);
         \coordinate [label=left:{$3$}] (1) at (.15+6,0.5);
         \coordinate [label=left:{$4$}] (1) at (.8+6,0.6);
         \coordinate [label=left:{$5$}] (1) at (.85+6,-0.1);
         \coordinate [label=left:{$6$}] (1) at (1.5+6,-0.2);
         \coordinate [label=left:{$7$}] (1) at (1.5+6,-1.0);
         \coordinate [label=left:{$8$}] (1) at (-1.2+6,-0.7);

    \end{tikzpicture}
    \caption{Two versions of the 16-points theorem.}\label{fig:16pts}
\end{figure}

The 16 Point Theorem actually can be interpreted as a theorem about about two collections of 4-lines and their 16 points of potential mutual intersections.
Consider Figure~\ref{fig:16pts} for the geometric configuration. 
\begin{theorem}
Let $l_1,\ldots,l_4,m_1,\ldots,m_4,$ be two sets of lines (such that that all 8 lines are distinct).
If 15 of the 16 pairs of lines $\{l_i,m_j\}$ are incident then also the final one is incident.
\end{theorem}

\begin{figure}
    \centering
    \begin{minipage}[c]{0.3\linewidth}
        \centering
        \begin{tikzpicture}[>=Stealth,scale=0.8]
 
            \node (r00) at (0,0) {$\4$};
            \node (r10) at (1,0) {$\8$};
            \node (r20) at (2,0) {$\1$};
            \node (r30) at (3,0) {$\5$};
            \node (r40) at (4,0) {$\4$};
           \node (r01) at (0,1) {$\3$};
            \node (r11) at (1,1) {$\7$};
            \node (r21) at (2,1) {$\2$};
            \node (r31) at (3,1) {$\6$};
            \node (r41) at (4,1) {$\3$};
           \node (r02) at (0,2) {$\1$};
            \node (r12) at (1,2) {$\5$};
            \node (r22) at (2,2) {$\4$};
            \node (r32) at (3,2) {$\8$};
            \node (r42) at (4,2) {$\1$};
            \coordinate (a11) at (-.3,0);
            \coordinate (a12) at (-.3,2);
           \coordinate (a21) at (4.3,0);
            \coordinate (a22) at (4.3,2);
            \coordinate (a31) at (0,-.4);
            \coordinate (a32) at (4,-.4);
            \coordinate (a41) at (2.1,2.4);
            \coordinate (a42) at (4,2.4);
            \coordinate (a51) at (0,2.4);
            \coordinate (a52) at (1.9,2.4);

            \draw[ thick] (r00) -- (r10);
            \draw[ thick] (r10) -- (r20);
            \draw[ thick] (r20) -- (r30);
            \draw[ thick] (r30) -- (r40);
            \draw[ thick] (r01) -- (r11);
            \draw[ thick] (r11) -- (r21);
            \draw[ thick] (r21) -- (r31);
            \draw[ thick] (r31) -- (r41);
            \draw[ thick] (r02) -- (r12);
            \draw[ thick] (r12) -- (r22);
            \draw[ thick] (r22) -- (r32);
            \draw[ thick] (r32) -- (r42);
            
            \draw[ thick] (r00) -- (r01);
            \draw[ thick] (r01) -- (r02);
            \draw[ thick] (r10) -- (r11);
            \draw[ thick] (r11) -- (r12);
            \draw[ thick] (r20) -- (r21);
            \draw[ thick] (r21) -- (r22);
            \draw[ thick] (r30) -- (r31);
            \draw[ thick] (r31) -- (r32);
            \draw[ thick] (r40) -- (r41);
            \draw[ thick] (r41) -- (r42);
            \draw[->, red, line width=1.2pt] (a11) -- (a12);
            \draw[->, red, line width=1.2pt] (a21) -- (a22);
            \draw[-, red, line width=1.2pt] ($(a11)+(.1,0)$) -- ($(a11)+(-.1,0)$);
            \draw[-, red, line width=1.2pt] ($(a21)+(.1,0)$) -- ($(a21)+(-.1,0)$);

            \draw[->, green!60!black, line width=1.2pt] (a31) -- (a32);
            \draw[->, green!60!black, line width=1.2pt] (a41) -- (a42);
            \draw[->, green!60!black, line width=1.2pt] (a51) -- (a52);
            \draw[-, green!60!black, line width=1.2pt] ($(a31)+(0,.1)$) -- ($(a31)+(0,-.1)$);
            \draw[-, green!60!black, line width=1.2pt] ($(a41)+(0,.1)$) -- ($(a41)+(0,-.1)$);

        \end{tikzpicture}
    \end{minipage}
    \hspace{0.9cm}
    \begin{minipage}[c]{0.6\linewidth}
        \centering
        \begin{tikzpicture}[>=Stealth,scale=0.8]
                  \node (r00) at (0,0) {$\2$};
            \node (r10) at (1,0) {$\5$};
            \node (r20) at (2,0) {$\8$};
            \node (r30) at (3,0) {$\3$};
            \node (r40) at (4,0) {$\6$};
            \node (r50) at (5,0) {$\1$};
            \node (r60) at (6,0) {$\4$};
            \node (r70) at (7,0) {$\2$};
            \node (r80) at (8,0) {$\3$};
           \node (r01) at (0,1) {$\1$};
            \node (r11) at (1,1) {$\4$};
            \node (r21) at (2,1) {$\7$};
            \node (r31) at (3,1) {$\2$};
            \node (r41) at (4,1) {$\5$};
            \node (r51) at (5,1) {$\8$};
            \node (r61) at (6,1) {$\3$};
            \node (r71) at (7,1) {$\6$};
            \node (r81) at (8,1) {$\1$};

            \coordinate (a11) at (-.3,0);
            \coordinate (a12) at (-.3,1);
           \coordinate (a21) at (8.3,0);
            \coordinate (a22) at (8.3,1);
            \coordinate (a31) at (0,-.4);
            \coordinate (a32) at (8,-.4);
            \coordinate (a41) at (3.1,1.4);
            \coordinate (a42) at (8,1.4);
            \coordinate (a51) at (0,1.4);
            \coordinate (a52) at (2.9,1.4);

            \draw[ thick] (r00) -- (r10);
            \draw[ thick] (r10) -- (r20);
            \draw[ thick] (r20) -- (r30);
            \draw[ thick] (r30) -- (r40);
            \draw[ thick] (r40) -- (r50);
            \draw[ thick] (r50) -- (r60);
            \draw[ thick] (r60) -- (r70);
            \draw[ thick] (r70) -- (r80);
            \draw[ thick] (r01) -- (r11);
            \draw[ thick] (r11) -- (r21);
            \draw[ thick] (r21) -- (r31);
            \draw[ thick] (r31) -- (r41);
            \draw[ thick] (r41) -- (r51);
            \draw[ thick] (r51) -- (r61);
            \draw[ thick] (r61) -- (r71);
            \draw[ thick] (r71) -- (r81);

            \draw[ thick] (r00) -- (r01);
            \draw[ thick] (r10) -- (r11);
            \draw[ thick] (r20) -- (r21);
            \draw[ thick] (r30) -- (r31);
            \draw[ thick] (r40) -- (r41);
            \draw[ thick] (r50) -- (r51);
            \draw[ thick] (r60) -- (r61);
            \draw[ thick] (r70) -- (r71);
            \draw[ thick] (r80) -- (r81);
                    
             \draw[->, red, line width=1.2pt] (a11) -- (a12);
            \draw[->, red, line width=1.2pt] (a21) -- (a22);
            \draw[-, red, line width=1.2pt] ($(a11)+(.1,0)$) -- ($(a11)+(-.1,0)$);
            \draw[-, red, line width=1.2pt] ($(a21)+(.1,0)$) -- ($(a21)+(-.1,0)$);

            \draw[->, green!60!black, line width=1.2pt] (a31) -- (a32);
            \draw[->, green!60!black, line width=1.2pt] (a41) -- (a42);
            \draw[->, green!60!black, line width=1.2pt] (a51) -- (a52);
            \draw[-, green!60!black, line width=1.2pt] ($(a31)+(0,.1)$) -- ($(a31)+(0,-.1)$);
            \draw[-, green!60!black, line width=1.2pt] ($(a41)+(0,.1)$) -- ($(a41)+(0,-.1)$);

        \end{tikzpicture}
    \end{minipage}
    \caption{The combinatorics of the two versions of the 16 point theorem.}\label{fig:tori}
\end{figure}
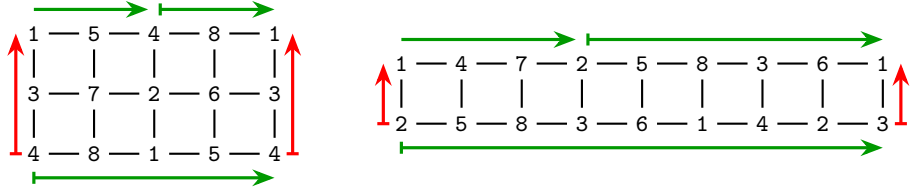

One can rephrase the theorem as an incidence theorem on 8 points and 8 coplanarity relations. This can be done in two combinatorially non-isomorphic ways. They are depicted in the two images of Figure~\ref{fig:16pts}.
The points of the configuration are shown as black dots. Each line is spanned by exactly two such points. Thus
each black point already encodes an incidence of two lines. 
The small red dots correspond to the 8 remaining incidences. The two lines $\I \vee \J$ and  $\K \vee \L$
are incident if and only if $\I, \J,\K , \L$ are coplanar. Hence in each of the two versions seven coplanarities imply an eighth one. The combinatorics of both types amazingly  can also be interpreted as a torus consisting  of 8 quadrilaterals. Figure~\ref{fig:tori} shows the combinatorics for the two cases. Each square represents a coplanarity between points.

In \cite{RG2} bi-quadratic final polynomials for both cases were given. We now translate these bi-quadratic final polynomials into CM-proofs and into quad-proofs. We start with the case associated to the structure in the left of Figures~\ref{fig:16pts} and~\ref{fig:tori}.
A final polynomial for this one  is given by (see \cite{RG2}):
\[
\begin{matrix}
[\1\2\5\3][\1\2\6\4] = [\1\2\5\4][\1\2\6\3] \quad\Leftarrow\quad h(\1\5\2\6)\\
[\2\3\6\1][\2\3\7\4] = [\2\3\6\4][\2\3\7\1] \quad\Leftarrow\quad h(\2\6\3\7)\\
[\3\4\5\1][\3\4\6\2] = [\3\4\5\2][\3\4\6\1] \quad\Leftarrow\quad h(\6\4\5\3)\\
[\1\4\6\3][\1\4\7\2] = [\1\4\6\2][\1\4\7\3] \quad\Leftarrow\quad h(\6\4\7\1)\\
[\3\4\7\1][\3\4\8\2] = [\3\4\7\2][\3\4\8\1] \quad\Leftarrow\quad h(\4\8\3\7)\\
[\1\4\5\2][\1\4\8\3] = [\1\4\5\3][\1\4\8\2] \quad\Leftarrow\quad h(\4\8\1\5)\\
[\1\2\7\3][\1\2\8\4] = [\1\2\7\4][\1\2\8\3] \quad\Leftarrow\quad h(\8\2\7\1)\\
[\2\3\5\4][\2\3\8\1] = [\2\3\5\1][\2\3\8\4]  \quad\Rightarrow\quad h(\8\2\5\3)
\end{matrix}
\]
    
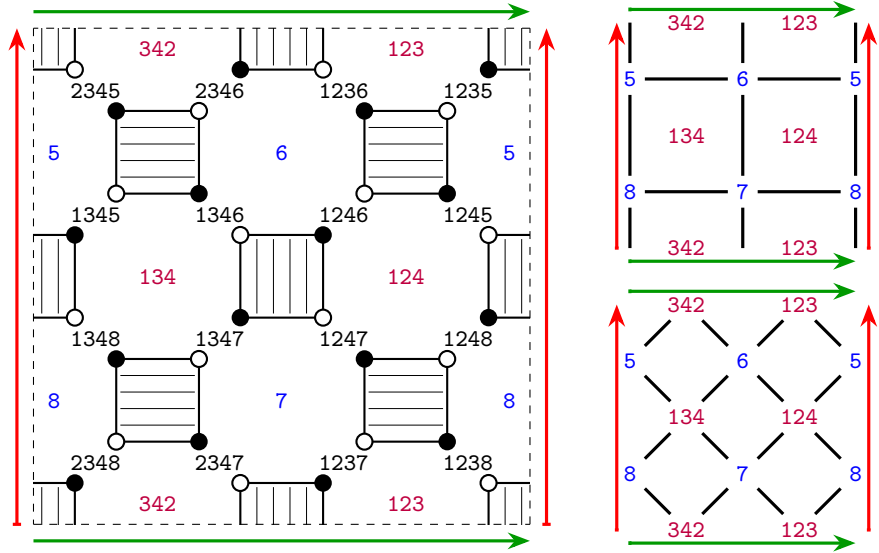
\begin{figure}
    \centering
    \begin{minipage}[c]{0.4\linewidth}
        \centering
        \begin{tikzpicture}[>=Stealth,scale=0.55]

            \coordinate (b00) at (0,0);
            \coordinate (b10) at (12,0);
            \coordinate (b11) at (12,12);
            \coordinate (b01) at (0,12);
         
            \draw[dashed, line width=.4pt] (b00) -- (b10);
            \draw[dashed, line width=.4pt] (b10) -- (b11);
            \draw[dashed, line width=.4pt] (b11) -- (b01);
            \draw[dashed, line width=.4pt] (b01) -- (b00);

           \begin{scope}
    \clip (0,0) rectangle (12,12);

    \foreach \i in {-2,-1,0,1,2,3,4} {

    \foreach \j in {-2,-1,0,1,2,3,4} {
            \draw[line width=.8pt] (\i*3+2,\j*6+\i*3+2) rectangle ++(2,2);
             \ifodd\i
              \draw[line width=.3pt] (\i*3+2+.4,\j*6+\i*3+2+1.9) --(\i*3+2+.4,\j*6+\i*3+2+.1) ;
              \draw[line width=.3pt] (\i*3+2+.8,\j*6+\i*3+2+1.9) --(\i*3+2+.8,\j*6+\i*3+2+.1) ;
              \draw[line width=.3pt] (\i*3+2+1.2,\j*6+\i*3+2+1.9) --(\i*3+2+1.2,\j*6+\i*3+2+.1) ;
              \draw[line width=.3pt] (\i*3+2+1.6,\j*6+\i*3+2+1.9) --(\i*3+2+1.6,\j*6+\i*3+2+.1) ;

            \fill [black,opacity=1] (\i*3+2,\j*6+\i*3+2) circle (6pt);
            \fill [black,opacity=1] (\i*3+2+2,\j*6+\i*3+2+2) circle (6pt);
            \path[fill=white,draw=black, thick] (\i*3+2,\j*6+\i*3+2+2) circle[radius=5.5pt];
            \path[fill=white,draw=black, thick] (\i*3+2+2,\j*6+\i*3+2) circle[radius=5.5pt];

            \else
            \draw[line width=.3pt] (\i*3+2+.1,\j*6+\i*3+2+.4) --(\i*3+2+1.9,\j*6+\i*3+2+.4) ;
            \draw[line width=.3pt] (\i*3+2+.1,\j*6+\i*3+2+.8) --(\i*3+2+1.9,\j*6+\i*3+2+.8) ;
            \draw[line width=.3pt] (\i*3+2+.1,\j*6+\i*3+2+1.2) --(\i*3+2+1.9,\j*6+\i*3+2+1.2) ;
            \draw[line width=.3pt] (\i*3+2+.1,\j*6+\i*3+2+1.6) --(\i*3+2+1.9,\j*6+\i*3+2+1.6) ;

            \fill [black,opacity=1] (\i*3+2+2,\j*6+\i*3+2) circle (6pt);
            \fill [black,opacity=1] (\i*3+2,\j*6+\i*3+2+2) circle (6pt);
                \path[fill=white,draw=black, thick] (\i*3+2,\j*6+\i*3+2) circle[radius=5.5pt];
            \path[fill=white,draw=black, thick] (\i*3+2+2,\j*6+\i*3+2+2) circle[radius=5.5pt];

            \fi

         }
       }     

  \end{scope}
  \node (r) at (1.5,1.5) {$\tt 2348$};
  \node (r) at (4.5,1.5) {$\tt 2347$};
  \node (r) at (7.5,1.5) {$\tt 1237$};
  \node (r) at (10.5,1.5) {$\tt 1238$};
  \node (r) at (1.5,4.5) {$\tt 1348$};
  \node (r) at (4.5,4.5) {$\tt 1347$};
  \node (r) at (7.5,4.5) {$\tt 1247$};
  \node (r) at (10.5,4.5) {$\tt 1248$};
  \node (r) at (1.5,7.5) {$\tt 1345$};
  \node (r) at (4.5,7.5) {$\tt 1346$};
  \node (r) at (7.5,7.5) {$\tt 1246$};
  \node (r) at (10.5,7.5) {$\tt 1245$};
  \node (r) at (1.5,10.5) {$\tt 2345$};
  \node (r) at (4.5,10.5) {$\tt 2346$};
  \node (r) at (7.5,10.5) {$\tt 1236$};
  \node (r) at (10.5,10.5) {$\tt 1235$};
  \node[blue] (r) at (0.5,3) {$\tt 8$};
  \node[blue] (r) at (6,3) {$\tt 7$};
  \node[blue] (r) at (11.5,3) {$\tt 8$};
 \node[blue] (r) at (0.5,9) {$\tt 5$};
  \node[blue] (r) at (6,9) {$\tt 6$};
  \node[blue] (r) at (11.5,9) {$\tt 5$};
  \node[purple] (r) at (3,6) {$\tt 134$};
  \node[purple] (r) at (9,6) {$\tt 124$};
  \node[purple] (r) at (3,.5) {$\tt 342$};
  \node[purple] (r) at (9,.5) {$\tt 123$};
    \node[purple] (r) at (3,11.5) {$\tt 342$};
  \node[purple] (r) at (9,11.5) {$\tt 123$};
     
       \coordinate (a31) at (0,-0.4);
       \coordinate (a32) at (12,-0.4);
           \coordinate (a41) at (0,12.4);
            \coordinate (a42) at (12.0,12.4);
            \coordinate (a11) at (-.4,0);
            \coordinate (a12) at (-.4,12);
            \coordinate (a21) at (12.4,0);
            \coordinate (a22) at (12.4,12);
            \draw[->, red, line width=1.2pt] (a11) -- (a12);
            \draw[->, red, line width=1.2pt] (a21) -- (a22);

             \draw[-, red, line width=1.2pt] ($(a11)+(.1,0)$) -- ($(a11)+(-.1,0)$);
            \draw[-, red, line width=1.2pt] ($(a21)+(.1,0)$) -- ($(a21)+(-.1,0)$);

            \draw[->, green!60!black, line width=1.2pt] (a31) -- (a32);
            \draw[->, green!60!black, line width=1.2pt] (a41) -- (a42);

               \end{tikzpicture}
    \end{minipage}
    \hfill
    \begin{minipage}[c]{0.4\linewidth}
        \centering
        \begin{tikzpicture}[>=Stealth,scale=0.25]

            \coordinate (b00) at (0,0);
            \coordinate (b10) at (12,0);
            \coordinate (b11) at (12,12);
            \coordinate (b01) at (0,12);

   \node[blue] (81) at (0,3) {$\tt 8$};
  \node[blue] (7) at (6,3) {$\tt 7$};
  \node[blue] (82) at (12,3) {$\tt 8$};
 \node[blue] (51) at (0,9) {$\tt 5$};
  \node[blue] (6) at (6,9) {$\tt 6$};
  \node[blue] (52) at (12,9) {$\tt 5$};
  \node[purple] (134) at (3,6) {$\tt 134$};
  \node[purple] (124) at (9,6) {$\tt 124$};
  \node[purple] (3421) at (3,0) {$\tt 342$};
  \node[purple] (1231) at (9,0) {$\tt 123$};
    \node[purple] (3422) at (3,12) {$\tt 342$};
  \node[purple] (1232) at (9,12) {$\tt 123$};
     
       \coordinate (a31) at (0,-0.7);
       \coordinate (a32) at (12,-0.7);
           \coordinate (a41) at (0,12.7);
            \coordinate (a42) at (12.0,12.7);
            \coordinate (a11) at (-.7,0);
            \coordinate (a12) at (-.7,12);
            \coordinate (a21) at (12.7,0);
            \coordinate (a22) at (12.7,12);
            \draw[->, red, line width=1.2pt] (a11) -- (a12);
            \draw[->, red, line width=1.2pt] (a21) -- (a22);

             \draw[-, red, line width=1.2pt] ($(a11)+(.1,0)$) -- ($(a11)+(-.1,0)$);
            \draw[-, red, line width=1.2pt] ($(a21)+(.1,0)$) -- ($(a21)+(-.1,0)$);

            \draw[->, green!60!black, line width=1.2pt] (a31) -- (a32);
            \draw[->, green!60!black, line width=1.2pt] (a41) -- (a42);
            \draw[-, green!60!black, line width=1.2pt] ($(a31)+(0,.1)$) -- ($(a31)+(0,-.1)$);
            \draw[-, green!60!black, line width=1.2pt] ($(a41)+(0,.1)$) -- ($(a41)+(0,-.1)$);
            \draw[green!0!black, line width=1.2pt] (7) -- (134);
            \draw[green!0!black, line width=1.2pt] (7) -- (1231);
            \draw[green!0!black, line width=1.2pt] (7) -- (3421);
            \draw[green!0!black, line width=1.2pt] (7) -- (124);
            \draw[green!0!black, line width=1.2pt] (6) -- (134);
            \draw[green!0!black, line width=1.2pt] (6) -- (1232);
            \draw[green!0!black, line width=1.2pt] (6) -- (3422);
            \draw[green!0!black, line width=1.2pt] (6) -- (124);
            \draw[green!0!black, line width=1.2pt] (51) -- (134);
            \draw[green!0!black, line width=1.2pt] (81) -- (134);
            \draw[green!0!black, line width=1.2pt] (81) -- (3421);
            \draw[green!0!black, line width=1.2pt] (51) -- (3422);
            \draw[green!0!black, line width=1.2pt] (52) -- (124);
            \draw[green!0!black, line width=1.2pt] (82) -- (124);
            \draw[green!0!black, line width=1.2pt] (82) -- (1231);
            \draw[green!0!black, line width=1.2pt] (52) -- (1232);

       \begin{scope}[shift={(0,15)}]

   \node[blue] (81) at (0,3) {$\tt 8$};
  \node[blue] (7) at (6,3) {$\tt 7$};
  \node[blue] (82) at (12,3) {$\tt 8$};
 \node[blue] (51) at (0,9) {$\tt 5$};
  \node[blue] (6) at (6,9) {$\tt 6$};
  \node[blue] (52) at (12,9) {$\tt 5$};
  \node[purple] (134) at (3,6) {$\tt 134$};
  \node[purple] (124) at (9,6) {$\tt 124$};
  \node[purple] (3421) at (3,0) {$\tt 342$};
  \node[purple] (1231) at (9,0) {$\tt 123$};
    \node[purple] (3422) at (3,12) {$\tt 342$};
  \node[purple] (1232) at (9,12) {$\tt 123$};
     
       \coordinate (a31) at (0,-0.7);
       \coordinate (a32) at (12,-0.7);
           \coordinate (a41) at (0,12.7);
            \coordinate (a42) at (12.0,12.7);
            \coordinate (a11) at (-.7,0);
            \coordinate (a12) at (-.7,12);
            \coordinate (a21) at (12.7,0);
            \coordinate (a22) at (12.7,12);
            \draw[->, red, line width=1.2pt] (a11) -- (a12);
            \draw[->, red, line width=1.2pt] (a21) -- (a22);

             \draw[-, red, line width=1.2pt] ($(a11)+(.1,0)$) -- ($(a11)+(-.1,0)$);
            \draw[-, red, line width=1.2pt] ($(a21)+(.1,0)$) -- ($(a21)+(-.1,0)$);

            \draw[->, green!60!black, line width=1.2pt] (a31) -- (a32);
            \draw[->, green!60!black, line width=1.2pt] (a41) -- (a42);
            \draw[-, green!60!black, line width=1.2pt] ($(a31)+(0,.1)$) -- ($(a31)+(0,-.1)$);
            \draw[-, green!60!black, line width=1.2pt] ($(a41)+(0,.1)$) -- ($(a41)+(0,-.1)$);
            \draw[green!0!black, line width=1.2pt] (7) -- (81);
            \draw[green!0!black, line width=1.2pt] (7) -- (82);
            \draw[green!0!black, line width=1.2pt] (7) -- (6);
            \draw[green!0!black, line width=1.2pt] (6) -- (51);
            \draw[green!0!black, line width=1.2pt] (6) -- (52);
            \draw[green!0!black, line width=1.2pt] (6) -- (6,12);
            \draw[green!0!black, line width=1.2pt] (51) -- (0,12);
            \draw[green!0!black, line width=1.2pt] (52) -- (12,12);
            \draw[green!0!black, line width=1.2pt] (7) -- (6,0);
            \draw[green!0!black, line width=1.2pt] (81) -- (0,0);
            \draw[green!0!black, line width=1.2pt] (82) -- (12,0);
            
             \draw[green!0!black, line width=1.2pt] (81) -- (51);
             \draw[green!0!black, line width=1.2pt] (82) -- (52);

          \end{scope}

               \end{tikzpicture}
    \end{minipage}
    \caption{Structure diagrams for the manifold proofs of the theorem.}\label{fig:structurev1}
\end{figure}

Inspecting the structure and arranging the bi-quadratic equations  along a manifold like structure  yields the picture in Figure~\ref{fig:structurev1} (left). The image has to be interpreted as a torus,  with opposite edges identified. Each shaded cell corresponds to a bi-quadratic equation with black vertices appearing on the left and white vertices appearing on the right of the equation. Obviously at each bracket a white and a black vertex meet yielding the decisive cancellation pattern. The regularity of this cancellation pattern was already mentioned in \cite{BRG2}.

Surrounded by the squares of the bi-quadratic equations there are eight other {\it empty} squares in this diagram. The labels in their center correspond to the points that are shared by the four adjacent brackets. There are cells  labelled by one  point that correspond to a 3D Ceva cycle and cells labelled by three points corresponding to a 3D Menelaus cycle. It can be easly checked, that these cycles indeed correspond to Ceva and Menelaus configurations in the original incidence configuration. In the Menelaus{\color{ForestGreen}-}like cells the three points span the cutting plane for the Menelaus configuration. In the Ceva cell the one letter in the inside characterizes the one additional point needed for a Ceva configuration.

\begin{figure}
    \centering
    \begin{minipage}[c]{0.4\linewidth}
        \centering
        \begin{tikzpicture}[>=Stealth,scale=0.2]

       \coordinate (a31) at (0,-0.9);
       \coordinate (a32) at (12,-0.9);
           \coordinate (a41) at (0,12.9);
            \coordinate (a42) at (12.0,12.9);
            \coordinate (a11) at (-.9,0);
            \coordinate (a12) at (-.9,12);
            \coordinate (a21) at (12.9,0);
            \coordinate (a22) at (12.9,12);
            \draw[->, red, line width=1.2pt] (a11) -- (a12);
            \draw[->, red, line width=1.2pt] (a21) -- (a22);
            \draw[->, green!60!black, line width=1.2pt] (a31) -- (a32);
            \draw[->, green!60!black, line width=1.2pt] (a41) -- (a42);
         
                     \fill [black,opacity=1] (0,0) circle (12pt);
                     \fill [black,opacity=1] (0,6) circle (12pt);
                     \fill [black,opacity=1] (0,12) circle (12pt);
                     \fill [black,opacity=1] (6,0) circle (12pt);
                     \fill [black,opacity=1] (6,6) circle (12pt);
                     \fill [black,opacity=1] (6,12) circle (12pt);
                     \fill [black,opacity=1] (12,0) circle (12pt);
                     \fill [black,opacity=1] (12,6) circle (12pt);
                     \fill [black,opacity=1] (12,12) circle (12pt);
            \draw[black, line width=1.2pt] ((0,0) -- (0,12);
            \draw[black, line width=1.2pt] ((6,0) -- (6,12);
            \draw[black, line width=1.2pt] ((12,0) -- (12,12);
            \draw[black, line width=1.2pt] ((0,0) -- (12,0);
            \draw[black, line width=1.2pt] ((0,6) -- (12,6);
            \draw[black, line width=1.2pt] ((0,12) -- (12,12);

         \node[purple]  at (3,3) {\large$ M$};
         \node[purple]  at (9,3) {\large$ M$};
         \node[purple]  at (3,9) {\large$ M$};
         \node[purple]  at (9,9) {\large$ M$};
               \end{tikzpicture}
    \end{minipage}
        \hspace{0.0cm}
    \begin{minipage}[c]{0.4\linewidth}
    
             \begin{tikzpicture}[>=Stealth,scale=0.2]

       \coordinate (a31) at (0,-0.9);
       \coordinate (a32) at (12,-0.9);
           \coordinate (a41) at (0,12.9);
            \coordinate (a42) at (12.0,12.9);
            \coordinate (a11) at (-.9,0);
            \coordinate (a12) at (-.9,12);
            \coordinate (a21) at (12.9,0);
            \coordinate (a22) at (12.9,12);
            \draw[->, red, line width=1.2pt] (a11) -- (a12);
            \draw[->, red, line width=1.2pt] (a21) -- (a22);
            \draw[->, green!60!black, line width=1.2pt] (a31) -- (a32);
            \draw[->, green!60!black, line width=1.2pt] (a41) -- (a42);
         
                     \fill [black,opacity=1] (0,0) circle (12pt);
                     \fill [black,opacity=1] (0,6) circle (12pt);
                     \fill [black,opacity=1] (0,12) circle (12pt);
                     \fill [black,opacity=1] (6,0) circle (12pt);
                     \fill [black,opacity=1] (6,6) circle (12pt);
                     \fill [black,opacity=1] (6,12) circle (12pt);
                     \fill [black,opacity=1] (12,0) circle (12pt);
                     \fill [black,opacity=1] (12,6) circle (12pt);
                     \fill [black,opacity=1] (12,12) circle (12pt);
            \draw[black, line width=1.2pt] ((0,0) -- (0,12);
            \draw[black, line width=1.2pt] ((6,0) -- (6,12);
            \draw[black, line width=1.2pt] ((12,0) -- (12,12);
            \draw[black, line width=1.2pt] ((0,0) -- (12,0);
            \draw[black, line width=1.2pt] ((0,6) -- (12,6);
            \draw[black, line width=1.2pt] ((0,12) -- (12,12);

         \node[purple]  at (3,3) {\large$ C$};
         \node[purple]  at (9,3) {\large$ C$};
         \node[purple]  at (3,9) {\large$ C$};
         \node[purple]  at (9,9) {\large$ C$};

               \end{tikzpicture}
    \end{minipage}
    \caption{Two proving schemes for the CM-proof. \label{fig:proof}}
\end{figure}
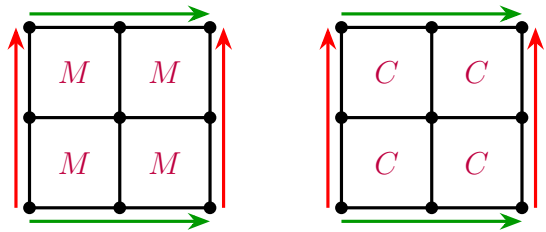

As before{\color{ForestGreen},} we may consider the binomial equations as glue between two Ceva or Menelaus cells. Each such bi-quadratic square can be used to encode the identity of two opposite edges and we can decide which pair
of edges is identified. In the drawing of  Figure~\ref{fig:structurev1} (left)
the squares  are shaded indicating a direction. We assume that the squares are used to identify the rations along the edges in the same direction as the shading.
For instance the upper left square whose vertices are ${\tt 2345}, {\tt 2346}, {\tt 1346}, {\tt 1345}$ is used to encode the ratio equivalence  ${[{\tt 2345}]\over [{\tt 2346}]}={[{\tt 1345}]\over [ {\tt 1346}]}$.
Doing all these identifications turns the diagram into a CM-proof that entirely consists of Menelaus configurations. It is shown in Figure~\ref{fig:structurev1} (right, top). This CM-proof simply consists of a $2\times 2$ grid with identified opposite edges. It has the topology of a torus. Each cell corresponds to a 3D Menelaus configuration. If we consequently identify the other edges of the bi-quadratic squares. We get a similar proof where each cell is a 3D Ceva configuration. Each such situation is already completely determined by the diagrams in Figure~\ref{fig:proof}. and the dimension in which the configuration should take place.

We next show how the structure can be used to create a quad-proof. The corresponding structure is shown
in Figure~\ref{fig:structurev1} (right, bottom). Again the underlying manifold is a torus. This time consisting of eight 4-gons. The vertices the manifold consist of are 4 points and 4 hyperplanes. The hyperplanes correspond to the red labels with 3  letters. The vertices are the blue labels with one letter. The same manifold scheme was also given in \cite{FP} (Theorem 5.5).

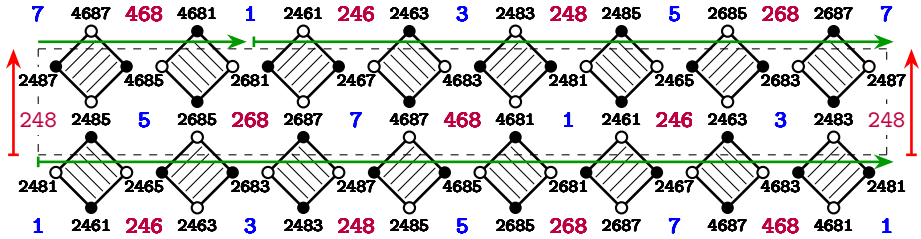
\begin{figure}
    \centering
    \begin{minipage}[c]{1\linewidth}
        \centering
        \begin{tikzpicture}[>=Stealth,scale=0.47]

       \coordinate (b00) at (0,0+2.5);
            \coordinate (b10) at (24,0+2.5);
            \coordinate (b11) at (24,3+2.5);
            \coordinate (b01) at (0,3+2.5);
         
            \draw[dashed, line width=.4pt] (b00) -- (b10);
            \draw[dashed, line width=.4pt] (b10) -- (b11);
            \draw[dashed, line width=.4pt] (b11) -- (b01);
            \draw[dashed, line width=.4pt] (b01) -- (b00);
              \node[fill=white,font=\scriptsize]  at (24,4.6) {$\tt 2487$};
              \node[fill=white,font=\scriptsize]  at (0,4.6) {$\tt 2487$};

           \begin{scope}

    \foreach \i in {0,1,2,3,4,5,6,7} {

    \foreach \j in {0,1} {
        \draw[line width=1pt] (\i*3+.5,\j*3+2) -- (\i*3+1.5,\j*3+1) ;
                     \draw[line width=1pt] (\i*3+.5,\j*3+2) -- (\i*3+1.5,\j*3+3) ;
                     \draw[line width=1pt] (\i*3+.5+2,\j*3+2) -- (\i*3+1.5,\j*3+1) ;
                     \draw[line width=1pt] (\i*3+.5+2,\j*3+2) -- (\i*3+1.5,\j*3+3) ;
    \pgfmathtruncatemacro{\s}{\i+\j}

             \ifodd\s
                      \draw[line width=.4 pt] (\i*3+.5+.2+.07,\j*3+2-.2+.07) -- (\i*3+1.5+.2-.07,\j*3+3-.2-.07) ;
                      \draw[line width=.4 pt] (\i*3+.5+.4+.07,\j*3+2-.4+.07) -- (\i*3+1.5+.4-.07,\j*3+3-.4-.07) ;
                      \draw[line width=.4 pt] (\i*3+.5+.6+.07,\j*3+2-.6+.07) -- (\i*3+1.5+.6-.07,\j*3+3-.6-.07) ;
                      \draw[line width=.4 pt] (\i*3+.5+.8+.07,\j*3+2-.8+.07) -- (\i*3+1.5+.8-.07,\j*3+3-.8-.07) ;

                      \fill [black,opacity=1] (\i*3+.5,\j*3+2) circle (5pt);
                      \fill [black,opacity=1] (\i*3+.5+2,\j*3+2) circle (5pt);
                   \path[fill=white,draw=black, thick] (\i*3+1.5,\j*3+1) circle[radius=4.5pt];
                   \path[fill=white,draw=black, thick] (\i*3+1.5,\j*3+3) circle[radius=4.5pt];

            \else
                       \draw[line width=.4 pt] (\i*3+.5+2-.07-.2,\j*3+2-.2+.07) -- (\i*3+1.5-.2+.07,\j*3+3-.2-.07);
                       \draw[line width=.4 pt] (\i*3+.5+2-.07-.4,\j*3+2-.4+.07) -- (\i*3+1.5-.4+.07,\j*3+3-.4-.07);
                       \draw[line width=.4 pt] (\i*3+.5+2-.07-.6,\j*3+2-.6+.07) -- (\i*3+1.5-.6+.07,\j*3+3-.6-.07);
                       \draw[line width=.4 pt] (\i*3+.5+2-.07-.8,\j*3+2-.8+.07) -- (\i*3+1.5-.8+.07,\j*3+3-.8-.07);

                    \fill [black,opacity=1] (\i*3+1.5,\j*3+1) circle (5pt);
                      \fill [black,opacity=1] (\i*3+1.5,\j*3+3) circle (5pt);
                   \path[fill=white,draw=black, thick] (\i*3+.5,\j*3+2) circle[radius=4.5pt];
                   \path[fill=white,draw=black, thick] (\i*3+.5+2,\j*3+2) circle[radius=4.5pt];

            \fi
  \node[font=\scriptsize] at (1.5,6.5) {$\tt 4687$};
  \node[font=\scriptsize]  at (1.5,3.5) {$\tt 2485$};
  \node[font=\scriptsize]  at (1.5,.5) {$\tt 2461$};

  \node[font=\scriptsize]  at (4.5,6.5) {$\tt 4681$};
  \node[font=\scriptsize]  at (4.5,3.5) {$\tt 2685$};
  \node[font=\scriptsize]  at (4.5,.5) {$\tt 2463$};

  \node[font=\scriptsize]  at (7.5,6.5) {$\tt 2461$};
  \node[font=\scriptsize]  at (7.5,3.5) {$\tt 2687$};
  \node[font=\scriptsize]  at (7.5,0.5) {$\tt 2483$};

  \node[font=\scriptsize]  at (10.5,6.5) {$\tt 2463$};
  \node[font=\scriptsize]  at (10.5,3.5) {$\tt 4687$};
  \node[font=\scriptsize]  at (10.5,0.5) {$\tt 2485$};

  \node[font=\scriptsize]  at (13.5,6.5) {$\tt 2483$};
  \node[font=\scriptsize]  at (13.5,3.5) {$\tt 4681$};
  \node[font=\scriptsize]  at (13.5,0.5) {$\tt 2685$};

  \node[font=\scriptsize]  at (16.5,6.5) {$\tt 2485$};
  \node[font=\scriptsize]  at (16.5,3.5) {$\tt 2461$};
  \node[font=\scriptsize]  at (16.5,.5) {$\tt 2687$};

  \node[font=\scriptsize]  at (19.5,6.5) {$\tt 2685$};
  \node[font=\scriptsize]  at (19.5,3.5) {$\tt 2463$};
  \node[font=\scriptsize]  at (19.5,.5) {$\tt 4687$};

  \node[font=\scriptsize]  at (22.5,6.5) {$\tt 2687$};
  \node[font=\scriptsize] at (22.5,3.5) {$\tt 2483$};
  \node[font=\scriptsize]  at (22.5,.5) {$\tt 4681$};

  \node[font=\scriptsize]  at (0,4.6) {$\tt 2487$};
  \node[font=\scriptsize]  at (0,1.6) {$\tt 2481$};
  
  \node[font=\scriptsize]  at (3,4.6) {$\tt 4685$};
  \node[font=\scriptsize]  at (3,1.6) {$\tt 2465$};
  
  \node[font=\scriptsize]  at (6,4.6) {$\tt 2681$};
  \node[font=\scriptsize]  at (6,1.6) {$\tt 2683$};

  \node[font=\scriptsize]  at (9,4.6) {$\tt 2467$};
  \node[font=\scriptsize]  at (9,1.6) {$\tt 2487$};

  \node[font=\scriptsize]  at (12,4.6) {$\tt 4683$};
  \node[font=\scriptsize]  at (12,1.6) {$\tt 4685$};

  \node[font=\scriptsize]  at (15,4.6) {$\tt 2481$};
  \node[font=\scriptsize]  at (15,1.6) {$\tt 2681$};

  \node[font=\scriptsize]  at (18,4.6) {$\tt 2465$};
  \node[font=\scriptsize]  at (18,1.6) {$\tt 2467$};

  \node[font=\scriptsize]  at (21,4.6) {$\tt 2683$};
  \node[font=\scriptsize]  at (21,1.6) {$\tt 4683$};

  \node[font=\scriptsize]  at (24,4.6) {$\tt 2487$};
  \node[font=\scriptsize]  at (24,1.6) {$\tt 2481$};
  \node[purple,fill=white] (r) at (0,3.5) {$\tt 248$};
  \node[purple] (r) at (3,0.5) {$\tt 246$};
  \node[purple] (r) at (3,6.5) {$\tt 468$};

  \node[purple] (r) at (6,3.5) {$\tt 268$};
  \node[purple] (r) at (9,0.5) {$\tt 248$};
  \node[purple] (r) at (9,6.5) {$\tt 246$};

  \node[purple] (r) at (12,3.5) {$\tt 468$};
  \node[purple] (r) at (15,0.5) {$\tt 268$};
  \node[purple] (r) at (15,6.5) {$\tt 248$};

  \node[purple] (r) at (18,3.5) {$\tt 246$};
  \node[purple] (r) at (21,0.5) {$\tt 468$};
  \node[purple] (r) at (21,6.5) {$\tt 268$};

  \node[purple,fill=white] (r) at (24,3.5) {$\tt 248$};
  
   \node[blue] (r) at (0,6.5) {$\tt 7$};
   \node[blue] (r) at (3,3.5) {$\tt 5$};
   \node[blue] (r) at (0,0.5) {$\tt 1$};

  \node[blue] (r) at (6,6.5) {$\tt 1$};
   \node[blue] (r) at (9,3.5) {$\tt 7$};
   \node[blue] (r) at (6,0.5) {$\tt 3$};

  \node[blue] (r) at (12,6.5) {$\tt 3$};
   \node[blue] (r) at (15,3.5) {$\tt 1$};
   \node[blue] (r) at (12,0.5) {$\tt 5$};

  \node[blue] (r) at (18,6.5) {$\tt 5$};
   \node[blue] (r) at (21,3.5) {$\tt 3$};
   \node[blue] (r) at (18,0.5) {$\tt 7$};

  \node[blue] (r) at (24,6.5) {$\tt 7$};
   \node[blue] (r) at (24,0.5) {$\tt 1$};

         }
       }     

  \end{scope}
        
                    \coordinate (a11) at (-.7,2.5);
            \coordinate (a12) at (-.7,5.5);
           \coordinate (a21) at (24.7,2.5);
            \coordinate (a22) at (24.7,5.5);
            \coordinate (a31) at (0,2.5-.2);
            \coordinate (a32) at (24.2,2.5-.2);
            \coordinate (a41) at (6.1,5.7);
            \coordinate (a42) at (24.2,5.7);
            \coordinate (a51) at (0,5.7);
            \coordinate (a52) at (5.9,5.7);
  \draw[->, red, line width=1pt] (a11) -- (a12);
            \draw[->, red, line width=1pt] (a21) -- (a22);
            \draw[-, red, line width=1pt] ($(a11)+(.15,0)$) -- ($(a11)+(-.15,0)$);
            \draw[-, red, line width=1pt] ($(a21)+(.15,0)$) -- ($(a21)+(-.15,0)$);

            \draw[->, green!60!black, line width=1pt] (a31) -- (a32);
            \draw[->, green!60!black, line width=1pt] (a41) -- (a42);
            \draw[->, green!60!black, line width=1pt] (a51) -- (a52);
            \draw[-, green!60!black, line width=1pt] ($(a31)+(0,.15)$) -- ($(a31)+(0,-.15)$);
            \draw[-, green!60!black, line width=1pt] ($(a41)+(0,.15)$) -- ($(a41)+(0,-.15)$);

               \end{tikzpicture}
    \end{minipage}
    
    \caption{Proof for the second version of the 16 point theorem. \label{fig:proof2}}
\end{figure}

\medskip 
Next we give the corresponding proving schemes for the incidence theorem shown in 
Figure~\ref{fig:tori} on the right, the second way to express the 16 point theorem.
Also here in \cite{BRG2} a bi-quadratic final polynomial was given. Instead of listing the equations, we  immediately start with
the graphical representation for the  bi-quadratic proof. It is shown in Figure~\ref{fig:proof2}. To make it easier to see the structure we show a region that corresponds to {\it two} full covers of the underlying torus. The rectangular dashed region corresponds to one full cover of the torus. To obtain the torus in that region the left and right side have to be identified and to bottom and top side are identified with a shift of 2 units. Again the shaded squares represent binomial equations and the black and white points meeting at the corners represent the same bracket (once on the left and once on the right of the bi-quadratic equations).

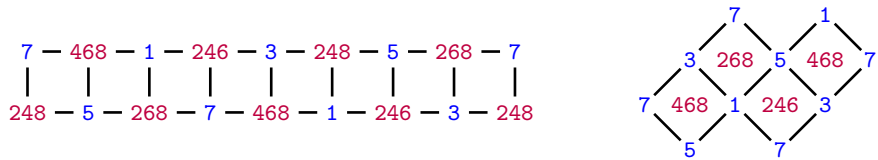
\begin{figure}
    \centering
    \begin{minipage}[c]{.55\linewidth}
        \centering
        \begin{tikzpicture}[>=Stealth,scale=0.27]

           \begin{scope}

  \node[purple] (248b) at (0,3.5) {$\tt 248$};
  \node[purple] (468t) at (3,6.5) {$\tt 468$};

  \node[purple] (268b) at (6,3.5) {$\tt 268$};
  \node[purple] (246t) at (9,6.5) {$\tt 246$};

  \node[purple] (468b) at (12,3.5) {$\tt 468$};
  \node[purple] (248t) at (15,6.5) {$\tt 248$};

  \node[purple] (246b) at (18,3.5) {$\tt 246$};
  \node[purple] (268t) at (21,6.5) {$\tt 268$};

  \node[purple] (248b2) at (24,3.5) {$\tt 248$};
  
   \node[blue] (7t) at (0,6.5) {$\tt 7$};
   \node[blue] (5b) at (3,3.5) {$\tt 5$};

  \node[blue] (1t) at (6,6.5) {$\tt 1$};
   \node[blue] (7b) at (9,3.5) {$\tt 7$};

  \node[blue] (3t) at (12,6.5) {$\tt 3$};
   \node[blue] (1b) at (15,3.5) {$\tt 1$};

  \node[blue] (5t) at (18,6.5) {$\tt 5$};
   \node[blue] (3b) at (21,3.5) {$\tt 3$};

  \node[blue] (7t2) at (24,6.5) {$\tt 7$};

\draw[-, line width=1pt] (7t) -- (468t);
\draw[-, line width=1pt] (1t) -- (468t);
\draw[-, line width=1pt] (1t) -- (246t);
\draw[-, line width=1pt] (3t) -- (246t);
\draw[-, line width=1pt] (3t) -- (248t);
\draw[-, line width=1pt] (5t) -- (248t);
\draw[-, line width=1pt] (5t) -- (268t);
\draw[-, line width=1pt] (7t2) -- (268t);

\draw[-, line width=1pt] (5b) -- (248b);
\draw[-, line width=1pt] (5b) -- (268b);
\draw[-, line width=1pt] (7b) -- (268b);
\draw[-, line width=1pt] (7b) -- (468b);
\draw[-, line width=1pt] (1b) -- (468b);
\draw[-, line width=1pt] (1b) -- (246b);
\draw[-, line width=1pt] (3b) -- (246b);
\draw[-, line width=1pt] (3b) -- (248b2);

\draw[-, line width=1pt] (7t) -- (248b);
\draw[-, line width=1pt] (1t) -- (268b);
\draw[-, line width=1pt] (3t) -- (468b);
\draw[-, line width=1pt] (5t) -- (246b);
\draw[-, line width=1pt] (7t2) -- (248b2);
\draw[-, line width=1pt] (5b) -- (468t);
\draw[-, line width=1pt] (7b) -- (246t);
\draw[-, line width=1pt] (1b) -- (248t);
\draw[-, line width=1pt] (3b) -- (268t);

  \end{scope}
              
               \end{tikzpicture}
    \end{minipage}
    \hfill
      \begin{minipage}[c]{.35\linewidth}
        \centering
        \begin{tikzpicture}[>=Stealth,scale=0.2]

           \begin{scope}

  \node[purple] (468t) at (3,6.5) {$\tt 468$};

  \node[purple] (268b) at (6,9.5) {$\tt 268$};
  \node[purple] (246t) at (9,6.5) {$\tt 246$};

  \node[purple] (468b) at (12,9.5) {$\tt 468$};

   \node[blue] (12) at (0,6.5) {$\tt 7$};
   \node[blue] (21) at (3,3.5) {$\tt 5$};
   \node[blue] (23) at (3,9.5) {$\tt 3$};

  \node[blue] (32) at (6,6.5) {$\tt 1$};
    \node[blue] (34) at (6,12.5) {$\tt 7$};

   \node[blue] (41) at (9,3.5) {$\tt 7$};
   \node[blue] (43) at (9,9.5) {$\tt 5$};

  \node[blue] (52) at (12,6.5) {$\tt 3$};
  \node[blue] (54) at (12,12.5) {$\tt 1$};

   \node[blue] (63) at (15,9.5) {$\tt 7$};

\draw[-, line width=1pt,shorten <=-3pt, shorten >=-3pt] (12) -- (21);
\draw[-, line width=1pt,shorten <=-3pt, shorten >=-3pt] (12) -- (23);
\draw[-, line width=1pt,shorten <=-3pt, shorten >=-3pt] (32) -- (21);
\draw[-, line width=1pt,shorten <=-3pt, shorten >=-3pt] (32) -- (23);
\draw[-, line width=1pt,shorten <=-3pt, shorten >=-3pt] (34) -- (23);
\draw[-, line width=1pt,shorten <=-3pt, shorten >=-3pt] (32) -- (41);
\draw[-, line width=1pt,shorten <=-3pt, shorten >=-3pt] (32) -- (43);
\draw[-, line width=1pt,shorten <=-3pt, shorten >=-3pt] (34) -- (43);

\draw[-, line width=1pt,shorten <=-3pt, shorten >=-3pt] (41) -- (52);
\draw[-, line width=1pt,shorten <=-3pt, shorten >=-3pt] (52) -- (43);
\draw[-, line width=1pt,shorten <=-3pt, shorten >=-3pt] (54) -- (43);
\draw[-, line width=1pt,shorten <=-3pt, shorten >=-3pt] (52) -- (63);
\draw[-, line width=1pt,shorten <=-3pt, shorten >=-3pt] (54) -- (63);

  \end{scope}
              
               \end{tikzpicture}
    \end{minipage}

    \caption{A quad-proof and Menelaus proof for the  second version of the 16 point theorem. \label{fig:quad2}}
\end{figure}

In a completely analogous way to the previous example we can extract CM-proofs and quad-proofs from this scheme. The corresponding schemes for the quad-proof and Menelaus proof is shown in Figure~\ref{fig:quad2}.

\begin{figure}[h]
    \centering
    \begin{tikzpicture}
        \draw (0, 0) node[inner sep=0] {\includegraphics[width=11cm]{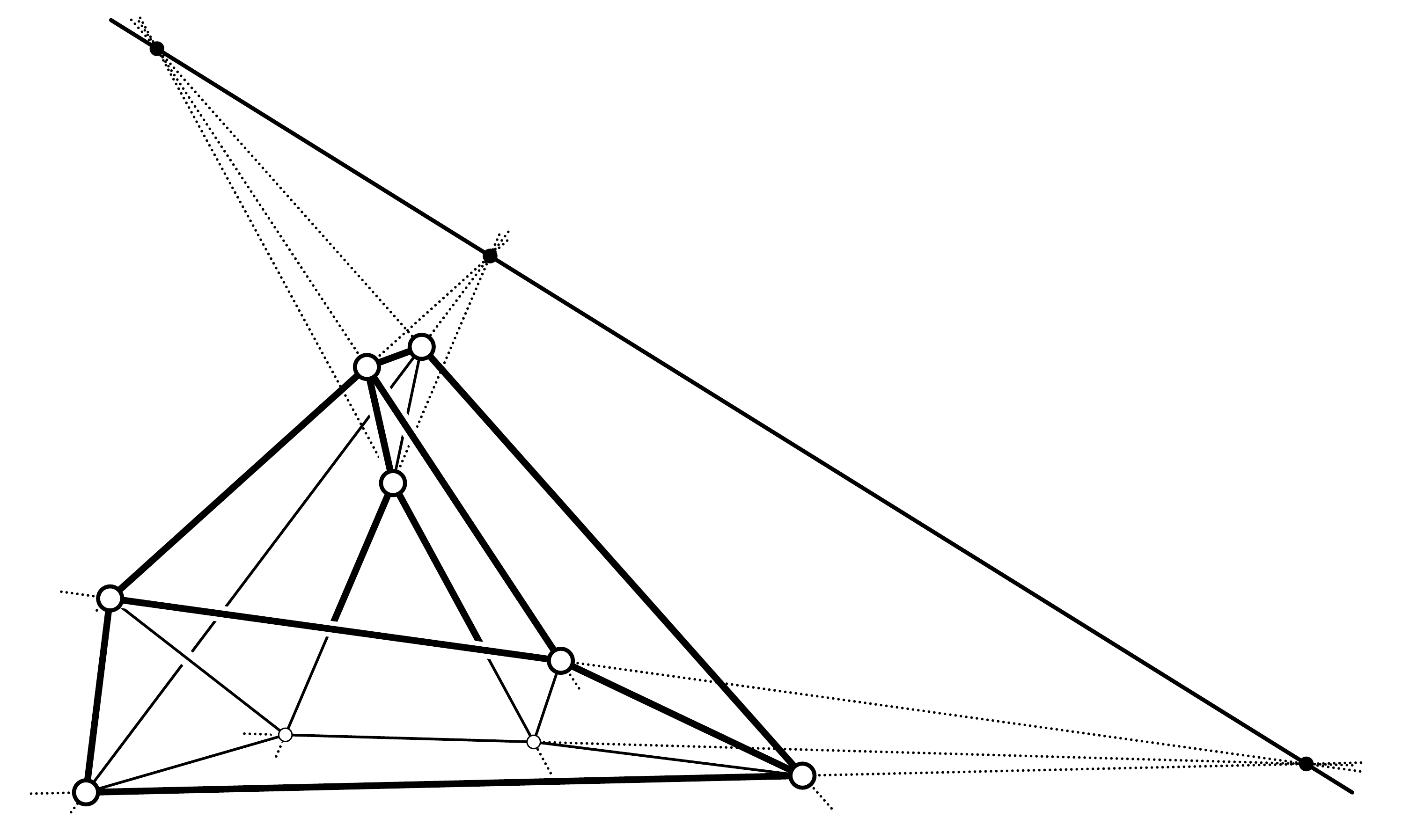}};
           \node[blue] at (-4.7,-1.1) {$\tt 1$};
           \node[blue] at (-2.6,.7) {$\tt 2$};
           \node[blue] at (-5,-2.7) {$\tt 4$};
           \node[blue] at (-3.3,-2.2) {$\tt 7$};
           \node[blue] at (-2.7,-.5) {$\tt 8$};
           \node[blue] at (-2.23,.95) {$\tt 5$};
           \node[blue] at (-1,-1.6) {$\tt 3$};
           \node[blue] at (-1.7,-2.3) {$\tt 9$};
           \node[blue] at (.7,-3.) {$\tt 6$};

 \node[blue] (20) at (2,0) {$\tt 1$};
 \node[blue] (21) at (2,1) {$\tt 4$};
 \node[blue] (22) at (2,2) {$\tt 7$};
 \node[blue] (23) at (2,3) {$\tt 1$};
 \node[blue] (30) at (3,0) {$\tt 2$};
 \node[blue] (31) at (3,1) {$\tt 5$};
 \node[blue] (32) at (3,2) {$\tt 8$};
 \node[blue] (33) at (3,3) {$\tt 2$};
 \node[blue] (40) at (4,0) {$\tt 3$};
 \node[blue] (41) at (4,1) {$\tt 6$};
 \node[blue] (42) at (4,2) {$\tt 9$};
 \node[blue] (43) at (4,3) {$\tt 3$};
 \node[blue] (50) at (5,0) {$\tt 1$};
 \node[blue] (51) at (5,1) {$\tt 4$};
 \node[blue] (52) at (5,2) {$\tt 7$};
 \node[blue] (53) at (5,3) {$\tt 1$};
        
\draw[-, line width=1pt,shorten <=-1pt, shorten >=-1pt] (20) -- (21);
\draw[-, line width=1pt,shorten <=-1pt, shorten >=-1pt] (21) -- (22);
\draw[-, line width=1pt,shorten <=-1pt, shorten >=-1pt] (22) -- (23);
\draw[-, line width=1pt,shorten <=-1pt, shorten >=-1pt] (30) -- (31);
\draw[-, line width=1pt,shorten <=-1pt, shorten >=-1pt] (31) -- (32);
\draw[-, line width=1pt,shorten <=-1pt, shorten >=-1pt] (32) -- (33);
\draw[-, line width=1pt,shorten <=-1pt, shorten >=-1pt] (40) -- (41);
\draw[-, line width=1pt,shorten <=-1pt, shorten >=-1pt] (41) -- (42);
\draw[-, line width=1pt,shorten <=-1pt, shorten >=-1pt] (42) -- (43);
\draw[-, line width=1pt,shorten <=-1pt, shorten >=-1pt] (50) -- (51);
\draw[-, line width=1pt,shorten <=-1pt, shorten >=-1pt] (51) -- (52);
\draw[-, line width=1pt,shorten <=-1pt, shorten >=-1pt] (52) -- (53);
\draw[-, line width=1pt,shorten <=-1pt, shorten >=-1pt] (20) -- (30);
\draw[-, line width=1pt,shorten <=-1pt, shorten >=-1pt] (30) -- (40);
\draw[-, line width=1pt,shorten <=-1pt, shorten >=-1pt] (40) -- (50);
\draw[-, line width=1pt,shorten <=-1pt, shorten >=-1pt] (21) -- (31);
\draw[-, line width=1pt,shorten <=-1pt, shorten >=-1pt] (31) -- (41);
\draw[-, line width=1pt,shorten <=-1pt, shorten >=-1pt] (41) -- (51);
\draw[-, line width=1pt,shorten <=-1pt, shorten >=-1pt] (22) -- (32);
\draw[-, line width=1pt,shorten <=-1pt, shorten >=-1pt] (32) -- (42);
\draw[-, line width=1pt,shorten <=-1pt, shorten >=-1pt] (42) -- (52);
\draw[-, line width=1pt,shorten <=-1pt, shorten >=-1pt] (23) -- (33);
\draw[-, line width=1pt,shorten <=-1pt, shorten >=-1pt] (33) -- (43);
\draw[-, line width=1pt,shorten <=-1pt, shorten >=-1pt] (43) -- (53);

    \end{tikzpicture}
    \caption{Geometry and combinatorics of the Toblerone Theorem.}\label{fig:Toblerone}
\end{figure}
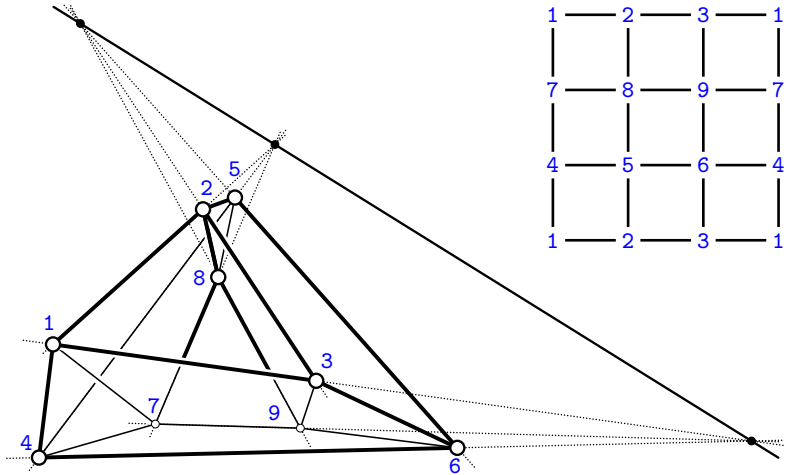
           
\subsection{The $3\times 3$ torus}            
Our final example is sometimes jokingly called the Toblerone Torus, since it can be assembled from three Toblerone boxes out of which one builds a triangular frame.
It is the torus  formed by a grid of $3\times 3$ quadrilaterals. It has 9 vertices. 
If eight of the quadrilaterals are flat the final one is flat automatically. Figure~\ref{fig:Toblerone} shows a projectively correct 2D drawing of this configuration along with its combinatorial structure.

Again in \cite{RG2} a bi-quadratic final polynomial for this incidence theorem was provided. We present this final polynomial here  in an pre-structured form. We separate it into three blocks. Within each block the red and the green parts cancel directly and we  have a chain of three bi-quadratic equations of which only the outer parts --the black brackets in the equations-- are relevant (similar to the chains we observed in Figure~\ref{fig: the big picture}). We show these blocks along with a graphical representation of the cancellation chains that belong to each block. 

\begin{minipage}[c]{.6\linewidth}
\[
\begin{matrix}
[\1\2\7\3]{\color{red}[\1\2\8\4]} = {\color{red}[\1\2\7\4]}[\1\2\8\3] \quad\Leftarrow\quad  h(\7\8\1\2)\\
{\color{ForestGreen}[\1\4\5\8]} {\color{red}[\1\4\2\7]}= {\color{red}[\1\4\2\8]}{\color{ForestGreen}[\1\4\5\7]} \quad\Leftarrow\quad h(\1\2\4\5)\\
{\color{ForestGreen}[\4\5\7\1]}[\4\5\8\6] = [\4\5\7\6]{\color{ForestGreen}[\4\5\8\1]} \quad\Leftarrow\quad h(\4\5\7\8)\\[2mm]

[\2\3\8\1] {\color{red}[\2\3\9\5]}= {\color{red}[\2\3\8\5]} [\2\3\9\1]\quad\Leftarrow\quad h(\8\9\2\3)\\
{\color{ForestGreen}[\2\5\6\9]}{\color{red}[\2\5\3\8]} = {\color{red}[\2\5\3\9]}{\color{ForestGreen}[\2\5\6\8]} \quad\Leftarrow\quad  h(\2\3\5\6)\\
{\color{ForestGreen}[\5\6\8\2]}[\5\6\9\4] = [\5\6\8\4] {\color{ForestGreen}[\5\6\9\2]}\quad\Leftarrow\quad  h(\5\6\8\9)\\[2mm]

[\3\1\9\2]{\color{red}[\3\1\7\6]} ={\color{red}[\3\1\9\6]}[\3\1\7\2] \quad\Leftarrow\quad h(\9\7\4\1)\\
{\color{ForestGreen}[\3\6\4\7]}{\color{red}[\3\6\1\9]} = {\color{red}[\3\6\1\7}]{\color{ForestGreen}[\3\6\4\9]} \quad\Leftarrow\quad h(\3\1\6\4)\\
{\color{ForestGreen}[\6\4\9\3]}[\6\4\7\5] = [\6\4\9\5]{\color{ForestGreen}[\6\4\7\3]}\quad\Rightarrow\quad h(\6\4\9\7)\\
\end{matrix}
\]
\end{minipage}
\hfill
\raisebox{-.25cm}{
\begin{minipage}[c]{.28\linewidth}
    \begin{tikzpicture}[>=Stealth,scale=0.6]
    
     \foreach \i in {0,1.5,3} {
      \foreach \j in {0,2.5,2.5*2} {
       \draw[-, line width=.7pt] (0+\i,1+\j) -- (1+\i,1+\j);
      \draw[-, line width=.7pt] (0+\i,1+\j) -- (0+\i,0+\j);
      \draw[-, line width=.7pt] (1+\i,0+\j) -- (1+\i,1+\j);
      \draw[-, line width=.7pt] (1+\i,0+\j) -- (0+\i,0+\j);
      \fill [black,opacity=1] (0+\i,1+\j) circle (4pt);
      \fill [black,opacity=1] (1+\i,0+\j) circle (4pt);
      \path[fill=white,draw=black, thick] (0+\i,0+\j) circle[radius=3.5pt];
      \path[fill=white,draw=black, thick] (1+\i,1+\j) circle[radius=3.5pt];

    }
    }
    
      \node[font=\scriptsize]  at (-.25,-0.4) {$\tt 1237$};
      \node[font=\scriptsize]  at (-.25,1.4) {$\tt 1239$};
      \node[font=\scriptsize]  at (1.25,-0.4) {$\tt 1367$};
      \node[font=\scriptsize]  at (1.25,1.4) {$\tt 1369$};
      \node[font=\scriptsize]  at (2.75,-0.4) {$\tt 3467$};
      \node[font=\scriptsize]  at (2.75,1.4) {$\tt 3469$};
      \node[font=\scriptsize]  at (4.25,-0.4) {$\tt 4567$};
      \node[font=\scriptsize]  at (4.25,1.4) {$\tt 4569$};
          
      \node[font=\scriptsize]  at (-.25,-0.4+2.5) {$\tt 1239$};
      \node[font=\scriptsize]  at (-.25,1.4+2.5) {$\tt 1238$};
      \node[font=\scriptsize]  at (1.25,-0.4+2.5) {$\tt 2359$};
      \node[font=\scriptsize]  at (1.25,1.4+2.5) {$\tt 2358$};
      \node[font=\scriptsize]  at (2.75,-0.4+2.5) {$\tt 2569$};
      \node[font=\scriptsize]  at (2.75,1.4+2.5) {$\tt 2568$};
      \node[font=\scriptsize]  at (4.25,-0.4+2.5) {$\tt 4569$};
      \node[font=\scriptsize]  at (4.25,1.4+2.5) {$\tt 4568$};
               
      \node[font=\scriptsize]  at (-.25,-0.4+5) {$\tt 1238$};
      \node[font=\scriptsize]  at (-.25,1.4+5) {$\tt 1237$};
      \node[font=\scriptsize]  at (1.25,-0.4+5) {$\tt 1248$};
      \node[font=\scriptsize]  at (1.25,1.4+5) {$\tt 1247$};
      \node[font=\scriptsize]  at (2.75,-0.4+5) {$\tt 1458$};
      \node[font=\scriptsize]  at (2.75,1.4+5) {$\tt 1457$};
      \node[font=\scriptsize]  at (4.25,-0.4+5) {$\tt 4568$};
      \node[font=\scriptsize]  at (4.25,1.4+5) {$\tt 4567$};
          
    \end{tikzpicture}
\end{minipage}
}
\medskip

\medskip

Since only the terminal edges of these chains can form  ratios that are needed for either CM-proofs or quad-proofs we can immediately collapse these chains and focus on the resulting equations. Doing so we are left with the following three equations (digits have been sorted while taking care of the signs of the determinants).

\begin{equation}
\begin{matrix}
[\1\2\3{\color{blue}\7}][\4\5\6{\color{blue}\8}] = [\1\2\3{\color{blue}\8}][\4\5\6{\color{blue}\7}]\\
[\1\2\3{\color{blue}\8}][\4\5\6{\color{blue}\9}] = [\1\2\3{\color{blue}\9}][\4\5\6{\color{blue}\8}]\\
[\1\2\3{\color{blue}\9}][\4\5\6{\color{blue}\7}] = [\1\2\3{\color{blue}\7}][\4\5\6{\color{blue}\9}]\\
\end{matrix} \label{eq:threeeqs}
\end{equation}

The amazingly difficult part of translating these equations into a CM-proof or a quad-proof comes from the fact that there es hardly any structure left that relates these equations to the geometry of the theorem. Let us first interpret the structure as a quad-proof.

\begin{figure}[h]
    \centering
    \begin{tikzpicture}
        \draw (0, 0) node[inner sep=0] {\includegraphics[width=6.5cm]{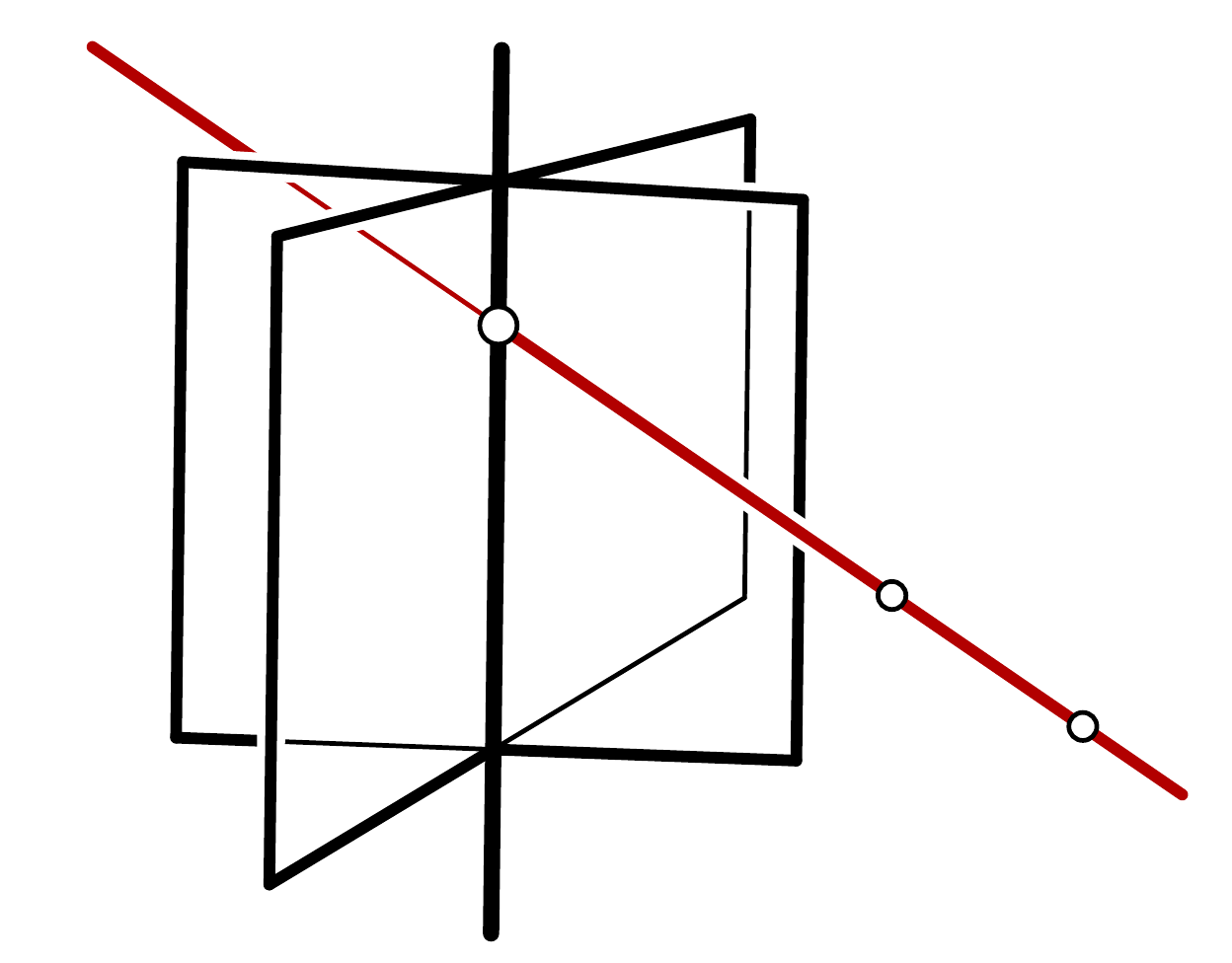}};
                \draw (6, 0) node[inner sep=0] {\includegraphics[width=5cm]{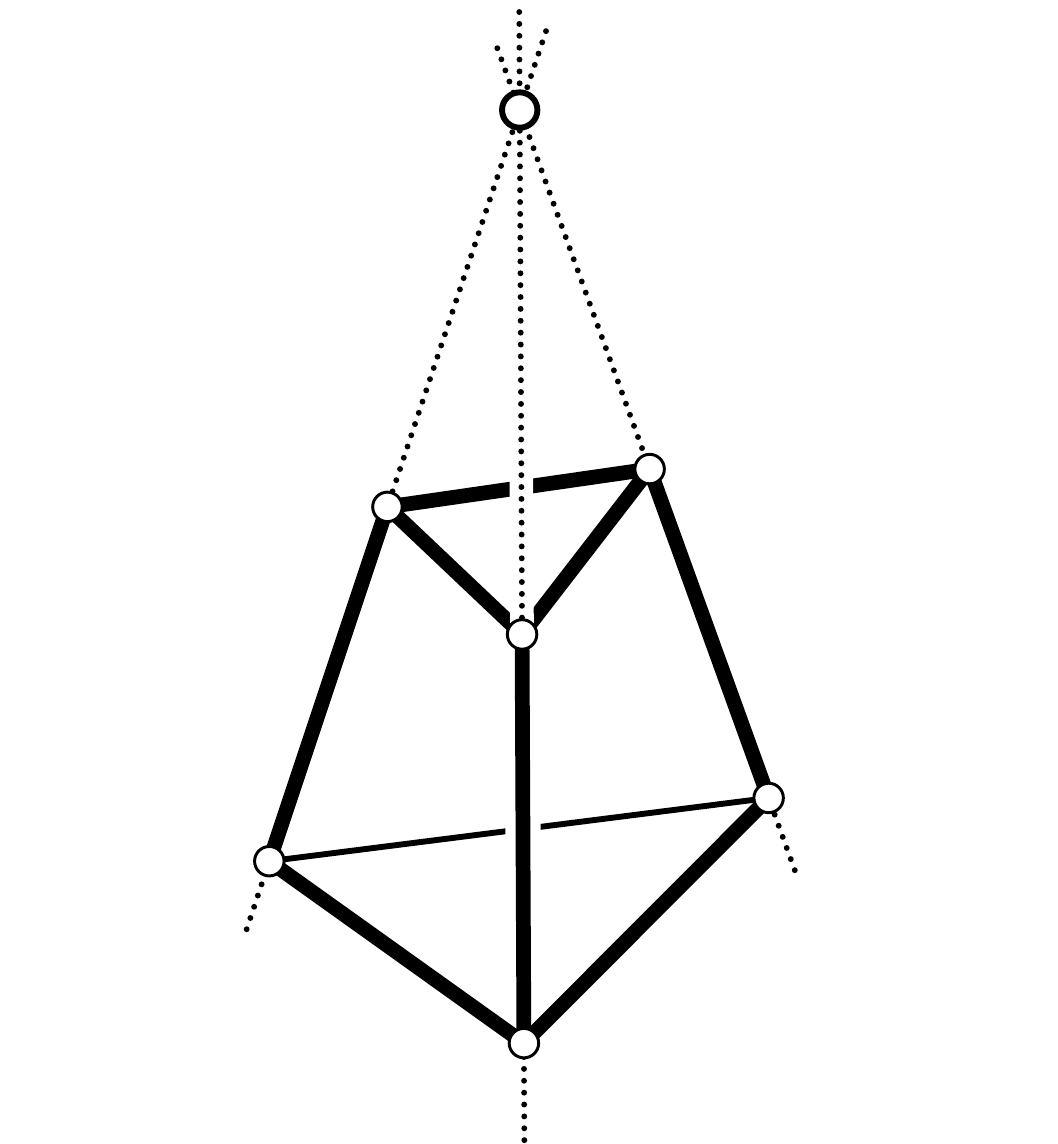}};
      \node  at (-1.25,-1.6) {$\tt 1$};
      \node  at (-1.5,-1.55) {$\tt 2$};
      \node  at (-1.65,-1.8) {$\tt 3$};
            \node  at (.2,-1.2) {$\tt 4$};
      \node  at (.5,-1.) {$\tt 5$};
      \node  at (.7,-1.3) {$\tt 6$};
      \node [blue] at (2.7,-1.2) {$\tt 8$};
      \node [blue] at (1.7,-0.5) {$\tt 7$};

      \node  at (4.5,-1.4) {$\tt 1$};
      \node  at (5.1,.4) {$\tt 2$};
      \node  at (5.8,-2.4) {$\tt 4$};
      \node  at (5.75,-.4) {$\tt 5$};
      \node  at (7.35,-1.0) {$\tt 7$};
      \node  at (6.8,.6) {$\tt 8$};

    \end{tikzpicture}
    \caption{Geometric relations to interpret equation (\ref{eq:bundle}). }\label{fig:bundle}
\end{figure}
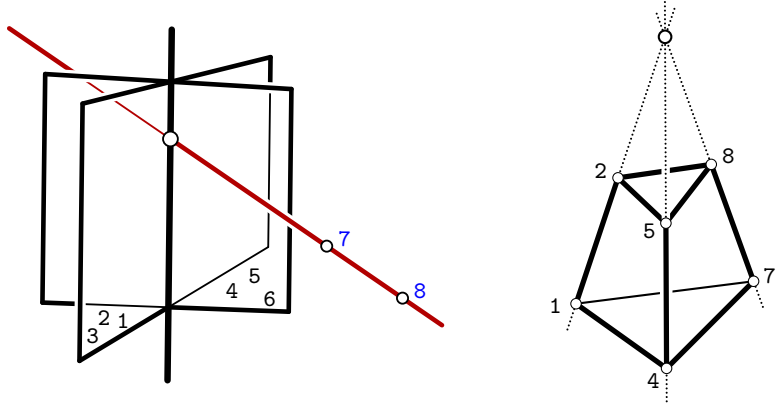

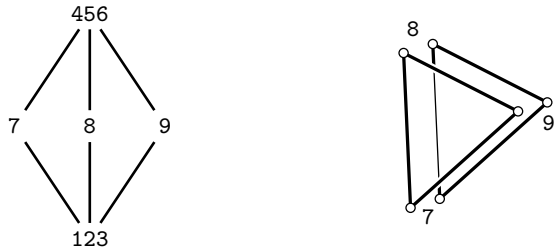
\begin{figure}[h]
    \centering
    \begin{tikzpicture}

      \node (123) at (0,0.5) {$\tt 123$};
      \node (456) at (0,3.5) {$\tt 456$};
      \node (7) at (-1,2) {$\tt 7$};
      \node (8) at (0,2) {$\tt 8$};
      \node (9) at (1,2) {$\tt 9$};
        \draw[-, line width=1pt] (123) -- (7);
        \draw[-, line width=1pt] (123) -- (8);
        \draw[-, line width=1pt] (123) -- (9);
     \draw[-, line width=1pt] (456) -- (7);
        \draw[-, line width=1pt] (456) -- (8);
        \draw[-, line width=1pt] (456) -- (9);

    \node  at (4.5,0.8) {$\tt 7$};
      \node  at (4.4-.1,1.1+2.2) {$\tt 8$};
      \node  at (4.8+1.3,0.7+1.35) {$\tt 9$};
     
          \draw (5, 2) node[inner sep=0] {\includegraphics[width=3cm]{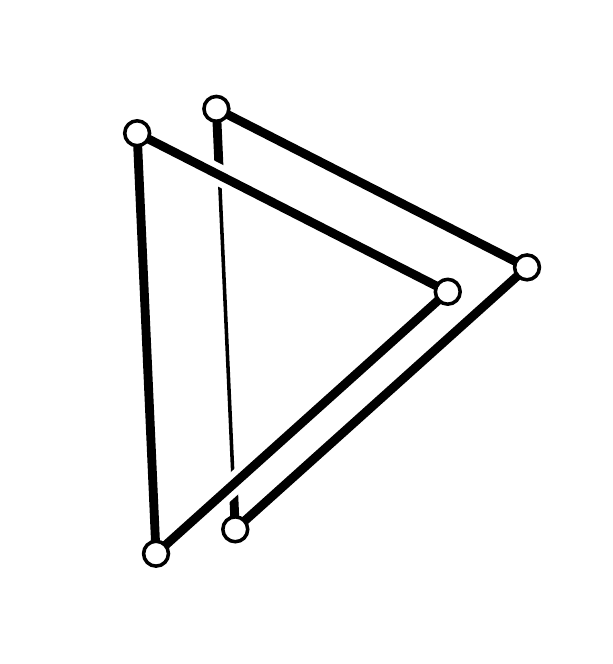}};

    \end{tikzpicture}
    \caption{Quad proof and CM-proof of for the Toblerone torus.}\label{fig:tobleronemanifold}
\end{figure}

\noindent
We first consider the equation 
\begin{equation}
[\1\2\3{\color{blue}\7}][\4\5\6{\color{blue}\8}] = [\1\2\3{\color{blue}\8}][\4\5\6{\color{blue}\7}]\label{eq:bundle}
\end{equation}
and its geometric relevance. It encodes a geometric relation  between the planes spanned by $\tt 123$, $\tt 456$ and the line spanned by $\tt 78$. It is exactly the type of equation given by  the equation (\ref{fomin eq}), the basic building block of a quad-proof. It holds if the join of  $\tt 7$ and $\tt 8$ hits the intersection line of the planes  $\tt 123$ and $\tt 456$. The situation is depicted in Figure~\ref{fig:bundle} on the left.
Thus our three equations can be directly interpreted as a quad-proof, representing a very simple configuration. The underlying manifold consists of 3 quadrangles that form a spherical structure as indicated in Figure~\ref{fig:tobleronemanifold}~(left). In that picture also the exterior face has to be considered as a quadrangle.

Why does this manifold structure proof our theorem? For this again consider the equation
 $[\1\2\3{\color{blue}\7}][\4\5\6{\color{blue}\8}] = [\1\2\3{\color{blue}\8}][\4\5\6{\color{blue}\7}]$. 
 and consider the three coplanarities $\tt 1245$, $\tt 1278$ and $\tt 4578$. They are part of our geometric configuration of the Toblerone torus. Since any three planes meet in a point (perhaps infinite) this implies that the lines $\overline{\tt 12}$, $\overline{\tt 45}$, $\overline{\tt  78}$ meet in a point. This in turn is exactly encoded by equation (\ref{eq:bundle}), where $\tt 3$ and $\tt 6$ play the role of auxiliary points in suitably general position. Conversely, if (\ref{eq:bundle}) holds and we know that 
 $\tt 1245$, $\tt 1278$ are coplanar we can conclude that  $\tt 4578$ is coplanar as well. These two statements immediately translate into the fact that the quad-manifold associated to (\ref{eq:threeeqs}) proofs the Toblerone torus theorem: If eight coplanarities are satisfied
 then two of the equations hold. This implies that the third equation holds this together with the two remaining coplanarities implies the conclusion of the theorem.

\medskip

Finally, it remains to interpret the equations in terms of a CM-proof. For this we again can consider the equations in (\ref{eq:threeeqs}) as glue between triangle edges. The manifold is as trivial as it could be. It just consists of two copies of the triangle $\tt 789$  as Menelaus triangles and glued along their edges. One instance is assumed to be cut by the plane $\tt 123$. The other instance is assumed to be cut by $456$. 
The interplay of the  geometric situation and the CM-proof can be understood by  reconsidering Figure~\ref{fig:Toblerone}. The three triangles  $\tt 123$, $\tt 456$ and  $\tt 789$ are boundary a cycle of the torus.
Cutting $\tt 789$ with a plane $\tt 123$ produces additional points on the lines that support the edges (one can see them in Figure~\ref{fig:Toblerone}). They are the three points of the perspective line indicated in this figure.  The almost trivial CM-proof manifold states that if one chooses plane $\tt 789$ such that it passes trough the intersection line of $\tt 123$ and $\tt 456$ then the following holds. If we have 
$\tt 78 \wedge 123=78 \wedge 456$ and 
$\tt 89 \wedge 123=89 \wedge 456$ the we automatically have
 $\tt 79 \wedge 123=79 \wedge 456$. This is exactly what we need for the incidence theorem to hold.

It is remarkable that CM-proofs that produce a void statement in a 2D interpretation produce an interesting incidence theorem when interpreted in a 3D framework. We do not claim that thinking of the proofs in that case in terms of quad-proofs or CM-proofs is conceptually simpler than using a bi-quadratic framework. However, it is nice to see that even these highly degenerate proving manifolds smoothly fit in our general interplay of the three different proving techniques.

\end{document}